\documentclass[12 pt]{article}
\usepackage[utf8]{inputenc}
\usepackage[T1]{fontenc}
\usepackage[top=2.5cm, bottom=2.5cm, left=3cm, right=3cm]{geometry}
\usepackage{amsmath}
\usepackage{amssymb}
\usepackage{amsthm}
\usepackage{mathtools}
\usepackage{mathrsfs}
\usepackage{braket}
\usepackage{enumitem}
\usepackage{comment}
\usepackage{xcolor}
\usepackage[pdftex,unicode=true,bookmarks,pdfpagelabels,pagebackref]{hyperref}
\usepackage{tabu}
\usepackage{stmaryrd}
\usepackage{authblk}

\newcommand\mb{\mathbf}
\newcommand\sco{\,;\,}
\newcommand\scr{\mathscr}
\newcommand\mr{\mathrm}

\DeclarePairedDelimiterX{\norm}[1]{\lVert}{\rVert}{#1}
\DeclarePairedDelimiterX{\nm}[1]{\lvert}{\rvert}{#1}
\DeclarePairedDelimiterX{\intbk}[1]{\llbracket}{\rrbracket}{#1}

\newtheoremstyle{mystyle}
{10pt}
{4pt}
{\itshape}
{}
{\bfseries}
{.}
{.5em}
{}
\theoremstyle{mystyle}

\newtheorem{The}{Theorem}[section]
\newtheorem{Lem}[The]{Lemma}
\newtheorem{Pro}[The]{Proposition}
\newtheorem{Cor}[The]{Corollary}

\newtheorem{Ass}{Assumption}
\renewcommand\theAss{\textbf{\Alph{Ass}}}

\newtheoremstyle{mystylebis}
{10pt}
{10pt}
{}
{}
{\bfseries}
{.}
{.5em}
{}
\theoremstyle{mystylebis}

\newtheorem{Rem}[The]{Remark}
\newtheorem{Exa}[The]{Example}

\newcommand\prf{\noindent\textbf{Proof.\ }}

\begin{document}

\title{Large population asymptotics for a multitype stochastic SIS epidemic model in randomly switched environment}

\author[1,2]{Adrien Prodhomme }

\author[3]{\'Edouard Strickler}

\affil[1]{ Institut Denis Poisson, Université de Tours, France }
\affil[2]{CMAP, Ecole Polytechnique, France}
\affil[3]{ Université de Lorraine, CNRS, Inria, IECL, Nancy, France }


\maketitle

	\begin{abstract}
		We consider an epidemic SIS model described by a multitype birth-and-death process in a randomly switched environment. That is, the infection and cure rates of the process depend on the state of a finite Markov jump process (the environment), whose transitions also depend on the number of infectives. The total size of the population is constant and equal to some $K\in\mb{N}^*$, and the number of infectives vanishes almost surely in finite time. We prove that, as $K \to \infty$, the process composed of the proportions of infectives of each type $X^K$ and the state of the environment $\Xi^K$, converges to a piecewise deterministic Markov process (PDMP) given by a system of randomly switched ODEs. The long term behaviour of this PDMP has been previously investigated by Benaïm and Strickler, and depends only on the sign of the top Lyapunov exponent $\Lambda$ of the linearised PDMP at $0$: if $\Lambda < 0$, the proportion of infectives in each group converges to zero, while if $\Lambda > 0$, the disease becomes endemic. In this paper, we show that the large population asymptotics of $X^K$ also strongly depend on the sign of $\Lambda$: if negative, then from fixed initial proportions of infectives the disease disappears in a time of order at most $\log(K)$, while if positive, the typical extinction time grows at least as a power of $K$. We prove that in the situation where the origin is accessible for the linearised PDMP, the mean extinction time of $X^K$ is logarithmically equivalent to $K^{p^*}$, where $p^* > 0$ is fully characterised. We also investigate the quasi-stationary distribution $\mu^K$ of $(X^K,\Xi^K)$ and show that, when $\Lambda < 0$, weak limit points of $(\mu^K)_{K > 0}$ are supported by the extinction set, while when $\Lambda > 0$, limit points belong to the (non empty) set of stationary distributions of the limiting PDMP which do not give mass to the extinction set. 
\end{abstract}
	
	\section{Introduction}

In the mid-seventies, in an influential paper \cite{LajYorke} Lajmanovich and Yorke developed and investigated a deterministic SIS (Susceptible-Infective-Susceptible) model of infection, describing the evolution of a disease that does not confer immunity (such as gonorrhea) in a population divided into $d$ groups.
The model is given by a differential equation on $[0,1]^d$ having the form
\begin{equation}
\label{eq:LY}
\frac{dx_i}{dt} = (1-x_i) (\sum_{j = 1}^d C_{ij} x_j) - D_i x_i\, ,\quad i \in \intbk{1,d},
\end{equation}
where  $C = (C_{ij})_{(i,j)\in\intbk{1,d}^2}$ is an irreducible matrix with nonnegative entries and $D_i > 0.$
Here $x_i\in[0,1]$ represents the proportion of infected individuals in group $i$. They are assumed to be infective and $C_{ij}$ represents the rate at which group $i$ transmits the infection to group $j$, while $D_i$ is the intrinsic cure rate in group $i$. The irreducibility of $C$ implies that each group indirectly affects the other groups. It is easily seen that the origin is an equilibrium of \eqref{eq:LY}, called in this framework the \emph{disease-free equilibrium} since it corresponds to the absence of the disease in each group. The dynamics of \eqref{eq:LY} is fully described by Lajmanovich and Yorke, and is surprisingly simple : either the disease-free equilibrium is linearly stable, in which case, every solution converges to $0$; or it is unstable, and there exists a unique positive equilibrium $x_*$, globally asymptotically stable on $[0,1]^d \setminus \{0\}$. When it exists, $x_*$ is called the \emph{endemic equilibrium}, and it represents the proportion of infectives in each group when the disease is eventually ingrained in the population. 
\medskip

In order to be more realistic, one should take into account the following features of real life : the populations are finite, and the events of contamination and cure involve randomness. A rigorous derivation of the Lajmanovich-Yorke equation from an individual-based stochastic model was performed by Bena\"{i}m and Hirsch in \cite{BH99} : \eqref{eq:LY} appears as the limit, when the size of the population goes to infinity, of a finite-state discrete-time Markov chain. Actually, Benaim and Hirsch consider an extended version, allowing for a exogenous source of infection, corresponding to the addition of constants $Q_i\geq 0$ in the right handside of \ref{eq:LY}. When $Q_i>0$ for all $i$, the Markov chain is irreducible and the authors show that when the size of the population goes to infinity, its stationary distribution concentrates near the unique stable equilibrium point of the ODE (which in turn converges, in the limit of vanishing $Q_i$'s, to the endemic equilibrium when it exists, or to $0$).

More generally, the links between stochastic population growth models and their mean-field ODE approximation are the topics of numerous papers, especially since the pioneer works of Kurtz (see e.g. \cite{K81}). A very important difference between stochastic and deterministic population models concerns the extinction phenomena : in \eqref{eq:LY}, if there are infectives initially, it remains so at all times, while in a stochastic individual-based model with no external source of infection, the disease will die out with probability one in finite time. Typically, in such a model the extinction time of the disease, starting from a fixed proportion of infectives, will grow with the size $K$ of the population. However, the speed of that growth drastically changes depending on whether the limiting ODE is endemic (i.e. the disease-free equilibrium is unstable) or not. For a classical monotype stochastic SIS model, Kryscio and Lefèvre \cite{KL89}, Andersson and Djehiche \cite{AD98}, and Doering, Sargsyan and Sander \cite{DSS05} proved that, if the disease-free equilibrium is linearly unstable for the ODE, the typical extinction time starting from a fixed proportion of infectives grows exponentially with $K$, while if it is linearly stable, this extinction time is of order $\log(K)$. In particular, in the endemic case, for large populations, it is very unlikely to observe the extinction of the disease. Instead of that, we have better chance to see the population reaching a metastable equilibrium, which can be related to the existence of a quasi-stationary distribution (QSD) for the Markov chain. Once again, the behaviour for large $K$ of the QSD strongly depends  on the nature of the limiting ODE (endemic or not). N{\aa}sell \cite{N96}, \cite{N99} studied this question for the classical monotype SIS model. More recently, in a series of papers, Chazottes, Collet, M\'el\'eard and Martinez study the behaviour of monotype \cite{CCM16} and multitype \cite{CCM17,CCMM20} birth-and-death processes involving a scaling parameter $K$ and such that 0 is linearly unstable for the limiting ODE. They get sharp asymptotics for the extinction rate, the rate of convergence to the QSD, and the total variation distance between the marginal law of the process and a mixture of the QSD and the Dirac mass at 0. They also show that, as $K$ goes to infinity, the QSD of the birth-and-death process approaches a Gaussian law centered on the unique non-zero equilibrium of the ODE (see \cite[Appendix C]{CCMM20}). Results on the extinction rate and the concentration of the QSD around the endemic equilibrium have also been obtained by Schreiber, Huang, Jiang and Wang in \cite{SHJW20} for discrete-time epidemiological stochastic models converging to discrete-time deterministic dynamical systems, based on arguments in the seminal paper of Faure and Schreiber \cite{FS14}.
\medskip

The randomness in the aforementioned Markov chain models is the reflect of \emph{demographic stochasticity}, given by the fact that infections and cures arrive at random time. Another way to add randomness is to take into account \emph{environmental stochasticity}. By this we mean that the population lives in an environment subject to random fluctuations that influences the rates of infections and cures (for more details on demographical and environmental stochasticity, we refer to the nice review of Schreiber \cite{Sch17}). Inspired by a paper by Ait Rami, Bokharaie, Mason and Wirth \cite{Masonetal14}, Bena\"{i}m and Strickler recently considered in \cite{BS19} the case where the matrix $C$ and the vector $D=(D_i)_{i\in \intbk{1,d}}$ in \eqref{eq:LY} are selected in finite families $(C^{\xi})_{\xi \in E}$ and $(D^{\xi})_{\xi \in E}$, and changed at random times. More precisely, they consider a stochastic process $(X(t))_{t \geq 0}=(X_1(t),\ldots,X_d(t))_{t\geq 0}$ with values in $[0,1]^d$ evolving according to 

\begin{equation}
\label{eq:LYR}
\frac{dX_i(t)}{dt} = (1-X_i(t)) (\sum_{j = 1}^d C_{ij}^{\Xi(t)} X_j(t)) - D_i^{\Xi(t)} X_i(t)\, ,\quad i\in \intbk{1, d},
\end{equation}
where $(\Xi(t))_{t \geq 0}$ is a continuous-time Markov chain on the finite set $E$ with some rate matrix $Q=(q(\xi,\xi'))_{(\xi,\xi')\in E^2}$.	

The process $\Xi$ represents an environment subject to abrupt random changes, for instance changes in the weather or in the public health policy (such as lock-down).
The process $\hat{X}=(X,\Xi)$ belongs to the class of \emph{Piecewise Deterministic Markov Processes} (PDMP), a term coined by Davis in one of the first general study on this kind of processes (see \cite{Dav84}). In the last decades, PDMPs generated by switching ODEs have been extensively studied by numerous authors in the context of population dynamics and epidemiology, see e.g , Takueshi et al. \cite{t06},  Du, Dang, and Feng \cite{DNW14},  Bena\"{i}m and Lobry \cite{BL16}, Costa \cite{C16},  Hening and Strickler \cite{HS17} for population dynamics or  Gray et al. \cite{GGMP12} and Li, Liu and Cui \cite{li17} for epidemiology.

After Gray et al. \cite{GGMP12} studied the case of a one dimensional SIS model in two environments, system \eqref{eq:LYR} was fully analysed in \cite{BS19},where it is shown that the behaviour of $X$ near $0$ only depends on the sign of a Lyapunov exponent $\Lambda\in\mb{R}$ of the linearised PDMP at $0$, satisfying
 \begin{equation}
 \label{eq:deflyapexpointro}
 \Lambda := \lim_{t \to \infty} \frac{1}{t} \log \|Y_t \|\quad \text{a.s.},
  \end{equation}
 where $Y_t$ solves $\dot Y_t = A^{\Xi_t} Y_t$, with $A^{\xi} = C^{\xi} - \rm{Diag} (D^{\xi})$ (see Section \ref{sec:PDMP} for a more explicit definition of $\Lambda$). The results in \cite{BS19}, completed by Nguyen and Strickler \cite{NS20} who handled the critical case $\Lambda=0$, can be summed up by the following theorem :
 
 \begin{The}[Bena\"{i}m and Strickler, \cite{BS19}; Nguyen and Strickler, \cite{NS20}]
 \label{th:BS19} Let $\Lambda$ be as defined by Equation \ref{eq:deflyapexpointro}.
 \begin{enumerate} 
 \item\label{item:persistent} If $\Lambda > 0$, there exists a unique invariant probability distribution $\mu^*$ of $\hat{X}$ on $[0,1]^d\times E$ such that $\mu^*(\{0\}\times E)=0$ and, provided that $X_0\neq 0$ a.s., $(X_t, \Xi_t)$ converges in law to $\mu^*$.
 \item\label{item:nonpersistent} If $\Lambda < 0$, then $X_t$ converges almost surely to $0$ exponentially fast. 
 \item If $\Lambda = 0$, then $X_t$ converges in probability to $0$ as $t\rightarrow +\infty$.
 \end{enumerate}
 \end{The}
In the first case ($\Lambda>0$), the PDMP $(X,\Xi)$ is said to be \emph{persistent}, while in the second case ($\Lambda<0$), it is said to be \emph{non-persistent}. These results can be seen as an extension of the Lajmanovich and Yorke dichotomy in random environment. Despite its apparent simplicity, Theorem \ref{th:BS19} can lead to striking results : indeed, in \cite{BS19}, examples are exhibited where $\Lambda$ is positive even though the disease-free equilibrium is globally asymptotically stable for each individual environment! In other words, if the environment was fixed in one of its possible states, the disease would disappear, but the disease is able to persist due to environmental changes. 
\medskip

Naturally, the next step forward is to combine both demographic and environmental stochasticity. That is, to consider a stochastic individual-based model evolving in a random environment, see e.g. the paper of Cogburn  and Torrez \cite{CT81}. The particular case of epidemiological birth-and-death processes in randomly switched environment has been considered by Artalejo, Economou and Lopez-Herrero \cite{AEL13}, and in several works of Baca{\"e}r, notably \cite{B16} (see the references therein). In \cite{B16}, Baca{\"e}r investigates a monotype SIS stochastic model, with a constant population of size $
K$, evolving in a Markovian environment switching between two states. The process giving the proportion of infectives admits, as a scaling limit when $K$ goes to infinity, a PDMP satisfying \eqref{eq:LYR} in dimension $d=1$. Bacaër focuses on the persistent case, which he calls supercritical, and makes the following conjecture, relying on various heuristic approximations as well as numerical simulations. If the two environments are favourable to the disease (strongly supercritical case), then the mean extinction time starting from a fixed initial proportion of infectives grows exponentially in $K$. In contrast, if one environment is unfavourable to the disease (weakly supercritical case), the time of extinction is of order $K^{\omega}$, for some explicit positive constant $\omega$.	
\medskip


	The aim of this paper is to investigate the behaviour of a family of individual-based multi-type SIS models in random environment whose scaling limit as the size of the population goes to infinity is precisely the PDMP $(X, \Xi)$ described above. Let $d\in\mb{N}^*$, $E$ be a finite non empty set, $(C^\xi)_{\xi\in E}$ be a family of nonnegative irreductible $d\times d$ matrices, $(D^\xi)_{\xi\in E}$ be a family of elements of $(\mb{R}_+^*)^d$, and $Q=(q(\xi,\xi'))_{(\xi,\xi')\in E^2}$ be a transition rate matrix on $E$. For all $K\geq d$, we consider a process $(N^K(t))_{t \geq 0}$ describing the evolution of the number of infectives in a population of size $K$ divided into $d$ groups, with group $i$ of size $K_i\geq 1$. The process $N^K$ takes values in $\prod_{i=1}^d\intbk{0,K_i}$, the $i$-th coordinate counting the number of infectives of type $i$. Together with $N^K$, we consider a process $(\Xi^K(t))_{t \geq 0}$ modelling the evolution of the environment, taking values in a finite set $E$. We assume that $(N^K,\Xi^K)$ is a continuous-time Markov chain on $\prod_{i=1}^d\intbk{0,K_i}\times E$, with three different kinds of transitions possible starting from a state $(n,\xi)$. First, a susceptible of group $i$ may become infected, at rate	
	\[
	(K_i - n_i)\sum_{j=1}^d C_{i,j}^{\xi} \frac{n_j}{K_j},
	\]
	In this formula, $(K_i-n_i)$ is the number of susceptibles in group $i$,  $C_{i,j}^{\xi}$ is the rate at which a given individual in group $i$ contacts people in group $j$, and $n_j/K_j$ is the proportion of infectives in that group. Second, an infective may cure : in group $i$, each infective cures at rate $D_i^{\xi}$, which amounts to a total cure rate equal to $n_i D_i^{\xi}$. Finally, the environment may switch to another state $\xi'$, at rate $q(\xi,\xi')$.
	
	Let us assume that, in the large population limit, each group represents a fixed non-zero fraction of the total population, i.e. $K_i/K\rightarrow \alpha_i>0$ as $K\rightarrow +\infty$ for all $i\in\intbk{1,d}$. Then, as we will show, the scaled process $\hat X^K= (X^K, \Xi^K)$ with
	\[
	X^K(t) = \left(\frac{N^K_1(t)}{K_1},\ldots, \frac{N^K_d(t)}{K_d}\right)
	\]
	converges, as $K$ goes to infinity, to the PDMP $\hat{X}=(X, \Xi)$ where $X$ evolves according to \eqref{eq:LYR}, and $\Xi$ is a continuous-time Markov chain on $E$ with rate matrix $Q$.\medskip

	We denote by $\tau^K_0$ the hitting time of $0$ by $N^K$ (or equivalently $X^K$), which we call the \emph{extinction time}. It corresponds to the extinction of the disease in all the groups. This time is almost surely finite, due to the finiteness of the state space of the Markov chain $(N^K,\Xi^K)$ and the accessibility of $\{0\}\times E$. Our first main result gives bounds on the extinction time starting from any initial condition, as summarised in the following theorem (see Theorems \ref{mintimeextpers}, \ref{majtimeextpers} and \ref{majtimeextnonpers} for more precise statements).
	
\newpage	 
\begin{The}
\label{th:intro1}
Let $\Lambda$ be defined by Equation \ref{eq:deflyapexpointro}. Then,
\begin{enumerate}
\item If $\Lambda > 0$, there exists $p^* \in (0, +\infty]$ such that, for all $p\in(0,p^*)$, there exists $C_1,C'_1 >0$ such that, for all $K$ large enough, all $n\in\mb{N}^*$ and all $t\geq 0$,
\begin{align}\label{lowboundsurvival}
\mb{P} \left( \tau_0^K > t\,\Big|\,\sum_iN^K_i(0)=n\right) \geq \exp\left(-\frac{C'_1}{n^{p}}\right)\exp\left(-\frac{C_1t}{K^{p}}\right).
\end{align}
In particular, conditional on $N^K(0)\neq 0$ we have
\[
\mb{E}(\tau_0^K) \geq \frac{e^{-C'_1}}{C_1}K^p.
\]
\item \label{weaklysup} If $\Lambda > 0$ and $p^* < + \infty$ then for all $p' > p^*$, there exists $C_2,C'_2>0$ such that for all $K$ large enough, all initial conditions and all $t\geq 0$, 
\[
\mb{P}( \tau_0^K > t) \leq C'_2 \exp\left(-\frac{C_2t}{K^{p'}}\right).
\]
In particular, in that case, for all $p < p^* < p'$, for all $K$ large enough, conditional on $N^K(0)\neq 0$ we have
\[
\frac{e^{-C'_1}}{C_1}K^{p} \leq \mb{E}( \tau_0^K) \leq \frac{C'_2}{C_2}K^{p'}.
\]
\item If $\Lambda < 0$, there exists $C_3,C'_3>0$ such that for all $K$ large enough, all initial conditions and all $t\geq 0$,
\[
\mb{P}( \tau_0^K > t) \leq C'_3\exp\left(-\frac{C_3t}{\log(K)}\right).
\]
In particular, 
\[
\mb{E}(\tau_0^K) \leq \frac{C'_3}{C_3}\log(K).
\]
\end{enumerate}
\end{The}
Let us describe and comment the above results. We discuss first the case $\Lambda> 0$, meaning that the limiting PDMP is persistent. The prefactor $\exp(-C'_1/n^p)$ in \eqref{lowboundsurvival} corresponds to a lower bound on the probability of a major outbreak of the epidemic starting from $n$ infectives, meaning that a macroscopic fraction of the population gets infected. We see that this bound does not depend on $K$. Then, starting from a large initial number of infectives, at least for $p>0$ small enough the extinction time $\tau^K_0$ is greater than $K^p$ with high probability.  In some situations ($p^*<\infty$), this extinction time is also smaller than $K^{p'}$ with high probability for $p'$ large enough. What's more, item \ref{weaklysup} shows that when $p^*<\infty$, the mean extinction time is logarithmically equivalent to $K^{p^*}$. The finiteness of $p^*$ depends on the possibility for the linearised PDMP $Y$ to get arbitrarily close to 0 (see Theorem \ref{plyap}) : for example, if the disease-free equilibrium is linearly asymptotically stable in at least one environment, then $p^*$ is finite. 
In the one-dimensional case, the conditions $\Lambda>0,$ $p^*<\infty$ correspond to the weakly supercritical case studied by Baca{\"e}r in \cite{B16} that we mentioned before. We will check in Example \ref{ex:baca} that our $p^*$ is equal to his constant $\omega$ ; hence, our results prove rigorously the conjecture of Baca{\"e}r on mean extinction times.  The fact that the extinction time grows as a power of $K$ when $\Lambda>0$ and $p^*<\infty$ contrasts with the case of a constant, supercriticical environment studied in \cite{CCM16, CCM17}, where the extinction time grows exponentially in $K$ (see \cite[Corollary 3.5]{CCM16} and \cite[Theorem 3.2]{CCM17}). Intuitively, this comes from the fact that when $p^*<\infty$, there is a scenario where $X^K$ is led to reach values of order $1/K$ by environmental switches and then goes extinct, which facilitates the extinction compared to the case where there is only one environment with an endemic equilibrium. Finally, when the scaling limit PDMP is non-persistent ($\Lambda < 0$), the extinction time is of order at most $\log(K)$, which corresponds to the time needed for the limiting process $X$ to reach values of order $1/K$ starting from a macroscopic initial condition.

Altough the disease eventually dies out with probability one, in the persistent case the typical extinction time is very long. From basic results on finite continuous-time Markov chains with an accessible absorbing set (here $\{0\}\times E$) which are irreducible on the complement, the distribution of $\hat{X}^K(t)$ conditional on $\tau^K_0>t$ converges, as $t$ goes to infinity and for any non-zero initial condition, to the unique QSD $\mu^K$ (see Section \ref{sec:limitQSD} for a recall of the definition of this notion). Provided the convergence of conditional distributions is quick with respect to the time scale of extinction, the QSD $\mu^K$ reflects a metastable behaviour of $\hat{X}^K$. The study of the speed of convergence to the QSD is an important issue, but we postpone its study to future works. Our primary focus, which is the object of our second main result, concerns the asymptotic behaviour of the family $(\mu^K)_{K \geq d}$. Recall that a basic property of the QSD is that if $\hat{X}^K(0)\sim\mu^K$ then $\tau_0^K$ has an exponential distribution with rate parameter $\lambda^K \in (0, +\infty)$. We let $\Longrightarrow$ denote the weak convergence of probability measures.
 
 \begin{The}
 \label{th:intro2}
 Let $\Lambda, \mu^*, p^*$ be as defined by Equation \eqref{eq:deflyapexpointro}, Theorem \ref{th:BS19} and Theorem \ref{th:intro1} respectively. 
 \begin{enumerate}
 \item If $\Lambda > 0$, then
 \[
 \mu^K \underset{K\rightarrow +\infty}{\Longrightarrow} \mu^*.
 \]
 Furthermore, 
 \[
 \frac{\log(\lambda^K)}{\log(K)}\underset{K\rightarrow +\infty}{\rightarrow}- p^*.
  \]
\item If $\Lambda < 0$, then
\[
 \mu^K( \cdot \times E)\underset{K\rightarrow +\infty}{\Longrightarrow} \delta_0. 
\]
 
 \end{enumerate}
 \end{The}
Briefly put, this theorem states that, when the limiting PDMP is non-persistent, the QSD $\mu^K$ degenerates into a measure concentrated on the extinction set $\{0\}\times E$ as $K$ goes to $+infty$. On the contrary, if the PDMP is persistent, then $\mu^K$ converges to $\mu^*$, which is concentrated on the survival set $([0,1]^d\setminus\{0\})\times E$. 
Moreover, in the case $p^* < + \infty$, the extinction rate $\lambda^K$ under the QSD is of order $K^{-p^*}$, while if $p^*=\infty$, $\lambda^K$ is rapidly decreasing with $K$, in the sense that $\lambda^K=\mathcal{O}(K^{-p})$ for all $p>0$.

The main difficulty in the proof of the above theorem is the proof of the convergence of $\mu^K$ to $\mu^*$ in the persistent case. A result of Strickler \cite{S19}[Theorem 6.1] shows that every weak limit point of $(\mu^K)_{K\geq d}$ is an invariant probability measure of the limiting PDMP $\hat{X}$. Since by Theorem \ref{th:BS19} $\mu^*$ is the unique invariant probability measure that gives no mass to $\{0\}\times E$, all the problem is to show that $(\mu^K)_{K\geq d}$ is tight on $([0,1]^d\setminus \{0\})\times E$. Compared to the case with a constant environment, this result is harder to establish in random environment  because even though $\hat{X}$ is persistent, it might go arbitrarily close to the extinction set.

%

To the best of our knowledge, this paper is among the first ones to provide a rigorous study of the large population asymptotics related to the extinction time and quasi-stationary distributions for non-branching birth-and-death processes in randomly switched environment.

\medskip

The paper is organised as follows. First, in Section \ref{sec:description}, we construct properly the processes $\hat{X}^K$ and $\hat X$ on the same probability space, give our precise assumptions as well as the convergence rate of $\hat{X}^K$ to $\hat X$. In Section \ref{sec:PDMP}, we slightly complete the analysis of the  PDMP $\hat X$ performed in \cite{BS19}. Our main results are stated precisely in Section \ref{sec:main} and proved in Section \ref{sec:proofs}. Finally, some annex results are proved in the Appendix, Section \ref{sec:appendix}.

\section{Description and construction of the processes}
\label{sec:description}
In this section, we give a precise description of the different processes, in a more general framework than the one presented in the introduction. Indeed, as noticed by Bena{\"i}m and Hirsch \cite{BH99}, the dichotomy exhibited by the Lajmanovich-Yorke model is valid for a wide range of vector fields with some monotonic properties (see below), which enables to choose more general infection and cure rates. We also allow the environmental transition rates to depend (continuously) on the vector of proportion of infectives.

Before proceeding further, let us introduce the following notations and conventions. We let $\mathcal{F}(\mathcal{Y})$ denote the set of real-valued measurable functions on a measurable space $(\mathcal{Y},\scr{Y})$, and $\mathcal{F}_b(\mathcal{Y})$ denote the subset of bounded functions. If $\mathcal{Y}$ is a topological space, then $\mathcal{C}(\mathcal{Y})$ stands for the set of real-valued continuous functions on $\mathcal{Y}$. Given $k,n,p\in\mb{N}^*$, $Z\subset\mb{R}^n$ and $f:Z\to\mb{R}^p$, we say that $f$ is of class $\mathcal{C}^k$ if it can be extended into a function $\tilde{f}$ defined on some open subset of $\mb{R}^n$ containing $Z$, which is of class $\mathcal{C}^k$ in the usual sense. In that case, we assume that such an extension is fixed and write $\partial_j f(z)$, $j\in\llbracket 1,n\rrbracket$ the partial derivatives of $\tilde{f}$ at $z\in Z$ and $Df(z)=\left(\partial_j f_i(z)\right)_{(i,j)\in\llbracket 1,p\rrbracket \times \llbracket 1,n\rrbracket}$ the Jacobian matrix of $\tilde{f}$ at $z$. Moreover, if $g$ is a function defined on $Z\times E$, for all $\xi\in E$ the notation $g^{\xi}$ stands for $g(\cdot,\xi)$ and we say that $g$ is of class $\mathcal{C}^k$ if each $g^{\xi}$ is of class $\mathcal{C}^k$. We let $\norm{\cdot}$ denote the $\ell^{1}$-norm on $\mb{R}^d$ i.e. $\norm{x}=|x_1|+\ldots+|x_d|$ for all $x\in\mb{R}^d$. Given a function $f:E_1\to E_2$, where $(E_2,\norm{\cdot}_{E_2})$ is some normed vector space, we let $\norm{f}_{\infty}:=\sup_{x\in E_1}\norm{f(x)}_{E_2}$. Finally, if $E_1$ is a subset of $\mb{R}^d$, we let $\norm{f}_{\mr{Lip}}$ denote the quantity $\sup\left\{\norm{f(y)-f(x)}_{E_2}/\norm{x-y}\sco x,y\in\mb{R}^d, x\neq y\right\}$.


\subsection{Description and assumptions}\label{sec:desc+ass}

Let $d\in\mb{N}^*$. We assume that we are given sequences of positive integers $(K_i(K))_{K\geq d}$, $i\in\intbk{1,d}$, and $(\alpha_1, \ldots, \alpha_d) \in (0,1)^d$ such that
\begin{equation}
\label{eq:Ki}
\sum_{i = 1}^d K_i(K) = K, \quad \lim_{K \to \infty} \frac{K_i(K)}{K} = \alpha_i > 0.
\end{equation}
In the sequel, to avoid cumbersome notations, we will write $K_i$ instead of $K_i(K)$. The integer $K$ is the total size of the population, $K_i$ is the size of group $i$. The limit in \eqref{eq:Ki} is to be interpreted as follows : for $K$ large enough, the proportion of individuals belonging to group $i$ is approximatively $\alpha_i \in (0,1)$. We set $\underline{\alpha}=\min_{i\in\intbk{1,d}}\alpha_i$.

We let $E=\left\{1,\ldots,|E|\right\}$ with $|E|\geq 2$ denote the set of possible states of the environment. Setting 	$\mathcal{X}=[0,1]^d$ and $\mathcal{X}^K=\mathcal{X}\cap \prod_{i=1}^d \left(K_i^{-1}\mb{N}^{d}\right)$, we consider the state spaces  $\hat{\mathcal{X}}=\mathcal{X}\times E$  and $\hat{\mathcal{X}}^K=\mathcal{X}^K\times E$. In addition, we set $\mathcal{X}_+=\mathcal{X}\setminus \{0\}$, $\mathcal{X}^K_+=\mathcal{X}^K\setminus \{0\}$, $\hat{\mathcal{X}}_+=\mathcal{X}_+\times E$ and $\hat{\mathcal{X}}^K_+=\mathcal{X}^K_+\times E$, which correspond to survival sets of the disease.
	
	We assume that $\hat{X}^K:=(X^K(t),\Xi^K(t))_{t\geq 0}$ is a continuous-time Markov chain on $\hat{\mathcal{X}}^K$ with the following transition rates, for all $x\in \mathcal{X}^K$, $i\in\llbracket 1,d\rrbracket$, and $(\xi,\xi')\in E^2$ with $\xi'\neq \xi$ :

	\begin{center}
	\def\arraystretch{1.5} 
	\begin{tabu}{c c}
	Transition & Rate \\
	$(x,\xi)\rightarrow (x+\cfrac{e_i}{K_i},\xi)$ & $K_i (1-x_i)b_{i}(x,\xi)$ \\
	$(x,\xi)\rightarrow (x-\cfrac{e_i}{K_i},\xi)$& $K_i x_i d_{i}(x,\xi)$ \\
	$(x,\xi)\rightarrow (x,\xi')$ & $q(x,\xi,\xi')$,	
	\end{tabu}
	\end{center}
	where $\left(e_1,\ldots,e_d\right)$ denotes the canonical basis of $\mb{R}^d$, $b_i,d_i:\hat{\mathcal{X}}\to\mb{R}_+$ and $q(\cdot,\xi,\xi'):\mathcal{X}\to\mb{R}_+$. In the following, we put $q(x,\xi,\xi)=-\sum_{\xi'\neq \xi} q(x,\xi,\xi')$ for all $\xi\in E$ and denote by $Q(x)=\left(q(x,\xi,\xi')\right)_{(\xi,\xi')\in E^2}$ the environmental transition rate matrix. Moreover, we denote by $L^K\in\mb{R}^{\mathcal{\hat{X}}^K\times \mathcal{\hat{X}}^K}$ the transition rate matrix of $\hat{X}^K$. Note that the special form of the transition rates does indeed imply that $\hat{X}^K$ cannot leave the state space  $\hat{\mathcal{X}}^K$. The number $b_i(x,\xi)$ represents the total rate of contact and transmission of the disease from infectives to a susceptible individual in group $i$, when the proportion of infectives in each group is $x_i$ and the environment is in state $\xi$. As for $d_i(x,\xi)$, it represents the individual cure rate of an infective in group $i$ in environment $\xi$. In the Lajmanovich-Yorke example given in the introduction, we have $b_i(x,\xi) = \sum_j C_{i,j}^{\xi} x_j$ and $d_i(x,\xi) = D_i^{\xi}$.
	
	We first make a regularity assumption, that stands in the rest of the paper.	
		
	\begin{Ass}[Standing assumption]	\label{ass:reg}
		 For all $\xi\in E$ and $i\in\intbk{1,d}$, the functions $b^{\xi}_i$ and $d^{\xi}_i$ are of class $\mathcal{C}^2$ on $\mathcal{X}$, and the function $x\mapsto Q(x)$ is continuous on $\mathcal{X}$.
	\end{Ass}

	We let $\mathcal{L}^K:\mathcal{F}(\hat{\mathcal{X}}^K)\to \mathcal{F}(\hat{\mathcal{X}}^K)$, $f\mapsto \mathcal{L}^Kf$, denote the infinitesimal generator of $\hat{X}^K$. Letting $\beta_{e_i},\beta_{-e_i}:\hat{\mathcal{X}}\to\mb{R}_+$ be defined by  $\beta_{e_i}(x,\xi)=(1-x_i)b_i(x,\xi)$ and $\beta_{-e_i}(x,\xi)=x_id_i(x,\xi)$, the generator writes
	\begin{align*}
	\mathcal{L}^Kf(x,\xi)=\!\!\sum_{\substack{i\in\llbracket 1,d\rrbracket \\ h\in\left\{-1,1\right\}}}  \!K_i\beta_{he_i}^{\xi}(x)\left[f^{\xi}\hspace{-2pt}\left(x+h\frac{e_i}{K_i}\right)-f^{\xi}(x)\right]\!+\!\sum_{\substack{\xi'\in E \\ \xi'\neq \xi}}q(x,\xi,\xi')\left[f^{\xi'}\hspace{-2pt}(x)-f^{\xi}(x)\right]\!.
	\end{align*}
	For all $f\in\mathcal{C}^1(\hat{\mathcal{X}},\mb{R})$ we have, with a slight abuse of notation,
	\begin{align}\label{cvggen}
	\mathcal{L}^Kf(x^K,\xi)\underset{\substack{K\rightarrow +\infty \\ x^K\rightarrow x}}{\longrightarrow}\mathcal{L} f(x):=\sum_{i=1}^d F_i^{\xi}(x)\partial_if^{\xi}(x)+\sum_{\substack{\xi'\in E \\ \xi'\neq \xi}}q(x,\xi,\xi')\left[f^{\xi'}(x)-f^{\xi}(x)\right], 
	\end{align}
	where the $F_i$ are the coordinates of the $\mathcal{C}^2$-vector field $F:\hat{\mathcal{X}}\to\mb{R}^d$ defined by
	\[
	F^{\xi}(x)=\sum_{i=1}^{d}\left(\beta^{\xi}_{e_i}(x)-\beta^{\xi}_{-e_i}(x)\right)e_i=\sum_{i=1}\left((1-x_i)b_i^{\xi}(x)-x_id_i^{\xi}(x)\right)e_i.
	\]
	Note that since  $F^{\xi}_i(x) \geq 0$ when $x_i = 0$ and   $F_i^{\xi}(x) \leq 0$ when $x_i=1$, the space $[0,1]^d$ is \emph{positively invariant under $F^\xi$}. By this, we mean that for all $x_0\in\mathcal{X}$, the Cauchy problem $\dot{x} = F^{\xi}(x)$, $x(0) = x_0$ admits a solution defined on $\mb{R}_+$ which takes values in $\mathcal{X}$. This solution is unique and denoted by $t\mapsto\psi_t^{\xi}(x_0)$. We call $\psi^{\xi}:\mb{R}_+\times \mathcal{X}\to\mathcal{X}$, $(t,x_0)\mapsto \psi_t^{\xi}(x_0)$ the \emph{semi-flow} induced by $F^\xi$.

	The right handside of \eqref{cvggen} coincides with the expression of the generator $\mathcal{L}$ of a PDMP $\hat{X}=(X(t),\Xi(t))_{t\geq 0}$ with state space $\hat{\mathcal{X}}$ satisfying
	\[
	\begin{cases}\dot{X}(t)=F(\hat{X}(t)) \\ \mb{P}\left(\Xi(t+h)=\xi'\,\big|\,(\hat{X}(s))_{0\leq s \leq t},\ \Xi(t)=\xi\right)=hq(X(t),\xi,\xi')+o(h)\end{cases} 
	\]
	for all $t\geq 0$ and $\xi,\xi' \in E, \xi' \neq \xi$. We refer to $((F^{\xi})_{\xi\in E},Q)$ as the local characteristics of the switched dynamical system $\hat{X}$. The convergence \eqref{cvggen} strongly suggests that $\hat{X}$ is a scaling limit of $\hat{X}^K$ as $K\rightarrow +\infty$. We will see later (Proposition \ref{LLN}) that this is indeed the case.
	
	\medskip
	
	In the following, we say that a square real matrix is \emph{Metzler} if it has non-negative off diagonal entries; such a matrix is said to be \emph{irreducible} if, adding a sufficiently large multiple of the identity, we obtain a non-negative irreducible matrix in the usual sense. We introduce a second set of standing assumptions.
	
	\begin{Ass}[Standing assumption] \label{ass:Metzler}	For all $\xi\in E$ :	 
		\begin{enumerate}[label=\upshape \theAss \textbf{\arabic*.}, ref=\upshape \theAss \textbf{\arabic*}] 
			\item\label{i} $b^{\xi}_i(0)=0$ for all $i\in\intbk{1,d}$;
			\item \label{ii}  $d^{\xi}_i(x) > 0$ for all $i\in\intbk{1,d}$ and $x\in\mathcal{X}$;
			\item \label{Metzler} $\left(\partial_jb^{\xi}_i(x)\right)_{1\leq i,j\leq d}$ is non-negative and irreducible for all $x\in \mathcal{X}$; 
			\item \label{irredenv} $Q(x)$ is irreducible for all $x\in\mathcal{X}$.				
		\end{enumerate}
	\end{Ass} 
	
	It is easily seen that all these hypotheses are satisfied by the model presented in the introduction. They all have a natural interpretation in an epidemiological context. Assumption \ref{i} reflects the fact that there is no external source of infection for the system. Assumption \ref{ii} means that infectives always have the ability to recover regardless of the situation. As for Assumption \ref{Metzler}, it expresses the fact that the apparition of new infectives in one group increases, directly or indirectly, the infection rate of all the other groups. Finally, Assumption \ref{irredenv} means that the environment may transit from any state to any other state while the number of infectives (for $X^K$) stays the same. 
	
	From the mathematical point of view, \ref{i} entails that the extinction set $\{0\}\times E$ is absorbing both for the Markov chain $\hat{X}^K$ and the PDMP $\hat{X}$. Moreover, \ref{ii} ensures that the extinction set is accessible from every state in $\hat{\mathcal{X}}^K$ for $\hat{X}^K$. Together, Assumptions \ref{ii}, \ref{Metzler} and \ref{irredenv} ensure that \emph{$\hat{\mathcal{X}}^K_+$ is an irreducible set of states for the Markov chain $\hat{X}^K$}, for all $K\geq d$. Let us explain why. Thanks to \ref{irredenv}, it is enough to show that $\mathcal{X}^K_+\times \{\xi\}$ is irreducible, for some fixed $\xi$. Assumption \ref{ii} entails that, starting from the state $(\mb{1},\xi)$ where $\mb{1}=(1,\ldots,1)$, one can access every other state in $\mathcal{X}^K_+\times \{\xi\}$. Moreover, for all $i,j\in \intbk{1,d}$ such that $\partial_jb_i^\xi(0)>0$, we have $b^\xi_i(\varepsilon e_j)>0$ for all $\varepsilon>0$ small enough, hence $b^\xi_i(x)>0$ for all $x\in\mathcal{X}$ such that $x_j>0$ since $b^\xi_i$ is non-decreasing with respect to all variables. This entails that, provided $x_j>0$ and $x_i<1$, a new infective may appear in group $i$. Using the irreducibility of the matrix $(\partial_jb_i^\xi(0))_{i,j\in \intbk{1,d}}$, we obtain that $(\mb{1},\xi)$ is accessible from every state in $\mathcal{X}^K_+\times \{\xi\}$, which yields the claimed irreducibility of $\hat{\mathcal{X}}^K_+$. Another important consequence of \ref{Metzler} is the fact that $A^{\xi}:=DF^\xi(0)$ is Metzler and irreducible, given that its off-diagonal entries coincide with the $\partial_jb_i(0)$, $i\neq j$.
	
	\medskip	
	
	Finally, we introduce an additional set of assumptions which are also satisfied by the Lajmanovich-Yorke vector fields, that we will sometimes need, in order to obtain global and ergodic properties for the PDMP $\hat{X}$. Contrary to Assumptions \ref{ass:reg} and \ref{ass:Metzler}, the following hypotheses will stand only when explicitly stated. We let $\mb{R}_{++}^d$ denote the set of elements of $\mb{R}^d$ with positive entries, and for $x,y\in\mb{R}^d$, we write $x\ll y$ if $y-x\in\mb{R}_{++}^d$ and $x\leq y$ if $y-x\in\mb{R}_+^d$.
	
\begin{Ass}
	\label{ass:monotonesubhomo}
	For all $\xi\in E$ :
	\begin{enumerate}[label=\upshape \theAss \textbf{\arabic*.}, ref=\upshape \theAss \textbf{\arabic*}]
		\item\label{it:coop} $F^{\xi}$ is cooperative, meaning that $DF^{\xi}(x)$ is Metzler for all $x \in [0,1]^d$ ;
		\item\label{it:irr} $F^{\xi}$ is irreducible on $[0,1)^d$, meaning that  $DF^{\xi}(x)$ is irreducible for all $x \in [0,1)^d$;
		\item\label{it:sub} $F^{\xi}$ is strongly sub-homogeneous on $(0,1)^d$, i.e. $F^{\xi}(\lambda x) \ll \lambda F^{\xi}(x)$ for all $\lambda > 1$ and all $x \in (0,1)^d$ such that $\lambda x\in(0,1)^d$.
	\end{enumerate}
\end{Ass}

	Note that given \ref{Metzler}, a sufficient condition to have \ref{it:coop} and \ref{it:irr} is that $\partial_jd_i\leq 0$ on $\mathcal{X}$ for all $i,j$, expressing the fact that the apparition of new infectives in one group may only degrade the individual recovery rate of infectives in each group (for instance due to the stress put on the healthcare system). Finally, condition \ref{it:sub} refers to the fact that increasing the proportion of infectives in each group by a factor $\lambda > 1$ results in a smaller increase in the rate of progress of the disease. This can be understood since when the number of infectives is increased by a factor $\lambda > 1$, the number of susceptibles is decreased by a factor $(1 - \lambda x)/(1 - x) < 1$, hence the total number of possible contacts between infectives and susceptibles is likely to be increased by a factor less than $\lambda$. Provided that this effect is not compensated by the decrease of the individual recovery rates (for instance if they are constant as in the Lajmanovich-Yorke model), this yields the strong sub-homogeneity of $F^\xi$. Such a phenomena is also known by the economists as a decreasing returns to scale, meaning e.g. that doubling some factors of  production (here, the number of infectives), does not double the production (the disease).
	
	Assumption \ref{ass:monotonesubhomo} implies that for all $\xi\in E$,  $\psi^{\xi}$ is a strongly monotone, strongly subhomogeneous semi-flow on $\mathcal{X}$, see e.g. the works of Hirsch \cite{Hirsch94,BH99} or Tak\'a\v{c} \cite{T90}. Strong monotonicity means that for all $t>0$ and all $x,y\in\mb{R}_+^d$ such that $x\leq y$, $x\neq y$, we have $\psi^\xi_t(x)\ll \psi^{\xi}_t(y)$; and strong sub-homogeneity means that \ref{it:sub} holds replacing $F^{\xi}$ by $\psi^\xi_t$, $t>0$. These two properties of the semi-flow are sufficient to obtain the same structure of equilibria as in the Lajamnovich-Yorke model : either the disease-free equilibrium is linearly stable for $F^\xi$ in which case it is globally asymptotically stable ; or it is unstable, and there exists an endemic equilibrium $x^\xi_*\in\mb(0,1)^d$, globally asymptotically stable on $\mathcal{X}_+$, see \cite[Theorem 4.2]{BS19}.
	


\subsection{Coupled construction of the processes}\label{sec:construction} For comparison purposes, it is useful to construct the processes $\hat{X}^{K}$ and $\hat{X}$ starting from all possible initial conditions on the same probability space. From now on, we generally write the time variable $t$ in subscript to increase readability. Let $\left(\Omega^{\circ},\scr{F}^{\circ},(\scr{F}^{\circ}_t)_{0\leq t \leq \infty},\mb{P}\right)$ be a filtered probability space satisfying the usual conditions, equipped with a $\left(\scr{F}^{\circ}_t\right)$-Poisson point measure $\scr{N}$ on $\mb{R}_+^2\times \left(\left(\llbracket 1,d \rrbracket \times\left\{-1,1\right\}\right)\sqcup E\right)$ of intensity $\mr{Leb}^{\otimes 2}\otimes \sum_{y\in\left(\llbracket 1,d \rrbracket \times\left\{-1,1\right\}\right)\sqcup E}\delta_{y}$. We let $\scr{N}_{X}$ and $\scr{N}_{\Xi}$ denote the traces of $\scr{N}$ on $\mb{R}_+^2\times \llbracket 1,d \rrbracket \times \left\{-1,1\right\}$ and $\mb{R}_+^2\times E$ respectively. For each $K\geq d$ and $\hat{x}=(x,\xi)\in \hat{\mathcal{X}}$, we let $\hat{X}^{K,\hat{x}}=\left(X^{K,\hat{x}}_t,\Xi^{K,\hat{x}}_t\right)_{t\geq 0}$ solve
\begin{align}			
	X^{K,\hat{x}}_t&= \lfloor x\rfloor_K+\int_{(0,t]\times \mb{R}_+\times \llbracket 1,d\rrbracket\times\left\{-1,1\right\}}\mb{1}_{\left\{u\leq K_i \beta_{he_i}\left(\hat{X}^{K,\hat{x}}_{s-}\right)\right\}}\frac{he_i}{K_i}\scr{N}_{X}(\mr{d}s,\mr{d}u,\mr{d}i,\mr{d}h) \label{eqXK}\\
	\Xi^{K,\hat{x}}_t&= \xi+\int_{(0,t]\times \mb{R}_+\times E}\mb{1}_{\left\{\xi'\neq \Xi^{K,\hat{x}}_{s-}\right\}}\left(\xi'-\Xi^{K,\hat{x}}_{s-}\right)\mb{1}_{\left\{u\leq q\left(\hat{X}^{K,\hat{x}}_{s-},\xi'\right)\right\}}\scr{N}_{\Xi}(\mr{d}s,\mr{d}u,\mr{d}\xi') \label{eqIK}
\end{align}		
		$\mb{P}$-almost surely for all $t\geq 0$, where $\lfloor x \rfloor _K:=(\lfloor K_1x\rfloor/K_1,\ldots,\lfloor K_d x \rfloor/K_d)$. Moreover, we let $\hat{X}^{\hat{x}}=(X^{\hat{x}}_t,\Xi^{\hat{x}}_t)_{t\geq 0}$ solve
		\begin{align}			
			X^{\hat{x}}_t&= x+\int_0^t F(\hat{X}^{\hat{x}}_s)\mr{d}s \label{eqX}\\
			\Xi^{\hat{x}}_t&= \xi+\int_{(0,t]\times \mb{R}_+\times E}\mb{1}_{\left\{\xi'\neq \Xi^{\hat{x}}_{s-}\right\}}\left(\xi'-\Xi^{\hat{x}}_{s-}\right)\mb{1}_{\left\{u\leq q\left(\hat{X}^{\hat{x}}_{s-},\xi'\right)\right\}}\scr{N}_{\Xi}\left(\mr{d}s,\mr{d}u,\mr{d}\xi'\right). \label{eqI}	
		\end{align}
	 The proof of existence and uniqueness for systems \eqref{eqXK}-\eqref{eqIK} and \eqref{eqX}-\eqref{eqI} is given in Section \ref{represtraj}. It is shown that $\hat{X}^{K,\hat{x}}$ is indeed a continuous-time Markov chain of transition rate matrix $L^K$ and that $\hat{X}^{\hat{x}}$ is a switched dynamical system of local characteristics $((F^{\xi})_{\xi\in E},Q)$, and that these two processes satisfy the strong Markov property with respect to $(\scr{F}^{\circ}_t)_{0\leq t \leq \infty}$. Moreover, we can choose the versions of $\hat{X}^{\hat{x}}$, $\hat{x}\in\hat{\mathcal{X}}$ in such a way that $(\hat{x},t,\omega^\circ)\mapsto \hat{X}^{\hat{x}}_t(\omega^{\circ})$ is measurable (with respect to $\scr{B}(\hat{\mathcal{X}})\otimes\scr{B}(\mb{R}_+) \otimes \scr{F}^\circ)$ and that $\mb{P}(\mr{d}\omega^\circ)$-almost surely,  the sample path $t\mapsto\hat{X}^{\hat{x}}_t(\omega^\circ)$ is càdlàg for all $\hat{x}\in\hat{\mathcal{X}}$.  
	
	 In the following, most of the time it will be notationally convenient to drop the exponent associated to the initial condition $\hat{x}$. One way to do this is to work when needed in the extended filtered space $(\Omega,\scr{F},(\scr{F}_t)_{ 0\leq t \leq \infty})$ where $\Omega=\hat{\mathcal{X}}\times \Omega^\circ$, $\scr{F}=\scr{B}(\hat{\mathcal{X}})\otimes\scr{F}^\circ$, $\scr{F}_t=\scr{B}(\hat{\mathcal{X}})\otimes\scr{F}^{\circ}_t
$, $0\leq t\leq \infty$, equipped with the family of probability measures $\mb{P}^{\mu}=\mu\otimes \mb{P}$, $\mu\in\mathcal{P}(\hat{\mathcal{X}})$, and with the processes $(\hat{X}^K_t)_{t\geq 0}$ and $(\hat{X}_t)_{t\geq 0}$ defined by $\hat{X}^K_t(\hat{x},\omega^\circ)=\hat{X}^{K,\hat{x}}_t(\omega^\circ)$ and $\hat{X}_t(\hat{x},\omega^\circ)=X^{\hat{x}}_t(\omega^\circ)$. As usual, we write $\mb{P}^{\hat{x}}=\mb{P}^{\delta_{\hat{x}}}$, $\hat{x}\in\hat{\mathcal{X}}$. For every random element $S$ defined on $(\Omega^{\circ},\scr{F}^\circ)$, we still denote by $S$ its natural extension to $(\Omega,\scr{F})$, writing $S(\hat{x},\omega^\circ)=S(\omega^\circ)$.
		
	Let $\tilde{\scr{N}}_X$ denote the compensated measure associated to $\scr{N}_X$, namely $\tilde{\scr{N}}_X=\scr{N}_X-\mr{Leb}^{\otimes 2}\otimes\left(\sum_{z\in \llbracket 1,d\rrbracket \times\left\{-1,1\right\}}\delta_z\right)$. Equation \eqref{eqXK} yields
	\begin{align}
	X^K_t=X^K_0+\int_0^tF(\hat{X}^K_s)\mr{d}s+M^K_t, \label{eqXKbis}
	\end{align}
	where, for all $\hat{x}\in\mathcal{\hat{X}}$, 
	\[
	M^K_t=\int_{(0,t]\times \mb{R}_+\times \llbracket 1,d\rrbracket\times \left\{-1,1\right\}}\mb{1}_{\left\{u\leq K_i \beta_{he_i}(\hat{X}^K_{s-})\right\}}\frac{he_i}{K_i}\tilde{\scr{N}}_X(\mr{d}s,\mr{d}u,\mr{d}i,\mr{d}h)
	\]
	$\mb{P}^{\hat{x}}$-almost surely for all $t\geq 0$. The process $M^K$ is a $(\scr{F}_t)$-martingale under each $\mb{P}^{\hat{x}}$.
	
	\subsection{Convergence to the PDMP and key estimates near the extinction set}

	As seen above, when $K\rightarrow+\infty$ the generator of $\hat{X}^K$ converges to the one of $\hat{X}$, see \eqref{cvggen}, which strongly suggests that we have convergence in distribution of $\hat{X}^{K,\hat{x}^K}$ to $\hat{X}^{\hat{x}}$ with respect to the Skorokhod topology when $\hat{x}^K\to\hat{x}$. Such a convergence was established in a similar setting by Crudu, Debussche, Muller and Radulescu \cite{CDMR12}, who investigated a set of chemical reactions involving a group of molecules in large number (of order $K$) and a group of molecules in small number (of order one), playing the role of the environment. Here, the coupled construction enables us to prove convergence in probability, with respect to the topology of uniform convergence, see Proposition \ref{LLN} below. 
	
	Although useful, such a convergence is, however, largely insufficient for our purposes. Indeed, when examining questions related to the extinction of $X^K$, one needs much finer estimates about its behaviour near zero. A first crucial estimate is given by Lemma \ref{lemcouplrel} below which controls $\norm{X^K-X}$ relatively to $\norm{X}$. Essentially, this control is good for $\hat{x}=(x,\xi)\in\mathcal{X}^K_+$ such that $\norm{x}\geq a/K$ with $a$ large. It is completed by a second key estimate, given by Lemma \ref{lemboundary}, about the behaviour of $\hat{X}^K$ on the complementary region $0<\norm{x}\leq a/K$.

Let us set $C_F=\max_{\xi\in E}\norm{F^{\xi}}_{\mr{Lip}}$. Let $\mb{1}=(1,\ldots,1)$ and let $\langle\cdot,\cdot\rangle$ denote the standard inner product on $\mb{R}^d$. Recall that $\norm{\cdot}$ denotes the $\ell^{1}$-norm on $\mb{R}^d$. We will repeatedly use the fact that since $F^{\xi}(0)=0$ for all $\xi\in E$, we have 
$|\langle \mb{1},F(\hat{x})\rangle| \leq C_F\norm{x}$ for all $\hat{x}=(x,\xi)\in\hat{\mathcal{X}}$. Since for all $t\in\mb{R}_+$, \eqref{eqX} yields
\[\norm{X^{\hat{x}}_t}=\langle \mb{1},X^{\hat{x}}_t \rangle=\norm{x}+\int_0^t\langle \mb{1},F(\hat{X}^{\hat{x}}_s)\rangle\mr{d}s
\]
this entails that 
\begin{align}\label{encadr}
\norm{x}e^{-C_Ft}\leq \norm{X^{\hat{x}}_t}\leq \norm{x}e^{C_Ft}.
\end{align}
We also define $\mb{\hat{d}}:\hat{\mathcal{X}}\times\hat{\mathcal{X}}\to\mb{R}_+$ by $\mb{\hat{d}}((x,\xi),(x',\xi'))=\norm{x-x'}+\mb{1}_{\xi\neq \xi'}$. This is a distance on $\hat{\mathcal{X}}$, which induces the usual topology. Note that for $\varepsilon \in(0,1)$, $\hat{\mb{d}}((x,\xi),(x',\xi'))>\varepsilon$ if and only if  $\xi\neq \xi'$ or $\norm{x-x'}>\varepsilon$.

The following lemma yields a general and precise control of the gap between $\hat{X}^K$ and $\hat{X}$, starting from the same environment but possibly from different proportions of infectives. It relies on the use of Grönwall's lemma and of Chernoff bounds given by Lemma \ref{lemChernoff} to control the martingale $M^K$ in \eqref{eqXKbis}. The proof is presented in Section \ref{prflemcoupl}.

\begin{Lem}\label{lemcoupl}	
	There exists $C_0,C'_0>0$ and $K_0\geq d$ such that for all $K\geq K_0$, $T>0$, $x\in\mathcal{X}$, $y\in\mathcal{X}^K$, $\xi\in E$ and $\varepsilon>\norm{y-x}$, 
	
	\begin{align}\label{couplgeneral}
		\hspace{-2pt}\mb{P}\hspace{-2pt}\left[\sup_{0\leq t \leq T}\mb{\hat{d}}\left(\hat{X}^{K,(y,\xi)}_t,\hat{X}^{(x,\xi)}_t\right)>\varepsilon\right]&\leq 2d\exp\left(-\frac{K\delta(\varepsilon,T,x,y)}{C_0}\hspace{-2pt}\left(\hspace{-2pt}\frac{\delta(\varepsilon,T,x,y)}{C'_0(Te^{C_FT}\norm{x}+\varepsilon)}\wedge 1\right)\hspace{-2pt}\right)\nonumber \\  &\ +T|E|\sup_{\xi_1\neq \xi_2,\ \norm{z-z'}\leq \varepsilon}\left|q(z,\xi_1,\xi_2)-q(z',\xi_1,\xi_2)\right|.	
	\end{align}	
	where $\delta(\varepsilon,T,x,y)=(\varepsilon e^{-C_F T}-\norm{y-x})_+$.
\end{Lem}

	A first consequence of the above lemma is that it enables to establish, as announced, the convergence in probability of $\hat{X}^K$ to $\hat{X}$ as $K\rightarrow +\infty$ with respect to the topology of locally uniform convergence, provided the initial condition converges. More precisely we prove the following proposition, which in addition yields uniformity of the convergence with respect to the initial condition.
	
	\begin{Pro}[Functional law of large numbers]\label{LLN}
	For all $T,\varepsilon>0$,
	\begin{align}\label{unifcvgproba}
	\sup_{\substack{\xi\in E,\,x\in\mathcal{X},\,y\in\mathcal{X}^K\\ \norm{y-x}\leq h}}\mb{P}\left[\sup_{0\leq t \leq T}\mb{\hat{d}}\left(\hat{X}^{K,(y,\xi)}_t,\hat{X}^{(x,\xi)}_t\right)>\varepsilon\right]\underset{(K,h)\rightarrow (+\infty,0)}{\longrightarrow} 0.
	\end{align}
	In particular, for every $\hat{x}\in\mathcal{\hat{X}}$ and every sequence $(\hat{x}^K)_{K\geq d}\in\prod_{K\geq d}\hat{\mathcal{X}}^K$ such that $\hat{x}^K\rightarrow \hat{x}$ as $K\rightarrow +\infty$,
	\begin{align}\label{cvgproba}
	\hat{X}^{K,\hat{x}^K}\overset{\mb{P}}{\underset{K\rightarrow +\infty}{\longrightarrow }} \hat{X}^{\hat{x}}
	\end{align}
	on $\mathcal{D}(\mb{R}_+,\hat{\mathcal{X}})$ equipped with the (metrisable) topology of locally uniform convergence.
\end{Pro}
\prf Let $T,\varepsilon>0$. For all $\varepsilon'\in(0,\varepsilon)$, 
\[\sup_{\substack{\xi\in E,\,x\in\mathcal{X},\,y\in\mathcal{X}^K\\ \norm{y-x}\leq h}}\delta(\varepsilon',T,x,y)\underset{h\rightarrow 0}{\longrightarrow} \varepsilon' e^{-C_F T}>0.\]
Hence, letting $p_T:\mb{R}_+^*\to[0,1]$ be defined by 
\[
p_{T}(r)=\limsup_{(K,h)\rightarrow (+\infty,0)}\sup_{\substack{\xi\in E,\,x\in\mathcal{X},\,y\in\mathcal{X}^K\\ \norm{y-x}\leq h}}\mb{P}\left[\sup_{0\leq t \leq T}\mb{\hat{d}}\left(\hat{X}^{K,(y,\xi)}_t,\hat{X}^{(x,\xi)}_t\right)>r\right],
\]
applying Lemma \ref{lemcoupl} and using that $\norm{x}\leq d$ for all $x\in\mathcal{X}$, we obtain
\[
p_{T}(\varepsilon)\leq p_T(\varepsilon')\leq T|E|\sup_{\xi_1\neq \xi_2,\ \norm{z-z'}\leq \varepsilon'}\left|q(z,\xi_1,\xi_2)-q(z',\xi_1,\xi_2)\right|.
\]
By uniform continuity of the functions $q(\cdot,\xi_1,\xi_2)$, $\xi_1\neq \xi_2$, letting $\varepsilon'\rightarrow 0$ yields $p_T(\varepsilon)=0$, which ends the proof of \eqref{unifcvgproba}. The convergence \eqref{cvgproba} follows immediately. 
\hfill $\square$

\

Although useful and intellectually satisfying, this functional law of large numbers is not fine enough for our purposes, which require a much more precise control of $X^K$ near $0$. The true aim of Lemma \ref{lemcoupl} is actually to prove the following estimate, which enables to control $\norm{X^K-X}$ relatively to $\norm{X}$.

\begin{Lem}\label{lemcouplrel}
	For all $T>0$, there exists $C_T,\varepsilon_T>0$ such that, for all $K$ large enough, all $\hat{x}=(x,\xi)\in\hat{\mathcal{X}}^K$ and all $\varepsilon\in (0,\varepsilon_T]$,
	\begin{align}
		&\ \mb{P}^{\hat{x}}\left[\sup_{0\leq t \leq T}\norm{X^{K}_t-X_t}>\varepsilon\norm{X_t}\right]\nonumber \\&\leq 2de^{-C_T K \norm{x}\varepsilon^2}+T|E|\sup_{\xi_1\neq \xi_2,\ \norm{z-y}\leq \varepsilon\, e^{-C_F T}\norm{x}}{\left|q(z,\xi_1,\xi_2)-q(z',\xi_1,\xi_2)\right|}.\label{ineqcoupl3}	
	\end{align}			
\end{Lem}	
\prf Let $C'_0,C_0>0$ and $K_0\geq d$ be given by Lemma \ref{lemcoupl}. Let $K\geq K_0$, $\hat{x}=(x,\xi)\in\hat{\mathcal{X}}^K$, $T>0$. Using \eqref{encadr} and Lemma \ref{lemcoupl}, we have, for all $\varepsilon>0$,
\begin{align*}
\mb{P}^{\hat{x}}\hspace{-3pt}\left[\sup_{0\leq t \leq T}\norm{X^{K}_t-X_t}>\varepsilon\norm{X_t}\right]&\hspace{-2pt}\leq \mb{P}^{\hat{x}}\hspace{-3pt}\left[\sup_{0\leq t \leq T}\norm{X^{K}_t-X_t}>\varepsilon e^{-C_F T}\norm{x}\right]\nonumber\\&\hspace{-2pt}\leq 2d\exp\!\left(\!-\frac{K\varepsilon e^{-2C_FT}\norm{x}}{C_0}\hspace{-2pt}\left(\hspace{-2pt}\frac{\varepsilon e^{-2C_F T}}{C'_0(Te^{C_FT}\!+\!\varepsilon e^{-C_F T})}\wedge 1\!\right)\!\right)\nonumber \\&+T|E|\sup_{\xi_1\neq \xi_2,\ \norm{z-y}\leq \varepsilon\, e^{-C_F T}\norm{x}}{\left|q(z,\xi_1,\xi_2)-q(z',\xi_1,\xi_2)\right|}. 
\end{align*}
If we choose $\varepsilon_T\in(0,1]$ small enough, then for all $\varepsilon\in(0,\varepsilon_T]$,  
\[\frac{\varepsilon e^{-2C_F T}}{C'_0(Te^{C_FT}+\varepsilon  e^{-C_F T})}\wedge 1=\frac{\varepsilon e^{-2C_F T}}{C'_0(Te^{C_FT}+\varepsilon  e^{-C_F T})}\geq \frac{\varepsilon e^{-2C_F T}}{C'_0(Te^{C_FT}+ e^{-C_F T})}. 
\]
Hence, we obtain \eqref{ineqcoupl3} with $C_T=e^{-4C_FT}\left(C_0C'_0\left(Te^{C_FT}+e^{-C_FT}\right)\right)^{-1}$.
\hfill $\square$

\

Let us define, for all $K\geq d$ and all $\rho\in\mb{R}_+$, the entrance times
\[
\underline{\tau}^K_{\rho}:=\inf\left\{t\geq 0 : \norm{X^K_t}\leq \rho\right\},\quad \overline{\tau}^{K}_{\rho}:=\inf\left\{t\geq 0 : \norm{X^K_t}\geq \rho\right\}.
\] When $x\mapsto Q(x)$ is not constant, we also need the following lemma, showing that for large $K$, with high probability $X^K$ satisfies bounds similar to \eqref{encadr}, i.e. the relative variation of its norm is not too big. We give the proof in Section \ref{prflemtauq}.

	\begin{Lem}\label{lemtauq}
		Let $T>0$. Set $M=2e^{C_FT}$ and $m=e^{-C_FT}/2$. There exists $C'>0$ and $K_0\in\mb{N}^*$ such that, for all $K\geq K_0$ and  $\hat{x}=(x,\xi)\in\hat{\mathcal{X}}^K_+$, 
		\begin{align}\label{ineqvar}
			\mb{P}^{\hat{x}}\left[\overline{\tau}^{K}_{M\norm{x}}\wedge\underline{\tau}^K_{m\norm{x}}\leq T\right]\leq 2e^{-C'K\norm{x}}.
		\end{align}
	\end{Lem}

The bounds \eqref{ineqcoupl3} and \eqref{ineqvar} involve the terms $2de^{-C_TK\norm{x}\varepsilon^2}$ and $2e^{-C'K\norm{x}}$. For fixed $T,\varepsilon>0$, they are small when $K\norm{x}$ is large. The following lemma, which we prove in Section \ref{prflemboundary}, enables to control the behaviour of $\hat{X}^{K,\hat{x}}$ when $0<K\norm{x}\leq a$.

	\begin{Lem}\label{lemboundary}
		Let $T,a>0$. There exist $c,c'>0$ and $K_0\in\mb{N}^*$ such that for all $K\geq K_0$,
		\begin{align}\label{boundext}
			\inf_{\hat{x}=(x,\xi)\in\hat{\mathcal{X}}^K,\ \norm{x}\leq a/K}\mb{P}^{\hat{x}}\left[\tau^K_0\leq T\right]\geq c
		\end{align}
		and
		\begin{align}\label{boundhita}
			\inf_{\hat{x}=(x,\xi)\in\hat{\mathcal{X}}^K_+,\ \norm{x}\leq a/K}\mb{P}^{\hat{x}}\left[\overline{\tau}^K_{a/K}\leq T\right]\geq c'.
		\end{align}
	\end{Lem}
	
	\section{Complementary study of the PDMP}
	\label{sec:PDMP}
	The  long-term behaviour of $\hat{X}$ has been investigated in \cite{BS19}. One of the main tools to achieve this study is the polar decomposition of $X$, together with the use of rsults of Stochastic Persistence theory in \cite{B18}. We first recall briefly this polar decomposition and the results obtained in \cite{BS19}. Then, in Section \ref{sec:plyap} we present a new result, Theorem \ref{plyap}. It enables to get precise information about the growth rates of $p$-th moments of $\norm{X}$ near the extinction set, for all $p\in\mb{R}$, yielding key Lyapunov and reverse Lyapunov conditions (Proposition \ref{prolyapPDMP}).
	
	\subsection{Polar decomposition and linearised PDMP} In order to understand the behaviour of the PDMP $\hat{X}$ near zero, we will need its polar decomposition. Here we use the $\ell^{1}$-norm, which is particularly well adapted for calculations. Note that results in \cite{BS19} are given for the $\ell^2$-norm, but by equivalence of the norms, they translate immediately for the $\ell^1$-norm. Set $\Delta=\left\{x\in\mb{R}_+^d:\norm{x}=1\right\}$, $\hat{\mathcal{U}}=\left\{(\rho,\theta,\xi)\in \mb{R}_+\times \Delta\times E: \rho\theta\in\mathcal{X}\right\}$ and $\hat{\mathcal{U}}_+=\hat{\mathcal{U}}\setminus \left(\left\{0\right\}\times \Delta\times E\right)$. Let $\hat{x}=(x,\xi)\in\hat{\mathcal{X}}_+$. Since $X^{\hat{x}}$ does not hit $0$, its polar decomposition remains defined for all $t\geq 0$, and \eqref{eqX} yields
	\[
	\norm{X^{\hat{x}}_t}=\norm{x}+\int_0^t\langle \mb{1},F(\hat{X}^{\hat{x}}_s)\rangle \mr{d}s,\qquad \frac{X^{\hat{x}}_t}{\norm{X^{\hat{x}}_t}}=\frac{x}{\norm{x}}+\int_0^t \left( \frac{F(\hat{X}^{\hat{x}}_s)}{\norm{X^{\hat{x}}_s}}-\frac{\langle \mb{1},F(\hat{X}^{\hat{x}}_s)\rangle X^{\hat{x}}_s}{\norm{X^{\hat{x}}_s}^2}\right)\mr{d}s.
	\]
	Let $u=(\rho,\theta)=(\norm{x},x/\norm{x})$, $\hat{u}=(u,\xi)$,  $U^{\hat{u}}_t=\left(R^{\hat{u}}_t,\Theta^{\hat{u}}_t\right)= \left(\norm{X^{\hat{x}}_t},X^{\hat{x}}_t/\norm{X^{\hat{x}}_t}\right)$, $\Xi^{\hat{u}}_t=\Xi^{\hat{x}}_t$ and $\hat{U}^{\hat{u}}_t=\left(U^{\hat{u}}_t,\Xi^{\hat{u}}_t\right)$, $t\geq 0$. The above set of equations can be written
	\begin{align}
	R^{\hat{u}}_t &=\rho+\int_0^t R^{\hat{u}}_s G(\hat{U}^{\hat{u}}_s) \mr{d}s \label{defR}\\
	\Theta^{\hat{u}}_t &=\theta+
	\int_0^t H(\hat{U}^{\hat{u}}_s)\mr{d}s \label{deftheta}
	\end{align}
	where $\tilde{F},G,H$ are defined on $\hat{\mathcal{U}}_+$ by 
	\[
	\tilde{F}(\hat{u})=F(\rho\theta,\xi)/\rho,\quad G(\hat{u})=\langle \mb{1},\tilde{F}(\hat{u})\rangle,\qquad
	H(\hat{u})=\tilde{F}(\hat{u})-G(\hat{u})\theta.\] 
	As a matter of fact, $\tilde{F}$ can be extended on $\hat{\mathcal{U}}$ by setting $\tilde{F}(0,\theta,\xi)=A^{\xi}\theta$. That way, for all $\hat{u}=(\rho,\theta,\xi)\in \hat{\mathcal{U}}$ we have
	\begin{align}\label{extF}
	\tilde{F}(\hat{u})=\int_0^1 \mr{D}F^{\xi}(u\rho\theta) \theta\,\mr{d}u.
	\end{align}
	Since $F$ is $\mathcal{C}^2$ on $\hat{\mathcal{X}}$, $\tilde{F}$ is $\mathcal{C}^1$ on $\hat{\mathcal{U}}$. 
	We extend $G$ and $H$ to $\hat{\mathcal{U}}$ accordingly. In other words, the functions $G_0=G(0,\cdot,\cdot)$ and $H_0=H(0,\cdot,\cdot)$  are given, for all $\hat{\theta}=(\theta,\xi)\in\Delta\times E$, by  
	\[
	G_0(\hat{\theta})=\langle \mb{1},A^{\xi}\theta\rangle,\qquad H_0(\hat{\theta})=A^{\xi}\theta-\langle \mb{1},A^{\xi}\theta\rangle\theta.
	\]
	
	This enables to define $\hat{U}^{(0,\hat{\theta})}$ for all $\hat{\theta}=(\theta,\xi)\in\Delta\times E$. We still set $R^{(0,\hat{\theta})}=\norm{X^{\hat{x}}}\equiv 0$ and $\Xi^{(0,\hat{\theta})}=\Xi^{\hat{x}}$ where $\hat{x}=(0,\xi)$. As for $\Theta^{(0,\hat{\theta})}$, we define it to be the unique solution of \eqref{deftheta}, i.e. 
	\begin{align}\label{thetaH0}
	\Theta^{(0,\hat{\theta})}_t=\theta+\int_0^t H_0\left(\Theta^{(0,\hat{\theta}) }_s,\Xi^{(0,\xi)}_s\right)\mr{d}s.
	\end{align}
	
	That way, \eqref{defR} and \eqref{deftheta} are satisfied for all $\hat{u}\in\hat{\mathcal{U}}$. In what follows we write $\Theta^{\hat{\theta}}$ and $\Xi^{\xi}$   instead of $\Theta^{(0,\hat{\theta})}$ and $\Xi^{(0, \xi)}$, and we set $\hat{\Theta}^{\hat{\theta}}=(\Theta^{\hat\theta},\Xi^{\xi})$. The process $\Xi^\xi$ is a continuous-time Markov chain of transition rate matrix $Q(0)$. The processes $\hat{X}^{\hat{x}}$, $\hat{U}^{\hat{u}}$ and $\hat{\Theta}^{\hat{\theta}}$ are PDMPs generated by switched flows, on the compact state spaces $\hat{\mathcal{X}},\hat{\mathcal{U}}$ and $\hat{\Delta}=\Delta\times E$ respectively. We denote by $(P^{\hat{X}}_t)_{t\in\mb{R}_+}$, $(P^{\hat{U}}_t)_{t\in\mb{R}_+}$ and $(P^{\hat{\Theta}}_t)_{t\in\mb{R}_+}$ the associated semi-groups. As shown\footnote{In \cite{BLMZ} the authors assume for simplicity that the vector fields driving the dynamics of the PDMP between the jumps are smooth. However their proof of the Feller property only requires $\mathcal{C}^1$-regularity, which is satisfied here.} in Proposition 2.1 of \cite{BLMZ}, these semi-groups are Feller. That is, for all $f\in\mathcal{C}(\hat{\mathcal{X}})$ and $t\geq 0$ we have $P^{\hat{X}}_t f\in \mathcal{C}(\hat{\mathcal{X}})$, and $\norm{P^{\hat{X}}_tf-f}_{\infty}\underset{t\rightarrow 0}\longrightarrow 0$, and the same holds replacing $(\hat{X},\hat{\mathcal{X}})$ by $(\hat{U},\hat{\mathcal{U}})$ and $(\hat{\Theta},\hat{\Delta})$.
	
	It turns out that $\Theta^{\hat{\theta}}$ is also the angular process of the linearised version of $X$ at $0$. That is, for $\hat y = (y, \xi) \in \mb{R}^d \times E$, let $(Y_t^{\hat y})_{t \geq 0}$ solve
	\begin{equation}
	\label{eq:pdmplinear}
	Y_t^{\hat y} = y+ \int_0^t A^{\Xi_s^{\xi}} Y_s^{\hat y} \, \mr{d}s.
	\end{equation}
	or, in compact form,
	\[
	\frac{d Y_t^{\hat y}}{dt } =  A^{\Xi_t^{\xi}} Y_t^{\hat y}, \quad  Y_0^{\hat y}=y. 
	\]
    The process $(Y^{\hat{y}},\Xi^\xi)$ is also a PDMP, generated by switched linear flows. It is easily seen that if $y \neq 0$ then $Y^{\hat y}/\|Y^{\hat y}=\Theta^{\hat{\theta}}$ where $\hat{\theta}=(y/\|y\|,\xi)$, and 
    \begin{align}\label{eq:normy}
    \norm{Y^{\hat{y}}_t}=\norm{y}+\int_0^t G_0(\hat{\Theta}^{\hat{\theta}}_s)\mr{d}s.
    \end{align}
   
   \subsection{Top Lyapunov exponent}
   
   The polar decomposition enables us to define rigorously the quantity $\Lambda$ from the introduction. As explained below, it corresponds to the mean growth rate of $\norm{Y}$ under the unique invariant probability measure of $(P^{\hat{\Theta}}_t)_{t\in\mb{R}_+}$.
   
    Assumption \ref{ass:Metzler} and \cite[Proposition 2.13]{BS19} imply the following :
	
	\begin{Pro}
	The semi-group $(P^{\hat{\Theta}}_t)_{t\in\mb{R}_+}$ admits a unique invariant probability measure $\pi_0$. 
	\end{Pro}
	 It follows from the Feller property $(P^{\hat{\Theta}}_t)_{t\in\mb{R}_+}$ and a classical compactness-uniqueness argument that for all $\hat{\theta}\in\Delta\times E$, the empirical measure $A\mapsto t^{-1}\int_0^t \mb{1}_A(\hat{\Theta}^{\hat{\theta}}_s)\mr{d}s$ converges almost surely as $t\rightarrow +\infty$ to $\pi_0$ for the weak topology. Hence, using \eqref{eq:normy} we obtain that, for all $\hat{y}\in\mb{R}_+^d\setminus\{0\}\times E$,
\begin{align}\label{PDMPlinLambda}
	\lim_{t \to \infty} \frac{1}{t} \log \| Y^{\hat y}_t \| = \lim_{t \to \infty} \frac{1}{t}\int_0^t G_0(\hat{\Theta}^{\hat{\theta} }_s)\mr{d}s = \Lambda \quad{a.s.},
\end{align}
where $\Lambda$ is defined by
\begin{align}\label{eq:defLambda}
\Lambda := \pi_0(G_0) = \int_{\hat{\Delta}} G_0(\hat{\theta})\pi_0(\mr{d}\hat{\theta}).
\end{align}
 It is proven in \cite{BS19} that $\Lambda$ coincides with the top Lyapunov exponent in the sense of Oseledets' Multiplicative Ergodic Theorem (see \cite[Theorem 3.4.1]{ards} and \cite[Proposition 2.5]{BS19}). 
The main results in \cite{BS19} state that the asymptotic behaviour of $\hat{X}$ depends only on the sign of $\Lambda$. More precisely, if $\Lambda < 0$, then $X_t$ converges to $0$ exponentially fast with positive probability, while if $\Lambda>0$, then $\hat X$ is stochastically persistent, meaning that $X$ spends an arbitrarily large proportion of time away from $0$. Assuming \ref{ass:monotonesubhomo}, these results can be strengthened to Theorem \ref{th:BS19gen} below (it corresponds to Theorem \ref{th:BS19} of the introduction). We say that a probability measure $\mu\in \mathcal{P}(\hat{\mathcal{X}})$ is \emph{persistent} if $\mu(\{0\}\times E)=0$ (in other words $\mu\in\mathcal{P}(\hat{\mathcal{X}}_+)$).

\begin{The}[Bena\"{i}m and Strickler \cite{BS19}]
\label{th:BS19gen}
Under Assumption \ref{ass:monotonesubhomo}, the following hold.
\begin{enumerate}
\item \label{convpers} If $\Lambda > 0$ and $x \mapsto Q(x)$ is constant,  then $\hat X$ admits a unique persistent stationary distribution and, provided that $X_0\neq 0$ a.s., $(X_t, \Xi_t)$ converges in law to $\mu^*$ as $t\rightarrow +\infty$.
\item \label{convnonpers}If $\Lambda < 0$, then for all $\hat x \in \hat{\mathcal{X}}_+$, 
\[\mb{P}^{\hat x}\left( \limsup_{t \to \infty} \frac{1}{t} \log \|X_t\| \leq \Lambda\right) = 1.
\]
\end{enumerate}
\end{The}

The first item is a consequence of \cite[Theorem 4.12]{BS19}. As for the second one, it essentially follows from Theorem 3.1 of the same paper. We fully justify it in Section \ref{sec:access} just below with Lemma \ref{lem:accessLineartoNonLinear}.

One can also associate a top Lyapunov exponent $\Lambda^{\xi}$ to each fixed environment $\xi\in E$, which corresponds to the principal eigenvalue of $A^\xi$. By principal eigenvalue, we mean that
\[
\Lambda^\xi\in\mr{Sp}(A^\xi)\cap\mb{R}\quad \text{and}\quad  \Lambda^{\xi}= \max\left\{\mr{Re}(\lambda),\ \lambda\in\mr{Sp}(A^\xi)\right\},
\]
which exists thanks to the Perron-Frobenius theorem, using that $A^\xi$ is Metzler. The disease-free equilibrium is linearly (asymptotically) stable for the semi-flow $\psi^\xi$ if and only if $\Lambda^{\xi}$ is (strictly) less than $0$.  What's more, by irreducibility of $A^{\xi}$ there exists a unique right eigenvector $\theta^\xi_*\in \Delta_{++}:=\Delta\cap \mb{R}_{++}^d$
, and for all $y\in(\mb{R}_+^d\setminus \{0\})$, we have 
$\log(\norm{e^{tA^{\xi}}y})/t\rightarrow \Lambda^{\xi}$ and $e^{tA^{\xi}}y/\norm{e^{tA^{\xi}}y}\rightarrow \theta^{\xi}_*$ as $t\rightarrow +\infty$. The analog of \eqref{eq:defLambda} is that $\Lambda^{\xi}=G_0(\theta^\xi_*,\xi)$.
As we mentioned in the introduction, except when $d=1$ there is no way in general to deduce the sign of $\Lambda=\pi_0(G_0)$ from the signs of the $\Lambda^{\xi}=G_0(\theta^\xi_*,\xi)$, $\xi\in E$. When $d\geq 2$, one can have $\Lambda^\xi>0$ for all $\xi$ and $\Lambda<0$, see \cite[Example 4.6]{BS19}; and one can have $\Lambda^\xi<0$ for all $\xi$ and $\Lambda>0$, see \cite[Example 4.7]{BS19}.

\subsection{Accessibility}\label{sec:access}

	Before proceeding further, we need to define the notion of \emph{accessible points}, that appears in the statement of Theorem \ref{plyap} below. Let $M$ be a closed subset of $\mb{R}^d$, let $\mathcal{F}=(\mathcal{F}^{\xi})_{\xi \in E}$ be family of $\mathcal{C}^1$-vector fields from $M$ to $\mb{R}^d$ leaving $M$ positively invariant, and let $(\phi^{\xi})_{\xi\in E}$ the associated family of semi-flows. For  $\mb{i}=(i_1, \ldots, i_m) \in E^m$ and $\mb{u} = (u_1,\ldots,u_m) \in \mb{R}_+^m$, we denote by $\mb{\phi}_{\mb{u}}^{\mb{i}}$ the composite flow : $\mb{\phi}_{\mb{u}}^{\mb{i}} = \phi_{u_m}^{i_m} \circ \ldots \circ \phi_{u_1}^{i_1}$.  For $x \in \mb{R}^d$ and $t \geq 0$, we denote by $\gamma^+_t(x)$ (resp. $\gamma^+(x)$)  the set of points that are reachable from $x$ at time $t$ (resp. at any nonnegative time) with a composite flow: $$\gamma^+_t(x)=\{ \mb{\phi}_{\mb{u}}^{\mb{i}} (x), \: (\mb{i},\mb{u}) \in E^m \times \mb{R}_+^m, m \in \mb{N}, u_1 + \ldots + u_m = t\},$$
$$ 
\gamma^+(x) = \bigcup_{t \geq 0} \gamma^+_t(x).
$$ 
Finally, given a subset $M'$ of $M$, we let the set of $\mathcal{F}$-\emph{accessible points from $M'$} be defined by 
\[
\Gamma(M',\mathcal{F}) = \cap_{x \in M'} \overline{\gamma^+(x)},
\]
where the topological closure is taken in the Alexandroff compactification of $M$. That is to say, if $\gamma^+(x)$ is unbounded for all $x\in M'$ then we write $\infty\in\Gamma(M',\mathcal{F})$. Note that by construction, $\Gamma(M',\mathcal{F})\setminus\{\infty\}$ is positively invariant by each $\mathcal{F}^{\xi}$.
Now, we define three different sets of accessible points by taking $\mathcal{F}$ to be the family of driving vector fields associated to the PDMPs $\hat{X}$, $\hat{\Theta}$ and $\hat{Y}$, by setting 
\[\Gamma(X)=\Gamma\!\left(\mathcal{X}_+,(F^{\xi})_{\xi\in E}\right)\!; \ \Gamma(Y)=\Gamma\!\left(\mb{R}_+^d\!\setminus\!\{0\},(y\mapsto A^{\xi}y)_{\xi\in E}\right)\!;\  \Gamma(\Theta)=\Gamma\!\left(\Delta,(H_0^{\xi})_{\xi\in E}\right)\!.
\]
Moreover, we let $\Gamma(\hat{X})=\Gamma(X)\times E$, $\Gamma(\hat{Y})=\Gamma(Y)\times E$ and $\Gamma(\hat{\Theta})=\Gamma(\Theta)\times E$. This choice of notation is motivated by the fact that $\Gamma(\hat{X})$ coincides with the set of points that are accessible from $\mathcal{X}_+\times E$ for $\hat X$ as a Markov process, see e.g. \cite[Lemme 3.2]{BLMZ}. That is, $\hat{y} =(y,\xi) \in \Gamma(\hat X)$ if and only if, for all neighbourhoods $\mathcal{O}$ of $y$ and all $\hat x \in \hat{\mathcal{X}}_+$, there exists $t\geq 0$ such that $\mb{P}( \hat{X}^{\hat{x}}_t \in \mathcal{O} \times \{ \xi \} ) > 0$. The same property holds replacing $\hat{X}$ by $\hat{Y}$ (resp. $\hat{\Theta}$)  and $\mathcal{X}_+$ by $\mb{R}_+^d\setminus\{0\}$ (resp. $\Delta$).

Next lemma will be very useful in the following.
	\begin{Lem}
		\label{lem:accessLineartoNonLinear}
		Assume \ref{ass:monotonesubhomo}. If $0 \in \Gamma(Y)$ then $0 \in \Gamma(X)$. In particular if $\Lambda<0$ then $0 \in \Gamma(X)$.
	\end{Lem}
	
	\prf
	Assume \ref{ass:monotonesubhomo}. Then for all $\xi\in E$, $\psi^\xi$ is monotone and sub-homogeneous. Hence, for all $x,y \in \mathcal{X}$ such that $x\leq y$ and for all $t\geq 0$, we have $\psi^\xi_t(x)\leq\psi^{\xi}_t(y) \leq h^{-1}\psi^\xi_t(hy)$ for all $0<h<1$, which yields, letting $h\rightarrow 0$, $\psi^\xi_t(x)\leq D\psi^\xi_t(0)y$. Now, by classical results on differential equations, $D\psi^\xi_t(0)y=\phi^\xi_t(y)$ where $\phi^{\xi}$ is the semi-flow of the linear ODE $\dot{z}=A^{\xi}z$. Using this inequality recursively entails that for all $\mb{i} \in E^m$ and $\mb{u} \in \mb{R}^m_+$, $\psi^{\mb{i}}_{\mb{u}}(x) \leq \phi^{\mb{i}}_{\mb{u}}(x)$. 
	
	Now assume what's more that $0 \in \Gamma(Y)$. Let $x \in \mathcal{X}$ and $\varepsilon > 0$. Then, there  exist $\mb{i} \in E^m$ and $\mb{u} \in \mb{R}^m_+$ such that $\| \phi^{\mb{i}}_{\mb{u}}(x)  \| \leq \varepsilon$. Thus, by the previous inequality on the composite flows, $\| \psi^{\mb{i}}_{\mb{u}}(x)  \| \leq \varepsilon$. This implies that $0 \in \Gamma(X)$.
	
	Finally, let us prove the final claim. We only assume that $\Lambda < 0$ and that \ref{ass:monotonesubhomo}holds. Then, it is clear from \eqref{PDMPlinLambda} that $0\in\Gamma(Y)$, and thus $0\in \Gamma(X)$.	
	\hfill $\square$
	\medskip
	
	This lemma enables to fully justify the second item of Theorem \ref{th:BS19gen}. Indeed, Assumption \ref{ass:monotonesubhomo} and $\Lambda<0$ imply that the condition $0\in\Gamma(X)$ is satisfied, which allows to apply \cite[Theorem 3.1]{BS19}.

\subsection{Moment Lyapunov exponents of the linearised PDMP }\label{sec:plyap}

The sign of the Lyapunov exponent $\Lambda$ gives the almost sure stability or unstability of the linear process $Y$. Another classical question is the $p$-moment stability, for some $p \in \mb R$. That is, the asymptotic behaviour of $\mb{E}( \|Y_t\|^p )$. It is classical in the multiplicative ergodic theory to consider the \emph{$p$-moment Lyapunov exponent}, starting from $\hat{y}=(y,\xi)\in(\mb{R}_+^d\setminus \{0\})\times E$,  given by
\[
g(p,\hat y) = \limsup_{ t \to \infty} \frac{1}{t} \log \mb{E}( \|Y_t^{\hat y}\|^p).
\] 
Using again the polar decomposition of $Y$, the $p$-moment Lyapunov can also be expressed as 
\begin{equation}
\label{eq:plyappolar}
g(p,\hat y) = \limsup_{ t \to \infty} \frac{1}{t} \log \mb{E}\left[ \exp \left( p\int_0^t G_0(\hat \Theta^{\hat{\theta} }_s)\mr{d}s \right) \right]
\end{equation}
where $\hat{\theta}=(y/\norm{y},\xi)$.
The $p$-moment Lyapunov exponent for linear stochastic systems has been notably studied by Arnold and his co-authors in the real and the white noise case (see \cite{A84}, \cite{AKO86} and \cite{AOP86}). Under irreducibility type conditions, they prove that $g$ does not depend on $\hat{y}$, that the map $p \mapsto g(p)$ is analytic, convex, and that its derivative at $p=0$ is the top Lyapunov exponent of the linear system under consideration. They use the fact that $g(p)$ can be interpreted as the principal eigenvalue of the generator of an irreducible, compact positive semigroup $(T_t^p)_{t \geq 0}$ which, for $p$ close to $0$ is a perturbation of the initial semigroup of the process on the sphere. 

This consideration was extended to jump linear systems by Leizarowitz \cite{L91} who showed that, in general, $g(p)$ can only be seen as an approximated eigenvalue. He has also proved that, under irreducibility conditions, if $g$ is differentiable at $0$ then its derivative must coincide with the top Lyapunov exponent. It was proven by Fang in its PhD thesis \cite{F94} that the derivative of $g$ from the right is the Lyapunov exponent (see \cite[Theorem 3.4.12]{F94}) and that, in general, $g$ is not differentiable at $0$ (see \cite[Example 3.4.11]{F94}). Moreover, he gives an irreducibility condition under which $g$ is this time differentiable at $0$ (see \cite[Proposition 3.4.11]{F94}). However, this condition is not satisfied in our context.

In the one-dimensional case ($d=1$), a general result from Bardet, Guérin and Malrieu (see  \cite[Propositions 4.1 and 4.2]{BGM}) implies that $g$ is convex, differentiable at $p=0$, and can be explicitly computed as $g(p)=\eta_p$, where $\eta_p$ is the principal eigenvalue of the matrix $Q_p := Q + p \rm{Diag}(A^1, \ldots, A^{|E|})$.

Here, we are able to prove the following theorem, which states that the aforementioned results still hold when switching between Metlzer, irreducible, $d \times d$ matrices. The key tool is the use of the Hilbert projective metric, with respect to which the dynamics of the angular process $\Theta$ is contractive.

%
%
%
%
	\begin{The}\label{plyap}
	\
		
		\begin{enumerate}
		\item \label{item1plyap} For all $p\in\mb{R}$, there exists $g(p)\in\mb{R}$ such that
		\begin{align}\label{plyapmoment}
		g(p) &=\lim_{t\rightarrow +\infty}\frac{1}{t}\log\sup_{\hat{\theta}\in\Delta\times E}\mb{E}\left[\exp\left(p\int_0^t G_0(\hat{\Theta}^{\hat{\theta}}_s)\mr{d}s\right)\right]\nonumber \\&=\lim_{t\rightarrow +\infty}\frac{1}{t}\log\inf_{\hat{\theta}\in\Delta\times E}\mb{E}\left[\exp\left(p\int_0^t G_0(\hat{\Theta}^{\hat{\theta}}_s)\mr{d}s\right)\right].
		\end{align}
		In particular, Equation \eqref{eq:plyappolar} yields that for all $\hat y \in (\mb{R}^d_+ \setminus \{0\}) \times E$,  \[g(p,\hat y) = \lim_{t \to \infty}\frac{1}{t} \log \mb{E}( \|Y_t^{\hat y}\|^p)= g(p).\]  
		\item\label{item2bisplyap} The function $p\mapsto g(p)$ is convex and satisfies 
			\[ \max_{\xi\in E}\Lambda^{\xi} \leq  \lim_{p\rightarrow +\infty}\frac{g(p)}{p}\leq  \max_{\hat \theta \in  \Gamma(\hat \Theta)} G_0(\hat \theta), \qquad \min_{\hat \theta \in \Gamma(\hat \Theta)} G_0(\hat \theta) \leq  \lim_{p\rightarrow-\infty}\frac{g(p)}{p}\leq \min_{\xi\in E}\Lambda^{\xi}.\]
		\item \label{item2plyap} The function $p\mapsto g(p)$ is differentiable at $p=0$ and satisfies
			\[g(0)=0,\qquad g'(0)=\Lambda.
			\]
			 
			\item \label{item3plyap} Assume that $\Lambda > 0$ and let $p^* = \inf\{ p > 0 : \: g(-p) > 0 \}$. Then $p^*$ is finite if and only if $0\in\Gamma(Y)$.
			\item \label{item4plyap} Assume that $\Lambda < 0$ and  let $p_* = \inf\{ p > 0 : \: g(p) > 0 \}$. Then $p_*$ is finite if and only if $\infty \in \Gamma(Y)$.
		\end{enumerate}
	\end{The}

	\begin{Exa}
	\label{ex:baca}
	Consider the one dimensional Lajmanovich and Yorke model in switched environment, studied by Baca{\"e}r and mentioned in the introduction : that is, $E=\{1,2\}$, $b^{\xi}(x) = b^{\xi} x$ and $d^{\xi} = \tilde{d}$ for some  $b^1, b^2, \tilde{d} > 0$, and $Q(x,1,2)=q_1$, $Q(x,2,1)=q_2$ for some $q_1,q_2>0$.  For the associated one dimensional PDMP, it not hard to show that $\Lambda = \frac{q_2}{q_1+q_2}( b^1 - \tilde{d}) +\frac{q_1}{q_2+q_1}( b^2 - \tilde{d})$. In particular, to have $\Lambda > 0$, one needs at least that one of the $b_j$'s, says $b_1$, is greater than $d$. This means that if the environment was fixed in state $1$, the disease would persist in the population (for  the deterministic associated model). In that persistent case ($\Lambda > 0$), Baca{\"e}r distinguishes two regimes : \emph{strongly supercritical}, with $b^1 > b^2 > \tilde{d}$, and  \emph{weakly supercritical}, with $b^1 > \tilde{d} > b^2$. In the strongly supercritical regime, the rate of infection is higher than the cure rate in both environments, while in the weakly supercritical regime, the rate of infection is lower that the cure rate in the second environment. It is thus straightforward that in the strongly supercritical case, $\Gamma(Y)$ does not contain $0$ while it does in the weakly supercritical case. Moreover, using the result of Bardet, Guérin and Malrieu discussed above, we have that $p^*$ is the unique positive solution to $\eta_p = 0$, where $\eta_p$ is the principal eigenvalue of the matrix 
	\[
	Q_p = \begin{pmatrix}
		- q_1 + p( \tilde{d} - b^1) & q_1\\
		q_2 & - q_2 + p(\tilde{d} - b^2)
	\end{pmatrix}.
	\]
	A quick computation gives that 
	\[
	\det(Q_p) = p \left[ p(\tilde{d} - b^1)(\tilde{d} - b^2) - q_1(\tilde{d}-b^2) - q_2(\tilde{d}-b^1) \right],
	\]
	which yields
	\[
	p^* = \frac{q_1(\tilde{d}-b^2) + q_2(\tilde{d}-b^1)}{ (\tilde{d} - b^1)(\tilde{d} - b^2)},
	\]
	or equivalently
	\[
	p^* = \frac{q_2}{\tilde{d}-b^2} - \frac{q_1}{b^1-\tilde{d}}
	\]
	which coincides with the $\omega$ given by Bacaer in \cite{B16}. Note that the assumptions $\Lambda > 0$ and $b^1 > \tilde{d} > b^2$ imply that $p^* > 0$.
\end{Exa}

	\subsection{Lyapunov condition for the nonlinear PDMP}
	
Another important conclusion in \cite{BS19} is the existence of a Lyapunov function when $\Lambda > 0$. Indeed, in that case, it has been proven that for $p > 0$ small enough, there exists $T, \eta > 0$ and $\gamma_0 \in (0,1)$ such that, for all $x$ with $\|x\|\leq \eta$, 
\[
\mb{E}^{\hat{x}}\left[\norm{X_T}^{-p}\right]\leq \gamma_0\norm{x}^{-p}.
\]		
Theorem \ref{plyap} and the Feller regularity of $(P^{\hat{U}}_t)_{t\in\mb{R}_+}$ enable to complete this result and yield the following proposition, which plays a key role in our study.
	\begin{Pro}\label{prolyapPDMP}
		\ 
	\begin{enumerate}
	\item \label{itemprolyap}Let $p \in \mb{R}$ such that $g(p) < 0$. There exist $T>0$, $\gamma_0\in (0,1)$ and $\eta>0$ such that for all $\hat{x}=(x,\xi)\in \hat{\mathcal{X}}_+$ such that $\norm{x}\leq \eta$, 
	\begin{align}\label{ineqlyapPDMP}
	\mb{E}^{\hat{x}}\left[\norm{X_T}^{p}\right]\leq \gamma_0\norm{x}^{p}.
	\end{align}
	\item \label{itemprolyapinv}Let $p \in \mb{R}$ such that $g(p) > 0$. There exist $T>0$, $\tilde{\gamma}_0>1$ and $\eta>0$ such that for all $\hat{x}=(x,\xi)\in \hat{\mathcal{X}}_+$ such that $\norm{x}\leq \eta$, 
		\begin{align}\label{ineqlyapinvPDMP}
			\mb{E}^{\hat{x}}\left[\norm{X_T}^{p}\right]\geq \tilde{\gamma}_0\norm{x}^{p}.
		\end{align}
	\end{enumerate}

	\end{Pro}
	\prf
We only prove the first point, the proof of the second one is similar.	Let $p \in \mb{R}$ such that $g(p) < 0$. For all $\hat{u}=(\rho,\theta,\xi)\in\hat{\mathcal{U}}$, \eqref{defR} yields, for all $t\geq 0$,
	\begin{align}
		R^{\hat{u}}_t= \rho\exp\left(\int_0^t G(\hat{U}^{\hat{u}}_s)\mr{d}s\right) \label{solveqR}
	\end{align}
	hence 
	\begin{align}\label{eqRp}
		\mb{E}\left[(R^{\hat{u}}_t)^{p}\right]=\rho^{p}\mb{E}\left[\exp\left(p\int_0^t G(\hat{U}^{\hat{u}}_s)\mr{d}s\right)\right]
	\end{align}

	We claim that for all $t\geq 0$ the function 
	\[
	W_{p,t}\colon\hat{u} \mapsto \mb{E}\left[\exp\left(p\int_0^t G(\hat{U}^{\hat{u}}_s)\mr{d}s\right)\right]
	\]
	is continuous on $\hat{\mathcal{U}}$. Indeed, the Feller property of $(P^{\hat{U}}_t)_{t\in\mb{R}_+}$ entails that the probability distribution of $\hat{U}^{\hat{u}}$ is continuous with respect to $\hat{u}$, where $\mathcal{P}\left(\mathcal{D}(\mb{R}_+,\hat{\mathcal{U}})\right)$ is equipped with the weak topology associated to the Skorokhod topology (this fact is a corollary of Theorem 2.5 in Chapter IV of \cite{EthKur}). Now, since convergence in $\mathcal{D}(\mb{R}_+,\hat{\mathcal{U}})$ implies convergence Lebesgue-almost everywhere, it follows from the continuity of $G$ and the dominated convergence theorem that for all $t\geq 0$, $u\mapsto \exp\left(p\int_0^t G(u_s)\mr{d}s\right)$ is a bounded continuous functional on $\mathcal{D}(\mb{R_+},\hat{\mathcal{U}})$, yielding the continuity of $W_{p,t}$. 
		
		Given that $g(p) < 0$,  it follows from Proposition \ref{plyap} that we can choose $T>0$ such that 
		\[
		\gamma(T):=\sup_{\hat{\theta}\in \Delta\times E}W_{p,T}(0,\hat{\theta})=\sup_{\hat{\theta}\in \Delta\times E}\mb{E}\left[\exp\left(p\int_0^tG_0(\hat{\theta}^{\hat{\theta}}_s)\mr{d}s\right)\right]<1.
		\]
		By continuity of $W_{p,T}$, we can find $\eta>0$ small enough so that
		\[
		\gamma_0:=\sup_{\hat{u}\in [0,\eta]\times \Delta\times E}W_{p,T}(\hat{u})<1.
		\]
		Combining this with \eqref{eqRp}, we obtain that \eqref{ineqlyapPDMP} holds for all $\hat{x}\in\hat{\mathcal{X}}$ such that $\norm{x}\leq \eta$.
	\hfill $\square$
	
\medskip

Proposition \ref{prolyapPDMP} enables to control the tail near $0$ of persistent stationary distributions of $\hat X$:	
	\begin{Cor}\label{cor:momentmu}
	Assume that $\Lambda > 0$, and let $\mu$ be a persistent stationary distribution of $\hat X$. Then :
	\begin{enumerate}
	\item For all $p\in(0,p^*)$, 
	\[
	\int_{\hat{\mathcal{X}}_+} \norm{x}^{-p} \, \mu(\mr{d}x,\mr{d}\xi) < + \infty.
	\]
	\item If $0 \in \Gamma(Y)$ and Assumption \ref{ass:monotonesubhomo} holds, then for all $p > p^*$,
	\[
	\int_{\hat{\mathcal{X}}_+} \norm{x}^{-p} \, \mu(\mr{d}x,\mr{d}\xi) = + \infty.
	\]
	\end{enumerate}
	\end{Cor}
\prf
For all $p\in(0,p^*)$ we have $g(-p)<0$, hence Proposition \ref{prolyapPDMP} entails that $(x,\xi)\mapsto \|x\|^{-p}$ is a Lyapunov function for $\hat X$. Then, it follows from classical arguments (see e.g. Proposition 4.24 in \cite{H06}) that this function is integrable with respect to $\mu$. 

Now, assume that $0 \in \Gamma(Y)$ (i.e. $p^*<\infty$) and that \ref{ass:monotonesubhomo} holds. Fix $p > p^*$. We have $g(-p)>0$ and we let $T, \eta >0$ and $\tilde{\gamma}_0>1$ be given by Proposition \ref{prolyapPDMP} applied to $-p$. Let $\varphi:\mathcal{X}\to\mb{R}_+$ be defined by $\varphi(x,\xi)= \mb{1}_{\norm{x}>0}(\|x\|^{-p} \vee \eta^{-p})$. A straightforward consequence of the definition of $\varphi$ and Proposition \ref{prolyapPDMP} is that for all $\hat x=(x,\xi) \in \hat{\mathcal{X}}$,
\[
\mb{E}^{\hat x} \left( \varphi(\hat{X}_T) \right) \geq \tilde{\gamma}_0 \varphi(\hat{x})\mb{1}_{\|x\|\leq \eta} +  \varphi(\hat{x})\mb{1}_{\|x\|> \eta}.
\]  
Integrating the above equation with respect to $\mu$ and using the invariance of $\mu$ yields
\[
\mu \varphi \geq \mu \varphi + (\tilde{\gamma}_0 - 1 ) \mu ( \varphi \mb{1}_{ \{ \hat x\in\hat{\mathcal{X}}_+ : \: \|x\| \leq \eta\}} ).
\]
Now, due to Lemma \ref{lem:accessLineartoNonLinear}, $0 \in \Gamma(X)$, and therefore $\{0\} \times E$ is in the topological support of $\mu$ by Proposition 3.17 in \cite{BLMZ}. Since $\tilde{\gamma}_0 > 1$ and $\varphi > 0$ on $\hat{\mathcal{X}}_+$, this implies that the second term in the right handside of the above inequality is strictly positive. Hence, we must have $\mu \varphi = + \infty$, which entails that $\int_{\hat{\mathcal{X}}_+} \norm{x}^{-p} \, \mu(\mr{d}x,\mr{d}\xi) = + \infty$.
\hfill $\square$

	\section{Main results}
	\label{sec:main}
	Our main results concern the behaviour of $\hat{X}^K$ for large $K$, in two aspects. In Section \ref{sec:mainexttime} we give stochastic bounds on the extinction time of the disease, starting from any initial condition. Then, in Section \ref{sec:limitQSD} we describe the asymptotic behaviour of the unique quasi-stationary distribution of $\hat{X}^K$, as $K$ goes to infinity. The results are very different depending on whether the limiting PDMP is persistent or not.		
	
	\subsection{Asymptotics for the extinction time}\label{sec:mainexttime}
	
	As mentioned in Section \ref{sec:desc+ass}, Assumption \ref{ass:Metzler} entails that $\{0\}\times E$ is accessible from every state in $\hat{\mathcal{X}}^K$ for the Markov chain $\hat{X}^K$.  Due to the finiteness of $\hat{\mathcal{X}}^K$, this ensures that the extinction time of the disease
	\[
	\tau_0^K = \inf \{ t \geq 0 : \: X^K_t = 0 \},
	\]
	is finite almost surely, under every $\mb{P}^\mu$. A natural question is to know what is the order of magnitude of $\tau_0^K$. Our first main results, Theorems \ref{mintimeextpers}, \ref{majtimeextpers} and \ref{majtimeextnonpers} below, show that the answer to that question depends on whether the limit process $\hat{X}$ is persistent or not. Together, these three theorems form Theorem \ref{th:intro1}, stated in the introduction in the simplified framework.	
	
\subsubsection{The persistent case}

We begin with the case where the limiting PDMP $\hat X$ is persistent, i.e. $\Lambda > 0$. Our first theorem yields two results. First, we obtain a lower bound, independent of $K$, on the probability that a macroscopic fraction of the population gets infected, starting from any initial condition with at least one infective. Second, our theorem yields a stochastic lower bound on the extinction time of the epidemic. When the  initial number of infectives is large i.e. $K\norm{X^K_0}\gg 1$, the extinction time is greater than $K^p$ with high probability, for all $p<p^*$.

\begin{The}\label{mintimeextpers}
		Assume that $\Lambda>0$, and let $p\in (0,p^*)$. There exist $\eta,C_1,C'_1,C''_1>0$ and $K_0\in\mb{N}^*$ such that, for all $K\geq K_0$, all $\hat{x}=(x,\xi)\in\hat{\mathcal{X}}^K_+$ and all $t\geq 0$,
		\[
		\mb{P}^{\hat{x}}\left(\overline{\tau}^K_\eta<\infty\right)\geq \exp\left(-\frac{C''_1}{(K\norm{x})^p}\right)\geq e^{-C''_1}.
		\]
		and
		\[
		\mb{P}^{\hat{x}}\left(\tau^K_0>t\right)\geq \exp\left(-\frac{C'_1}{(K\norm{x})^p}\right)\exp\left(-\frac{C_1t}{K^p}\right).
		\]
		In particular, for all $K\geq K_0$, 
		\[	\inf_{\hat{x}\in\hat{\mathcal{X}}^K_+}\mb{E}^{\hat{x}}(\tau^K_0) \geq\frac{e^{-C'_1}}{C_1}K^p.
		\]	
		
	\end{The}

We prove this theorem in Section \ref{sec:mintimeextpers}. Let us explain the main ingredients of the proof. The first crucial step consists in transfering to $\hat X^K$ the Lyapunov property of $\hat x \mapsto \| x \|^{-p}$ for $p\in(0,p^*)$ given by Proposition \ref{prolyapPDMP} for $\hat X$. More precisely, we have the following proposition.
	\begin{Pro}\label{prolyapXK:1}
	Let $p  \in \mb{R}$ such that $g(p) < 0$. Let $T>0$ and $\eta>0$ be given by Proposition \ref{prolyapPDMP}. There exist $\gamma\in(0,1)$, $a>0$ and $K_0\in\mb{N}^*$ such that for all $K\geq K_0$ and all $\hat{x}=(x,\xi)\in\hat{\mathcal{X}}_+^K$ satisfying $a/K\leq \norm{x} \leq \eta$, we have
	\begin{align}
		\mb{E}^{\hat{x}}\left[\norm{X^K_T}^{p}\mb{1}_{\left\{T<\tau_0^K\right\}}\right]\leq \gamma\, \norm{x}^{p}.
	\end{align}
\end{Pro}

The proof relies on the use of Lemma \ref{lemcoupl}, which entails that $\norm{X^K}/\norm{X}$ stays close to $1$ with high probability on $[0,T]$, as well as the control of $X^K$ given by Lemma \ref{lemtauq}. As we see below, it is crucial that this Lyapunov condition holds for $\norm{x}$ as small as $a/K$.  

Thanks to the above proposition, we can use a supermartingale argument to bound from above the probability for $\norm{X^K}$ to reach low barriers before higher ones. 

\begin{Lem}\label{hittingtime1:1}
	Assume that $\Lambda>0$ and let $p\in (0,p^*)$. Let $\eta>0$ be given by Proposition \ref{prolyapPDMP}, item \ref{itemprolyap}. There exist $a'>0$, $C_5\geq 1$ and $K_0\in\mb{N}^*$ such that for all $K\geq K_0$, $\rho>0$ and $\hat{x}=(x,\xi)\in\hat{\mathcal{X}}^K_+$ such that $a'/K\leq \rho \leq \norm{x}\leq \eta$, we have
	\[
	\mb{P}^{\hat{x}}\left(\underline{\tau}^K_{\rho}<\overline{\tau}^{K}_{\eta}\right)\leq C_5 (\norm{x}/\rho)^{-p}.
	\]
\end{Lem}

Let us fix some $\eta'<\eta$. Intuitively, assuming that $\norm{X^K_0}\geq \eta$, the trajectory of $X^K$ before its norm goes below some level $\rho$ astisfying $a/K\leq \rho<\eta'$ can be decomposed into successive excursions of the following form. First, wait for $\norm{X^K}$ to go below $\eta'$. Then, wait for it to either exceed $\eta$ again, in which case we call it a failed excursion and repeat the procedure, or to go below $\rho$, in which case we call it successful excursion and the decomposition is over. By means of the above lemma, we can obtain a lower bound on the probability of a failed excursion. In addition, the duration of a failed excursion can be bounded from below with high probability by some constant $T'$, using Lemma \ref{lemtauq} which shows that $\norm{X^K}$ varies slowly. This decomposition yields a lower bound on $\underline{\tau}^K_{\rho}$, as stated by the following lemma.

\begin{Lem}\label{hittingtime2:1}
	Assume that $\Lambda>0$, and let $p\in (0,p^*)$. Let $\eta>0$  be given by Proposition \ref{prolyapPDMP}. There exist $T', a'>0, c_1>0$, $\delta\in(0,\eta)$ and $K_0\in\mb{N}^*$ such that for all $K\geq K_0$, $\rho\in[a'/K,\delta]$, $\hat{x}=(x,\xi)\in \hat{\mathcal{X}}^K_+$ such that $\norm{x}\geq \eta$ and $t\geq 0$, we have
	\[
	\mb{P}^{\hat{x}}\left(\underline{\tau}^K_{\rho}>t\right)\geq e^{-c_1\rho^p\lceil t/T'\rceil}.
	\]
\end{Lem}

This yields a lower bound on the extinction time starting from a macroscopic fraction of infectives, using that $\tau^K_0\leq \underline{\tau}^K_{a/K}$. Here we see that in order to get the time scale $K^p$ in Theorem \ref{mintimeextpers}, it is crucial to be able to take $\rho$ of order $1/K$. Finally, Theorem \ref{mintimeextpers} follows from the above lemmas together with Lemma \ref{lemboundary}, which yields a lower bound on the probability that the total number of infectives reaches some prescribed value $a''$.

\begin{Rem}
\label{rem:subexp}
	In case $\Lambda >0$ and  $0 \notin \Gamma(Y)$, we have $p^* = + \infty$ hence Theorem \ref{mintimeextpers} implies that the typical extinction time starting from a large initial number of infectives (i.e. $K\norm{x}\gg 1$) grows faster than any power of $K$. Thus, one may wonder if the extinction time grows exponentially fast with $K$ in that case, as for the classical stochastic SIS epidemic model in constant environment in the supercritical case, and numerous other models where 0 is an unstable equilibrium of the limiting ODE, see e.g. \cite{DSS05,CCM16,CCM17}. However, \emph{this is not case in general}.  Consider the context of Example \ref{ex:baca}, that is, $E=\{1,2\}$, $b^{\xi}(x) = b^{\xi} x$ and $d^{\xi} = \tilde{d}$ for some  $b^1, b^2, \tilde{d} > 0$, and $Q(x,1,2)=q_1$, $Q(x,2,1)=q_2$ for some $q_1,q_2>0$. Assume that $b^1 > \tilde{d}$ and that $b^2 = \tilde{d}$. This implies that $\Lambda = \frac{q_2}{q_1+q_2}( b^1 - \tilde{d})> 0$. Moreover, in that case, $A^1 = b^1 - \tilde{d} > 0$ and $A^2 = 0$. In particular, $0$ is not in $\Gamma(Y)$. Let $\hat{x}=(x,\xi)\in\hat{\mathcal{X}}^K_+$ be an initial condition for $\hat{X}^K$. First, the environment $\Xi^K$ reaches state $2$ before time $1$ with probability at least $1-e^{-q_1}$. Then, as long as it stays in state $2$, $KX^K$ coincides with a critical stochastic SIS process, whose mean extinction time, starting from a maximal initial number of infectives, is asymptotically equivalent to $C\sqrt{K}$ for some $C>0$, see \cite[Section 2.2]{DSS05}. What's more, the probability that the environment stays in state $2$ during a period of time longer than $t$ is equal to $e^{-q_2t}$. Therefore, for $K$ large enough we have 
	\begin{align*}
		&\ \mb{P}^{\hat{x}}\left(\tau^K_0\leq 2C\sqrt{K}\right)\\
		&\geq (1-e^{-q_1})\mb{P}^{(x,2)}\left(\tau^K_0\leq 2C\sqrt{K}-1\right)\\
		&\geq (1-e^{-q_1})e^{-q_2(2C\sqrt{K}-1)}\mb{P}^{(x,2)}\left(\tau^K_0\leq 2C\sqrt{K}-1\,|\,\forall t\in[0,2C\sqrt{K}-1],\, \Xi^K_t= 2\right)\\
		&\geq C'e^{-2q_2C\sqrt{K}}
	\end{align*}
	where $C'=(1-e^{-q_1})e^{2q_2C}/3$, using Markov's inequality. Then the use of the Markov property shows that for all $\hat{x}\in\hat{\mathcal{X}}^K_+$, $\mb{P}^{\hat{x}}(\tau^K_0>t)\leq \left(1-C'e^{-2q_2C\sqrt{K}}\right)^{\left\lfloor t/(2C\sqrt{K})\right\rfloor}$, which entails that $\mb{E}^{\hat{x}}(\tau^K_0)\leq 2CC'^{-1}\sqrt{K}e^{2q_2C\sqrt{K}}$. Hence, it is possible to have $\Lambda>0$, $p^*=+\infty$ and a mean time of extinction which is subexponential in $K$.

\end{Rem}

We complete the lower bound on the extinction time given by Theorem \ref{mintimeextpers} by an upper bound. In the case where $0 \in \Gamma(Y)$ i.e. $p^*<\infty$, we prove that for all $p>p^*$, the extinction time is smaller than $K^p$ with high probability, under Assumption \ref{ass:monotonesubhomo}.
\begin{The}\label{majtimeextpers}
Assume that $\Lambda>0$, $0 \in \Gamma(Y)$ and that Assumption \ref{ass:monotonesubhomo} holds. Let $p>p^*$. There exist $C_2,C'_2>0$ and $K_0\in\mb{N}^*$ such that for all $K\geq K_0$, all $\hat{x}\in\hat{\mathcal{X}}^K_+$ and all $t\geq 0$,
	\[
	\mb{P}^{\hat{x}}\left(\tau^K_0>t\right)\leq C'_2\exp\left(-\frac{C_2t}{K^p}\right).
	\]
In particular, for all $K\geq K_0$,
\[
\sup_{\hat{x}\in\hat{\mathcal{X}}^K_+}\mb{E}^{\hat{x}}(\tau^K_0) \leq\frac{C'_2}{C_2}K^p.
\]

\end{The}
The two above theorems yield the following control on the mean extinction time, uniform in the initial condition. Under Assumption \ref{ass:monotonesubhomo}, the mean extinction time is logarithmically equivalent to $K^{p^*}$ if $p^*<\infty$, and grows faster than any power of $K$ if $p^*=\infty$.
\begin{Cor}\label{theordertimeext}
	Assume that $\Lambda>0$. Then, 
		\begin{align}\label{liminftimeext}
		\liminf_{K\rightarrow +\infty}\frac{\log \inf_{\hat{x}\in \hat{\mathcal{X}}^K_+}\mb{E}^{\hat{x}}(\tau^K_0)}{\log(K)}\geq p^*.
		\end{align}
	If, in addition, Assumption \ref{ass:monotonesubhomo} holds, then
		\begin{align*}
		\liminf_{K\rightarrow +\infty}\frac{\log \inf_{\hat{x}\in \hat{\mathcal{X}}^K_+}\mb{E}^{\hat{x}}(\tau^K_0)}{\log(K)}=\limsup_{K\rightarrow +\infty}\frac{\log \sup_{\hat{x}\in \hat{\mathcal{X}}^K_+}\mb{E}^{\hat{x}}(\tau^K_0)}{\log(K)}=p^*.
		\end{align*}
\end{Cor}

\prf Assume that $\Lambda>0$. Applying Theorem \ref{mintimeextpers}, letting $K\rightarrow +\infty$ and then $p\rightarrow p^*$ yields \eqref{liminftimeext}.
If $0\in \Gamma(Y)$, which is equivalent to $p^*=+\infty$, then the proof is complete. If $0\notin \Gamma(Y)$ and \ref{ass:monotonesubhomo} holds, then Theorem \ref{majtimeextpers} entails
\[
\limsup_{K\rightarrow +\infty}\frac{\log \sup_{\hat{x}\in \hat{\mathcal{X}}^K_+}\mb{E}^{\hat{x}}(\tau^K_0)}{\log(K)}\leq p^*,
\]
which concludes the proof.
\hfill $\square$

The proof of Theorem \ref{majtimeextpers} is given in Section \ref{sec:majtimeextpers}. Its ingredients are very similar to those of the proof of Theorem \ref{mintimeextpers}. This time, we transfer to the process $\hat{X}^K$ the "reverse" Lyapunov property verified by $\hat x \mapsto \|x\|^{-p}$ with respect to $\hat X$ when $p>p^*$.
\begin{Pro}\label{prolyapinvXK:1}
	Assume that $\Lambda>0$ and $0\in\Gamma(Y)$, and let $p>p^*$. Let $T>0$ and $\eta>0$ be given by item \ref{itemprolyapinv} of Proposition \ref{prolyapPDMP} applied to $-p$. There exist $\tilde{\gamma}>1$, $a>0$ and $K_0\in\mb{N}^*$ such that for all $K\geq K_0$ and all $\hat{x}=(x,\xi)\in\hat{\mathcal{X}}_+^K$ satisfying $a/K\leq \norm{x}\leq \eta$,
	\begin{align}
		\mb{E}^{\hat{x}}\left[\norm{X^K_T}^{-p}\mb{1}_{\left\{T<\tau_0^K\right\}}\right]\geq \tilde{\gamma}\, \norm{x}^{-p}.
	\end{align}
\end{Pro}

Using a submartingale argument, combined with Lemma \ref{lemboundary} and the fact that $0\in\Gamma(X)$ (Lemma \ref{lem:accessLineartoNonLinear}), we can show that the probability of extinction before a time of order $\log(K)$ is of order at least $K^{-p}$, uniformly on the starting point in $\hat{\mathcal{X}}^K_+$. Theorem \ref{majtimeextpers} then follows from the Markov property.

\subsubsection{The non-persistent case}
We conclude on the extinction time with the non-persistent case. In this situation, the extinction time is of order at most $\log(K)$. 
\begin{The}\label{majtimeextnonpers}
	Assume that $\Lambda<0$ and that \ref{ass:monotonesubhomo} holds. There exist $K_0\in\mb{N}^*$, $C_3,C'_3>0$ such that for all $K\geq K_0$, all $\hat{x}\in\hat{\mathcal{X}}^K_+$ and all $t\geq 0$, 
	\[
	\mb{P}^{\hat{x}}\left(\tau^K_0>t\right)\leq C'_3\exp\left(-\frac{C_3t}{\log(K)}\right).
	\]
	In particular, for all $K\geq K_0$,
	\[	\sup_{\hat{x}\in\hat{\mathcal{X}}^K_+}\mb{E}^{\hat{x}}(\tau^K_0) \leq\frac{C'_3}{C_3}\log(K).
	\]
\end{The}

We give the proof in Section \ref{sec:majtimeextnonpers}. Since $g'(0)=\Lambda<0$ (see Theorem \ref{plyap}), we can choose $p>0$ such that $g(p)<0$ and make use of Proposition \ref{prolyapXK:1}. By a supermartingale argument and using once that $0\in\Gamma(X)$ thanks to Lemma \ref{lem:accessLineartoNonLinear}, we are able to bound from below by a constant the probability of extinction before a time of order $\log(K)$, uniformly on the starting point in $\hat{\mathcal{X}}^K_+$. Theorem \ref{majtimeextnonpers} then follows from the Markov property.
	
	\subsection{Scaling limit of quasi-stationary distributions}\label{sec:limitQSD}
	
	As we know, with probability one the disease eventually dies out. A good way to gain understanding on the behaviour of $\hat{X}^K$ before extinction is to study its \textit{quasi-stationary distribution(s)} (QSD). A probability measure $\mu^K$ on $\hat{\mathcal{X}}^K_+$ is a QSD for $\hat X^K$ if and only if, for all $t \geq 0$,
	\[
	\mb{P}^{\mu^K} \left( \hat{X}_t^K \in \cdot \, \left|  \, \tau_0^K > t \right.\right) = \mu^K( \cdot).
	\]
	Due to the finiteness of the state space $\hat{\mathcal{X}}^K$ and the irreducibility of $\hat{X}^K_+$ for the continuous-time Markov chain $\hat{X}^K$, there exists a unique QSD $\mu^K$ for $\hat{X}^K$, by classical arguments based on the Perron-Frobenius theorem. Furthermore, the distribution of $\hat X^K_t$ conditional on non-extinction converges exponentially fast to $\mu^K$, in the sense that for some positive constants $C(K), \gamma(K)$,
	\begin{align}\label{cvgtoqsd}
	\left\| \mb{P}^{\nu} \left( \hat{X}_t^K \in \cdot \, \left|  \, \tau_0^K > t \right.\right) - \mu^K \right\|_{TV} \leq C(K) e^{ - \gamma(K) t},
	\end{align}
	for every initial distribution $\nu \in  \mathcal{P}( \hat{ \mathcal{X}}^K_+)$, where $\norm{\cdot}_{TV}$ denotes the total variation norm. We refer to the nice survey \cite{MV12} for the proof of these basic results and much more about quasi-stationary distributions. An interesting question, not treated here and left for future works, would be to find bounds on $C(K)$ and $\gamma(K)$. This was done by Chazottes, Collet and M\'el\'eard for the multitype birth-and-death processes they studied in \cite{CCM17}, which are similar to the multitype SIS model considered here in a constant  environment $\xi$ in the persistent case. They obtained $C(K)=\mathcal{O}(1)$ and  $\gamma(K)$ of order $1/\log(K)$, see their Theorem 3.1. Hence, in that context convergence to the QSD occurs on a time scale $\log(K)$, which is very small with respect to the time scale of extinction, exponential in $K$. We expect that this separation of time scales is also true in our context if $\Lambda>0$, \ref{ass:monotonesubhomo} holds and $x\mapsto Q(x)$ is constant. Indeed, in that case, on the one hand the typical extinction time starting from a large initial number of infectives (i.e. $K\norm{X^K_0}\gg 1$) grows at least as a power of $K$; one the other hand, the limiting PDMP $\hat{X}$ has good mixing properties, see Theorem \ref{th:BS19gen}, item \ref{convpers} (and more precisely \cite[Theorem 4.12]{BS19}), which makes us think that convergence to the QSD for the conditional marginals of $\hat{X}^K$ takes place on a much smaller time scale than extinction. If this is valid, then at least in that case, the QSD $\mu^K$ yields a good approximation of the marginal distributions of $\hat{X}^K$ for a large period of time, which motivates its study. 
	
	Here, we are interested in the behaviour of $\mu^K$ as $K$ goes to infinity. Once again, this behaviour strongly depends on the sign of $\Lambda$. Our main result is Theorem \ref{thm:qsdpers} below, which treats the persistent case. We also consider the non-persistent case in Proposition \ref{pro:qsdnonpersistent} below. Together with Proposition \ref{pro:OKp}, these three results extend Theorem \ref{th:intro2} from the introduction to our more general framework.
	
	Before presenting our results, let us recall that a few basic facts about the QSD (see e.g. \cite{MV12}). Under $\mb{P}^{\mu^K}$, the extinction time $\tau_0^K$ has an exponential distribution with a rate parameter $\lambda^K \in (0, +\infty)$, called the \emph{extinction rate}. In particular, $\mb{E}^{\mu^K}[\tau^K_0]=1/\lambda^K$. For all $t\in\mb{R}_+$, let us define the submarkovian kernel $\tilde{P}^K_t$ on $\mathcal{X}^K$ by
	\[
	\tilde{P}^K_t(\hat{x},A)=\mb{P}^{\hat{x}}\left(\hat{X}^K_t\in A, \tau^K_0>t\right)
	\]
	for all $\hat{x}\in \mathcal{X}^K$ and $A\subset \hat{X}^K$. The family $(\tilde{P}^K_t)_{t\in\mb{R}_+}$ constitutes a submarkovian semi-group on $\mathcal{X}^K$, called the \emph{killed semi-group} of $\hat{X}^K$. By definition of the QSD and the fact that $\mb{P}^{\mu^K}(\tau^K_0>t)=e^{-\lambda^Kt}$, we have $\mu^K\tilde{P}^K_t=e^{-\lambda^Kt}\mu^K$ for all $t\in\mb{R}_+$.
	
	\paragraph{The persistent case.} First, we consider the case where $\Lambda>0$. From Corollary \ref{theordertimeext} and the fact that $\mb{E}^{\mu^K}[\tau^K_0]=1/\lambda^K$, we immediately deduce the following result.	
	
	\begin{Pro}
		\label{pro:OKp}
		Assume that $\Lambda > 0$. Then $\lambda^K\underset{K\rightarrow +\infty}{\longrightarrow}0$, and more precisely
		\begin{align*}
			\limsup_{K\rightarrow +\infty}\frac{\log (\lambda^K)}{\log(K)}\leq -p^*.
		\end{align*}
		If, in addition, Assumption \ref{ass:monotonesubhomo} holds, then
		\begin{align*}
			\frac{\log (\lambda^K)}{\log(K)}\underset{K\rightarrow +\infty}{\longrightarrow}-p^*.
		\end{align*}
	\end{Pro}
	
	In order to study the asymptotic behaviour of the sequence $(\mu^K)_{K\geq d}$, we see these QSDs as probability distributions on the compact state space $\hat{\mathcal{X}}$, with the advantage that $\mathcal{P}(\hat{\mathcal{X}})$ is compact for the weak topology. In his PhD thesis \cite[Lemma 6.3]{S19}, Strickler proved that the convergence of $\lambda^K$ to $0$ entails that all weak limit points of $(\mu^K)_{K\geq d}$ in $\mathcal{P}(\hat{\mathcal{X}})$ are stationary probability distributions of the limiting PDMP $\hat{X}$. However, it remains to know whether part of the mass could escape to $\{0\}\times E$. Roughly speaking, in a large population $(K\gg 1)$, knowing that the disease is still present in the population at some large time $t$, is it likely that the proportion of infectives in the population is very low ? Or on the contrary, is it likely that a significant proportion of the population is still infected ? Intuitively, due to the tendency of the epidemic to invade the population from low levels of infectives in the persistent case, one is naturally inclined to choose the second option. In other words, it is quite intuitive to think that $\scr{L}\subset\mathcal{P}(\hat{\mathcal{X}}_+)$, where $\scr{L}$ denotes the set of weak limit points of $(\mu^{K})_{K\geq d}$. This is what happens for several models in constant environment, such as the birth-and-death processes studied in \cite{CCM16,CCMM20} and the supercritical stochastic SIS model (in discrete time) studied in \cite{SHJW20}, where the QSD converges to the unique non-zero stable equilibrium point of the limiting ODE. 
	
	Our main result, Theorem \ref{thm:qsdpers} below, states that this is indeed the case. In the following, we say that a collection (or a family) of elements of $\mathcal{P}(\hat{\mathcal{X}}_+)$ is \emph{persistent} if it is tight on $\mathcal{X}_+$. By Prohorov's theorem,  $\scr{L}\subset\mathcal{P}(\hat{\mathcal{X}}_+)$ is equivalent to saying that $(\mu^{K})_{K\geq d}$ is persistent.

	\begin{The}\label{thm:qsdpers}
		Assuming that $\Lambda>0$, the following hold.
		\begin{enumerate}
			\item \label{itemqsdpers} The sequence $(\mu^K)_{K\geq d}$ is persistent. Furthermore,  for all $p\in(0,p^*)$, there exists $C_4>0$ such that, for all $\mu\in\scr{L}$, 
			\[
			\int_{\hat{\mathcal{X}}_+}\norm{x}^{-p}\mu(\mr{d}x,\mr{d}\xi)\leq C_4.
			\]
			\item \label{lppers}Every element of $\scr{L}$ is a persistent stationary distribution of $\hat{X}$.
			\item \label{cvgqsd} If Assumption \ref{ass:monotonesubhomo} holds and $x\mapsto Q(x)$ is constant, then 
			\[
			\mu^K  \underset{K\rightarrow +\infty}{\Longrightarrow} \mu^*,
			\]
			where $\mu^*$ is the unique persistent stationary distribution for $\hat{X}$ given by Theorem \ref{th:BS19gen}.		
		\end{enumerate} 
	\end{The}
	
	The second item is a consequence of the first one and the fact that $\scr{L}$ is included in the set of invariant probability measures. As we just said, this fact follows from \cite[Lemma 6.3]{S19}, but for the sake of completeness we provide the proof in Section \ref{sec:prfqsdpers} (see Proposition \ref{pro:lambdapositifstationary}). As for item \ref{cvgqsd}, it follows from item \ref{lppers} and the unicity of the persistent stationary distribution when \ref{ass:monotonesubhomo} holds and $x\mapsto Q(x)$ is constant, given by Theorem \ref{th:BS19gen}. Note that in the case of a constant environment, this says that $\mu^K$ converges to the Dirac mass on the endemic equilibrium, as in \cite{SHJW20,CCM16,CCMM20}.
	
	The difficult part of Theorem \ref{th:BS19gen} is item \ref{itemqsdpers}. Although we argued that the result is somewhat intuitive, the mathematical proof requires a fine analysis of the microscopical behaviour of the process near the extinction set. We stress that \emph{the convergence of $\hat{X}^K$ to a persistent process $\hat{X}$ is not enough.} Let us explain why, on a simple example. Consider the classical supercritical stochastic SIS model, corresponding to $d=|E|=1$, $X^K=N^K/K$ where $N^K$ is a birth-and-death process with state space $\intbk{0,K}$, birth rate $bn(1-n/K)$ and death rate $\tilde{d}n$ in state $n$, with $b>\tilde{d}>0$.
	Consider also the slightly modified Markov process $\tilde{X}^K=\tilde{N}^K/K$ where $\tilde{N}^K$ is another birth-and-death process, with the same transition rates as $N^K$ except starting from state $1$, where we set the birth rate and the death rate equal to some very small $\varepsilon_K\in(0,1]$. The value of $\varepsilon_K$ will be fixed later. We know that if $X^K_0\underset{K\rightarrow +\infty}{\longrightarrow} x_0$, then $X^K$ converges in distribution to the solution of the Cauchy problem
	\begin{align}\label{cauchySIS}
		\dot{x}=x(b-\tilde{d}-bx),\quad x(0)=x_0,
	\end{align}  in the Skorokhod space $\mathcal{D}([0,T],\mb{R})$, for all $T>0$. It is not hard to see that the same goes for $\tilde{X}^K$ : if $\tilde{X}^K_0\rightarrow 0$, then $\tilde{X}^K\rightarrow 0$ in distribution on $\mathcal{D}([0,T],\mb{R})$; while if $\tilde{X}^K_0\rightarrow x_0>0$, $\tilde{X}^K$ has the same distribution as $X^K$ until the hitting time of $1/K$, which tends to $+\infty$, hence $\tilde{X}^K$ converges in distribution in $\mathcal{D}([0,T],\mb{R})$ to the solution of \ref{cauchySIS}. However the process $\tilde{X}^K$ has the particularity that the state $1/K$ is almost absorbing, in the sense that the total jump rate from $1/K$ is very small, equal to $2\varepsilon_K$. Due to the irreducibility of $K^{-1}\intbk{1,K}$ and the accessibility of $0$ for the process $\tilde{X}^K$, it admits a unique QSD $\tilde{\mu}^K$. We claim that by choosing $\varepsilon_K$ small enough, we can make $\tilde{\mu}^K(\{1/K\})$ arbitrarily close to $1$. The QSD $\tilde{\mu}^K$ satisfies the eigenequation $\tilde{\mu}^K\tilde{L}^K=-\tilde{\lambda}^K\tilde{\mu}^{K}$ for some $\tilde{\lambda}^K>0$, where $\tilde{L}^K$ is the infinitesimal generator of the killed semi-group of $\tilde{X}^K$, i.e. the transition rate matrix of $\tilde{X}^K$ restricted to the states $K^{-1}\intbk{1,K}$. Setting $m^K_i=\tilde{\mu}^K(\{i/K\})$ and $\ell^{K}_{i,j}=\tilde{L}^K(i/K,j/K)$, $i,j\in\intbk{1,K}$, we have
	\[
	\tilde{\mu}^K\tilde{L}^K(\{1/K\})=-m^K_1(2\varepsilon_K)+m^K_2(2\tilde{d})=-\tilde{\lambda}^Km^K_1
	\] 
	and for all $i\in \intbk{2,K-1}$, 
	\[
	\tilde{\mu}^K\tilde{L}^K(\{i/K\})=m^K_{i-1}\ell^K_{i-1,i}-m^K_i|\ell^K_{i,i}|+m^K_{i+1}\ell^K_{i+1,i}=-\tilde{\lambda}^Km^K_i.
	\]
	Thus, we have
	\[
	m^K_2\leq m^K_1(\varepsilon_K/\tilde{d})\quad \text{and} \quad \forall i\in \intbk{2,K-1},\  m^K_{i+1}\leq m^K_i|\ell^K_{i,i}|/\ell^K_{i+1,i}. 
	\]
	Since the rates $|\ell^K_{i,i}|, \ell^K_{i+1,i}$, $i\in \intbk{2,K-1}$, do not depend on $\varepsilon_K$, these inequalities entail that we can fix the value of $\varepsilon_K$ small enough so that  $\sum_{i=2}^{K}m^K_i\leq C_K\varepsilon_K m^K_1$ for some constant $C_K$ which does not depend on $\varepsilon_K$. Given that $\sum_{i=1}^Km^K_i=1$, we obtain $m^K_1\geq 1/(1+C_K\varepsilon_K)$. Hence, if we fix $\varepsilon_K=1\wedge(KC_K)^{-1}$, we obtain $m^K_1\geq 1-1/(K+1)$. With such a choice, we see that the QSD $\tilde{\mu}^K$ converges weakly to $\delta_0$ as $K\rightarrow +\infty$, even though $\tilde{X}^K$ converges to the solution of an ODE repelled by $0$.
	
	\medskip
	Let us go back to the general framework. As we just explained, the proof of the first item of Theorem \ref{thm:qsdpers} requires fine control of the microscopical behaviour of the chain $\hat{X}^K$ near the extinction set $\{0\}\times E$. This is what the following proposition provides.
	
	\begin{Pro}\label{lemphitheta:1}
		Assume that $\Lambda>0$ and let $p\in (0,p^*)$. Let $T,a>0$  be given by Proposition \ref{prolyapXK:1} applied to $-p$, and for all $K\geq d$, let $\varphi^K_1:\hat{\mathcal{X}}^K_+\to \mb{R}_+^*$ be defined by $\varphi^K_1(x,\xi)=\norm{x}^{-p}\wedge(a/K)^{-p}$. There exist $\theta\in (0,1)$ and $C_6>0$ such that :
		
		\begin{enumerate}[label=\roman*)]
			\item \label{itemphitheta1:1}for all $K$ large enough,
			\begin{align}\label{phitheta1:1}
				\tilde{P}^K_T\varphi^K_1\leq \theta \varphi^K_1 + C_6\ ;
			\end{align}
			\item \label{tightqsd:1} $\limsup_{K\rightarrow +\infty}\mu^K\varphi_1^K\leq C_6/(1-\theta)$.
		\end{enumerate}
	\end{Pro}
	
	In order to obtain the global Lyapunov condition \eqref{phitheta1:1} for $\varphi^K_1$ with respect to the killed kernel $\tilde{P}^K_T$, there are two main ingredients. Obviously, the first one is given by Proposition \ref{prolyapXK:1} which yields the Lyapunov condition on the set $\{(x,\xi)\in\hat{X}^K : \norm{x}\geq a/K\}$. The second main ingredient is given by Lemma \ref{lemboundary}, Equation \ref{boundext}, which yields a lower bound, independent of $K$, on the probability of extinction of  $\hat{X}^K$ in constant time starting from $\{(x,\xi)\in\hat{X}^K : 1\leq \norm{x}< a/K\}$. This notably prevents that $\hat{X}^K$ stays too long in that region, which was the problem for the "pathological" process $\tilde{X}^K$ constructed above. 
	We stress that in order to be able to glue these two ingredients and get a constant $\theta$ independent of $K$, it proved crucial that the Lyapunov condition given by Proposition \ref{prolyapXK:1} be able to cover a region which is at distance of order only $1/K$ of the extinction set. 
	Finally, item \ref{tightqsd:1} easily follows from the integration of \eqref{phitheta1} with respect to $\mu^K$ and the fact that $\lambda^K$ converges to $0$.
	
	Let us mention that the idea of using a global Lyapunov condition with respect to the killed semi-group in order to control moments of the QSD was inspired to us by the reading of the works of Champagnat and Villemonais, see e.g. \cite[Lemma 9.6]{CV17gen}. Here, we do not need a reverse Lyapunov condition because we already have an upper bound on $\lambda^K$.

	\paragraph{The non-persistent case.}
	
	To complete the picture, we give the following result for the non-persistent case, which is a direct consequence of \cite[Theorem 6.2]{S19}. 
	\begin{Pro}
		\label{pro:qsdnonpersistent}
		Assume that $\Lambda < 0$ and that Assumption \ref{ass:monotonesubhomo}holds. Then
		\[
		\mu^K( \cdot \times E) \underset{K\rightarrow +\infty}{\Longrightarrow} \delta_0.
		\]
	\end{Pro}

	\section{Perspectives}
	
	A first important question that we have left for future work concerns the rate of convergence to the QSD of the marginal distributions of $\hat X^K$ conditional on non-extinction, in the persistent case. We believe that the environmental variation should not prevent a rate $\gamma(K)$ of order $1/\log(K)$ from being obtained in \eqref{cvgtoqsd}, as for the multitype birth-and-death processes in constant environment studied in \cite{CCM17}.
	
	Second, although we have focused on SIS epidemic models in this paper, it seems clear that several of the techniques that we have developed in this paper can be used or adapted for other classes of birth-and-death processes in random environment. We are in particular thinking of Lotka-Volterra competitive or prey-predator models in random environment. For these processes, when we are interested in coexistence of all the species, extinction corresponds to at least one of the species being extinct. Conditions for coexistence of the limiting PDMP are given by signs of average growth rates of the species when rare, see e.g. \cite{BL16}. It is also possible to build Lyapunov functions of the form $(x,\xi)\mapsto |x_i|^{-p}$ for the PDMP near the extinction set of species $i$ and thus, probably, for the associated birth-and-death process. One important issue in the birth-and-death Lotka-Volterra models, that was missing in the epidemic model we investigated, is the non finiteness of the state space : the population can grow to infinity and is no more bounded by $K$, which is in that case a scaling parameter representing the typical size of the population at equilibrium. However, if the limiting PDMP as an attracting compact set, as it the cases for most models, we believe that the behaviour at infinity can be handled with a bit of work. Moreover, using techniques borrowed to Champagnat and Villemonais \cite{CV17gen}, it is certainly possible to prove the existence and uniqueness of a quasi-stationary distribution for each $K$. We leave this extension for future research.	
	
	Finally, it would be interesting to be even more precise and general in our asymptotics. To begin with, the critical case $\Lambda=0$ remains open. In addition, for the persistent case, when $p^*<\infty$ we have obtained a logarithmic equivalent of the extinction time but not an equivalent. When $p^*=+\infty$, Remark \ref{rem:subexp} showed that it can lead both to subexponential and exponential order for the extinction time, hence it would be nice to determine under what conditions the extinction time is exponential in $K$. Finally, in the non-persistent case $\Lambda<0$, we conjecture that the typical extinction time from fixed proportions of infectives is of order $\log(K)/|\Lambda|$, as suggested by \ref{PDMPlinLambda} and the fact that $\hat{X}^K$ is likely to die in a time of order one as soon as it gets of order $1/K$.	

	\section{Proofs}
	\label{sec:proofs}
	
	\subsection{Proof of Proposition \ref{lemcoupl}}\label{prflemcoupl}
	
	Let us fix $K$ large enough so that for all $i\in\intbk{1,d}$, $K_i(K)/K\geq \underline{\alpha}/2$. Let $T>0$, $x\in\mathcal{X}^K$, $y\in\mathcal{X}$ and $\xi\in E$, and set $\hat{x}=(x,\xi)$, $\hat{y}=(y,\xi)$. Finally, let $\varepsilon>\norm{x-y}$, and introduce the following stopping times : 
	\[
	\sigma^{K,\hat{x},\hat{y}}_{\varepsilon}:=\inf\{t\geq 0 : \norm{X^{K,\hat{y}}_t-X^{\hat{x}}_t}>\varepsilon\},\quad \eta^{K,\hat{x},\hat{y}}:=\inf\{t\geq 0 : \Xi^{K,\hat{y}}_t\neq\Xi^{\hat{x}}_t\}.
	\]
	
	We claim that for some constants $C_0,C'_0>0$, which do not depend on $K,x,y,\varepsilon$ nor $T$, the following hold :

	\begin{gather}
		\mb{P}\left[\sigma^{K,\hat{x},\hat{y}}_{\varepsilon}\leq T\wedge \eta^{K,\hat{x},\hat{y}} \right]\leq 2d\exp\left(-\frac{K\delta(\varepsilon,T,x,y)}{C_0}\left(\frac{\delta(\varepsilon,T,x,y)}{C'_0(Te^{C_FT}\norm{x}+\varepsilon)}\wedge 1\right)\right) \label{ineqcoupl1}\\
		\mb{P}\left[\eta^{K,\hat{x},\hat{y}}\leq T \wedge \sigma^{K,\hat{x},\hat{y}}_{\varepsilon}\right]\leq T|E|\sup_{\xi_1\neq \xi_2,\ \norm{z-z'}\leq \varepsilon}{\left|q(z,\xi_1,\xi_2)-q(z',\xi_1,\xi_2)\right|} \label{ineqcoupl2}.	
	\end{gather}
	
	Let us start by the proof of \eqref{ineqcoupl1}. In the following, we simplify the notations by writing $\sigma^K_{\varepsilon}=\sigma^{K,\hat{x},\hat{y}}_{\varepsilon}$ and $\eta^K=\eta^{K,\hat{x},\hat{y}}$. Combining \eqref{eqXKbis} and \eqref{eqX}, we obtain that $\mb{P}$-almost surely for all $t\leq \eta^K$,  
	\[
	\norm{X^{K,\hat{y}}_t-X^{\hat{x}}_t}\leq \norm{y-x}+ C_F\int_0^t\norm{X^{K,\hat{y}}_s-X^{\hat{x}}_s}\mr{d}s+ \norm{M^{K,\hat{y}}_t}
	\]
	where 
	\[
	M^{K,\hat{y}}_t=\int_{(0,t]\times \mb{R}_+\times \llbracket 1,d\rrbracket\times \left\{-1,1\right\}}\mb{1}_{\left\{u\leq K_i \beta_{he_i}(\hat{X}^{K,\hat{y}}_{s-})\right\}}\frac{he_i}{K_i}\tilde{\scr{N}}_X(\mr{d}s,\mr{d}u,\mr{d}i,\mr{d}h)
	\]
	Thus Grönwall's lemma yields, almost surely for all $t\leq \eta^K$,
	\[
	\sup_{0\leq s \leq t}\norm{X^{K,\hat{y}}_s-X^{\hat{x}}_s}\leq e^{C_F t}\left(\norm{y-x}+\sup_{0\leq s \leq t}\norm{M^{K,\hat{y}}_{s}}\right).
	\]
	It follows that 
	\begin{align}\label{sigmaepsilon}
	\left\{\sigma^K_{\varepsilon}\leq T \wedge \eta^K\right\}&\subset \left\{\sup_{0\leq s \leq T \wedge \eta^K\wedge \sigma^K_{\varepsilon}}\norm{X^{K,\hat{y}}_s-X^{\hat{x}}_s}\geq \varepsilon\right\}\nonumber\\&\subset \left\{ \sup_{0\leq s \leq T}\norm{M^{K,\hat{y}}_{s\wedge \eta^K\wedge \sigma^K_{\varepsilon}}}\geq \delta(\varepsilon,T,x,y) \right\}
	\end{align}
	where $\delta(\varepsilon,T,x,y):=(\varepsilon e^{-C_FT}-\norm{x-y})_+$. Let $M^{K,\hat{y},i}$, $i\in \llbracket 1,d\rrbracket$, denote the coordinates of $M^{K,\hat{y}}$. We have
	\[
	M^{K,\hat{y},i}_{t\wedge\eta^K\wedge \sigma^K_{\varepsilon}}=\int_{(0,t]\times \mb{R}_+\times \left\{-1,1\right\}}G^{K,\hat{y},i}_{s,u}\tilde{\scr{N}_X}(\mr{d}s,\mr{d}u,\left\{i\right\},\mr{d}h), 
	\]
	with
	\begin{align}\label{G}
	G^{K,\hat{y},i}_{s,u}=\mb{1}_{\left\{s\leq \eta^K\leq \sigma^K_{\varepsilon}\right\}}\mb{1}_{\left\{u\leq K_i\beta_{he_i}\left(\hat{X}^{K,\hat{y}}_{s-}\right)\right\}}\frac{h}{K_i}.
	\end{align}
	Let $C_{\beta}=\sup_{i\in\llbracket 1,d\rrbracket, h\in\left\{-1,1\right\},\zeta\in E}\norm{\beta^{\zeta}_{he_i}}_{\mr{Lip}}$. For all $s\leq T\wedge\eta^K$, $i\in\llbracket 1,d\rrbracket$ and $h\in\left\{-1,1\right\}$ we have, using that $\beta^\zeta_{he_i}(0)=0$ for all $\zeta\in E$,
	\begin{align}\label{beta}
	\beta_{he_i}(\hat{X}^{K,\hat{y}}_{s-})\leq C_{\beta}(\norm{X^{\hat{x}}_{s-}}+\varepsilon)\leq C_{\beta}(\norm{x}e^{C_FT}+\varepsilon).
	\end{align}
	Using \eqref{G}, \eqref{beta} and applying Lemma \ref{lemChernoff}, we obtain
	\begin{align*}
	&\mb{P}\left(\sup_{0\leq s \leq T}\norm{M^{K,\hat{y}}_{s\wedge \eta^K\wedge \sigma^K_{\varepsilon}}}\geq \delta(\varepsilon,T,x,y)\right) \\
	&\leq \sum_{i\in\llbracket 1,d \rrbracket}\mb{P}\left(\sup_{0\leq s \leq T}\nm{M^{K,\hat{y},i}_{s\wedge \eta^K\wedge \sigma^K_{\varepsilon}}}\geq \delta(\varepsilon,T,x,y)/d\right) \\
	&\leq 2d\exp\left(-\frac{K \delta(\varepsilon,T,x,y) }{4d\underline{\alpha}}\left(\frac{\delta(\varepsilon,T,x,y) }{4dT(C_{\beta}e^{C_FT}\norm{x}+\varepsilon)}\wedge\log(2)\right)\right)	
	\end{align*}
    which yields \eqref{ineqcoupl1} with $C_0=4d\underline{\alpha}$ and $C'_0=4d(C_{\beta}\vee 1)/\log(2)$. 
 
	
	Now, we prove \eqref{ineqcoupl2}. Let us introduce the counting process $(D^K_t)_{t\geq 0}$ defined by
	\[
	D^K_t=\int_{(0,t]\times\mb{R}_+\times E}\mb{1}_{\left\{\xi'\neq \Xi^{\hat{x}}_{s-}\right\}}\left|\mb{1}_{\left\{u\leq q(X^{K,\hat{y}}_{s-},\Xi^{\hat{x}}_{s-},\xi')\right\}}-\mb{1}_{\left\{u\leq q(X^{\hat{x}}_{s-},\Xi^{\hat{x}}_{s-},\xi')\right\}}\right|\scr{N}_{\Xi}\left(\mr{d}s,\mr{d}u,\mr{d}\xi'\right).
	\]
	Recalling Equations \eqref{eqIK} and \eqref{eqI}, and given that $\Xi^{K,\hat{y}}$ and $\Xi^{\hat{x}}$ coincide before $\eta^K$, we see that $D^K$ has the useful property that
	\[
	\left\{\eta^K\leq T\wedge\sigma^K_{\varepsilon}\right\}=\left\{D^K_{T\wedge\eta^K\wedge\sigma^K_{\varepsilon}}= 1\right\}.
	\]
	It follows that 
	\begin{align*}
	\mb{P}\left[\eta^K\leq T\wedge\sigma^K_{\varepsilon}\right]&\leq \mb{E}\left((D^K_{T\wedge\sigma^K_{\varepsilon}})^2\right)\\
	&=\mb{E}\left[\int_0^{T\wedge \eta^K\wedge \sigma^K_{\varepsilon}} \mb{1}_{\left\{s\leq \sigma^K_{\varepsilon}\right\}}\sum_{\xi'\neq \Xi^{\hat{x}}_s }\left|q(X^{K,\hat{y}}_{s},\Xi^{\hat{x}}_{s},\xi')-q(X^{\hat{x}}_{s},\Xi^{\hat{x}}_s,\xi')\right|\mr{d}s\right] \\
	&\leq T|E|\sup_{\norm{z-z'}\leq \varepsilon,\,\xi_1\neq \xi_2,}\left|q(z,\xi_1,\xi_2)-q(z',\xi_1,\xi_2)\right|,	
	\end{align*}
	which ends the proof of \eqref{ineqcoupl2}.
	Now we can conclude : 
	\begin{align*}
	\hspace{-2pt}\mb{P}\left[\sup_{0\leq t \leq T}\mb{\hat{d}}\left(\hat{X}^{K,\hat{y}}_t,\hat{X}^{\hat{x}}_t\right)>\varepsilon\right]&\leq \mb{P}\left[\sigma^{K}_{\varepsilon}\wedge \eta^{K}\leq T\right]\\&\leq \mb{P}\left[\sigma^{K}_{\varepsilon}\leq T\wedge\eta^{K}\right]+\mb{P}\left[\eta^{K}\leq T\wedge\sigma^{K}_{\varepsilon}\right] \\&\leq  2d\exp\left(-\frac{K\delta(\varepsilon,T,x,y)}{C_0}\hspace{-2pt}\left(\hspace{-2pt}\frac{\delta(\varepsilon,T,x,y)}{C'_0(Te^{C_FT}\norm{x}+\varepsilon)}\wedge 1\right)\hspace{-2pt}\right)\nonumber \\  &\ +T|E|\sup_{\xi_1\neq \xi_2,\ \norm{z-z'}\leq \varepsilon}\left|q(z,\xi_1,\xi_2)-q(z',\xi_1,\xi_2)\right|.
	\end{align*}\hfill $\square$

	\subsection{Proof of Lemma \ref{lemtauq}}\label{prflemtauq}
	Let $K\geq d$ and $\hat{x}=(x,\xi)\in\hat{\mathcal{X}}^K_+$. We work under $\mb{P}^{\hat{x}}$. Equation \eqref{eqXKbis} yields
	\[
	\norm{X^K_{t}}=\langle \mb{1},X^K_t\rangle= \norm{x}+ \int_0^{t} \langle \mb{1},F(\hat{X}^K_{s})\rangle\mr{d}s+\langle \mb{1},M^K_t \rangle.
	\]
	Since for all $s\geq 0$,  $-C_F\norm{X^K_s}\leq\langle \mb{1},F(\hat{X}^K_{s})\rangle\leq C_F\norm{X^K_s}$, applying generalized Ito's formula to $t\mapsto e^{C_Ft}\norm{X^K_t}$ and $t\mapsto e^{-C_Ft}\norm{X^K_t}$ yields
	\begin{align}\label{ineqX}
		e^{C_Ft}\norm{X^K_t} \geq \norm{x}-\underline{Z}^K_t\quad \text{and} \quad e^{-C_Ft}\norm{X^K_t}\leq \norm{x} +\overline{Z}^K_t
	\end{align}
	where
	\begin{align*}
		\underline{Z}^K_t &:=-\sum_{i=1}^d\int_{0+}^te^{C_Fs}\mr{d}M^K_i(s)\nonumber\\&=-\int_{(0,t]\times \mb{R}_+\times\intbk{1,d}\times\left\{-1,1\right\}}\mb{1}_{\left\{u\leq K_i\beta_{he_i}(\hat{X}^K_{s-})\right\}}\frac{h}{K_i}e^{C_Fs}\tilde{\scr{N}}_{X}(\mr{d}s,\mr{d}u,\mr{d}i,\mr{d}h), \\
		\overline{Z}^K_t &:=\int_{0+}^te^{-C_Fs}\mr{d}M^K_i(s)\nonumber\\&=\int_{(0,t]\times \mb{R}_+\times\intbk{1,d}\times\left\{-1,1\right\}}\mb{1}_{\left\{u\leq K_i\beta_{he_i}(\hat{X}^K_{s-})\right\}}\frac{h}{K_i}e^{-C_Fs}\tilde{\scr{N}}_{X}(\mr{d}s,\mr{d}u,\mr{d}i,\mr{d}h).
	\end{align*}
	Using \eqref{ineqX}, we get
	\begin{align*}
		\left\{\overline{\tau}^{K}_{M\norm{x}\wedge\underline{\tau}^K_{m\norm{x}}}\leq T\right\}& \subset \left\{\underline{\tau}^K_{m\norm{x}}\leq T\wedge \overline{\tau}^{K}_{M\norm{x}}\right\}\cup
		\left\{\overline{\tau}^{K}_{M\norm{x}}\leq T \right\} \\
		& \subset \left\{\sup_{0\leq t \leq T}\underline{Z}^K_{t\wedge \overline{\tau}^{K}_{M\norm{x}}}\geq \norm{x}/2\right\}\cup\left\{\sup_{0\leq t \leq T}\overline{Z}^K_{t\wedge \overline{\tau}^{K}_{M\norm{x}}}\geq \norm{x}\right\}.
	\end{align*}
	
	We can bound the probability of each event of this union using Lemma \ref{lemChernoff}. Set $C_{\beta}=\sup_{i\in\intbk{1,d}, h\in \left\{-1,1\right\}}\norm{\beta^{\xi}_{he_i}}_{\mr{Lip}}$. We suppose that $K$ is large enough so that $K_i(K)/K\geq \underline{\alpha}/2$. Since 
	\[
	\underline{Z}^K_{t\wedge \overline{\tau}^{K}_{M\norm{x}}}=\int_{(0,t]\times \mb{R_+}\times\intbk{1,d}\times \left\{-1,1\right\}}\underline{G}_{s,u,i,h}\tilde{\scr{N}}_X(\mr{d}s,\mr{d}u,\mr{d}i,\mr{d}h)
	\]
	with 
	\begin{align*}
	\left|\underline{G}_{s,u,i,h}\right|=\left|-\mb{1}_{\left\{s\leq \overline{\tau}^{K}_{M\norm{x}}\right\}}\mb{1}_{\left\{u\leq K_i\beta_{he_i}(\hat{X}^K_s)\right\}}\frac{h}{K_i}e^{C_Fs}\right|\leq \frac{2e^{C_FT}}{\underline{\alpha}K}\mb{1}_{\left\{u\leq  KC_{\beta}M \norm{x}\right\}}
	\end{align*}
	almost surely for all $(s,u,i,h)\in [0,T] \times \mb{R}_+\times\intbk{1,d}\times\left\{-1,1\right\}$, Lemma \ref{lemChernoff} yields
	\[
	\mb{P}^{\hat{x}}\left(\sup_{0\leq t \leq T}\underline{Z}^K_{t\wedge \overline{\tau}^{K}_{M\norm{x}}}\geq \norm{x}/2\right)\leq \exp\left(-\underline{C}'K\norm{x}\right)
	\]
	where 
	\[
	\underline{C}'=\frac{\underline{\alpha}^2}{128dC_{\beta}Me^{2C_FT}T}\wedge\frac{\underline{\alpha}\log(2)}{8e^{C_FT}}>0.
	\]		
	Similarly, we obtain
	\[
	\mb{P}^{\hat{x}}\left(\sup_{0\leq t \leq T}\overline{Z}^K_{t\wedge \overline{\tau}^{K}_{M\norm{x}}}\geq \norm{x}\right)\leq \exp\left(-\overline{C}'K\norm{x}\right)
	\]
	for some $\overline{C}'>0$.	Hence, 
	\[
	\mb{P}^{\hat{x}}(\overline{\tau}^{K}_{M\norm{x}}\wedge\underline{\tau}^{K}_{m\norm{x}}\leq T)\leq e^{-\underline{C'} K\norm{x}}+e^{-\overline{C'} K\norm{x}}\leq 2e^{-C' K\norm{x}}
	\]
	with $C'=\underline{C'}\wedge\overline{C'}>0$.
	\hfill $\square$
	
	\subsection{Proof of Lemma \ref{lemboundary}}\label{prflemboundary}
	
	Let $T,a>0$. We first prove the existence of a constant $c>0$ satisfying \eqref{boundext}. By Assumption \ref{ass:Metzler}, we have 
	\[
	C_b:=\sup_{\substack{i\in \intbk{1,d}\\\xi \in E}}\norm{\beta^{\xi}_{e_i}}_{\mr{Lip}}<\infty\quad \text{and}\quad \underline{d}:=\min_{\substack{i\in \intbk{1,d}\\ \xi \in E}}d^{\xi}_i(0)>0.		\]
	By continuity of the $d^\xi_i$, there exists $K_0\in\mb{N}^*$ such that for all $K\geq K_0$,
	\[
	\min_{i\in \intbk{1,d}}\inf_{\substack{(x,\xi)\in \hat{\mathcal{X}} \\ \norm{x}\leq a/K}}d^{\xi}_i(x)\geq \underline{d}/2
	\]
	Let $K\geq K_0$ and $\hat{x}=(x,\xi) \in\hat{\mathcal{X}}^K$ such that $\norm{x}\leq a/K$. We can get a very rough lower bound on the probability of extinction before time $T$ starting from $\hat{x}$ by giving a condition on $\scr{N}_X$ which ensures that :
	\begin{enumerate}[label=(\alph*)]
		\item \label{itema} there is no infection on $[0,T]$; 
		\item \label{itemb} there is at least $K_ix_i$ recoveries in each group $i$.
	\end{enumerate}
	Recalling \eqref{eqXK}, we see that \ref{itema} holds $\mb{P}^{\hat{x}}$-a.s. on the event 
	\[
	I:=\left\{\scr{N}_X\left([0,T]\times [0,C_b a]\times \intbk{1,d}\times \left\{1\right\}\right)=0\right\},
	\]
	since $K_i\beta^\xi_{he_i}(y)\leq KC_b\norm{y}$ for all $i\in\intbk{1,d}$ and $(y,\zeta)\in\mathcal{X}^K$. Moreover, \ref{itemb} holds $\mb{P}^{\hat{x}}$-a.s. on the event
	\[
	R:=\bigcap_{i=1}^d\left\{\scr{N}_X\left([0,T]\times [0,\underline{d}]\times \left\{i\right\}\times \left\{-1\right\}\right)\geq a\right\}.	
	\]
	Indeed : (i) $\norm{x}\leq a/K$ entails that in each group $i$, the initial number of infected individuals $K_ix_i$ is less than or equal to $a$; (ii) each atom of $\scr{N}_X$ in $[0,T]\times [0,\underline{d}]\times \left\{i\right\}\times \left\{-1\right\}$ causes a new recovery in group $i$ as long are there are still infected individuals in that group, since $K_i\beta_{-e_i}^{\xi}(y)=K_iy_id_i^{\xi}(y)\geq \underline{d}$ for all $y\in\mathcal{X}^K_+$ such that $\norm{y}\leq a/K$. It follows that $\left\{I\cap R\right\}\subset \left\{\tau_0^K\leq T\right\}$ $\mb{P}^{\hat{x}}$-almost surely. Hence, \eqref{boundext} is satisfied with $c:=\mb{P}^{\hat{x}}(I\cap R)$, which is positive and does not depend on $\hat{x}$.
	
	\medskip		
	By a similar argument, we can prove the existence of a constant $c'>0$ satisfying \eqref{boundhita}. By Assumption \ref{ass:Metzler}, we have
	\[\underline{b}:=\min_{\substack{j\in \intbk{1,d} \\\xi \in E}}\,\max_{i\in \intbk{1,d}} \partial_jb^{\xi}_i(0)>0,
	\]
	and for all $i\in\intbk{1,d}$ and $\xi\in E$, 
	\begin{align}\label{Taylorinf}
		\beta_{e_i}(x,\xi)=(1-x_i)b_i^{\xi}(x)=\sum_{j=1}^d\partial_jb_i^\xi(0)x_j+o(\norm{x}).
	\end{align}
	Fix $j\in\intbk{1,d}$, $\xi\in E$, and let $i=i_{j,\xi}\in\intbk{1,d}$ be such that $\partial_jb_i^\xi(0)\geq \underline{b}$. Expansion \eqref{Taylorinf} entails that we may choose $K_0$ large enough so that for all $K\geq K_0$, the rate of infection in group $i$ satisfies 
	\begin{align}\label{minrateinf}
		\inf_{x\in\mathcal{X}^K,\, \norm{x}\leq a/K,\, x_j>0}	K_i\beta_{e_i}(x,\xi)\geq K_i(\underline{b}x_j+o(1/K))\geq (\underline{\alpha}\underline{b})/2,
	\end{align}
	since $x_j>0$ implies $x_j\geq 1/K$.  
	Let us introduce the events
	\begin{gather*}
		I'_{j,\xi}=\left\{\scr{N}_X\left([0,T]\times [0,(\underline{\alpha}\underline{b})/2]\times \{i_{j,\xi}\}\times \left\{1\right\}\right)\geq a\right\},\\
		R'=\left\{\scr{N}_X\left([0,T]\times [0,\overline{d}]\times \left\{1,d\right\}\times \left\{-1\right\}\right)=0\right\},\\
		E'=\left\{\scr{N}_{\Xi}\left([0,T]\times [0,\overline{q}]\times E\right)=0\right\}
	\end{gather*}
	where $\overline{d}=\max_{i\in\intbk{1,d}}\norm{d_i}_{\infty}$ and
	$\overline{q}=\sup_{x\in\mathcal{X},\, \zeta\neq \zeta'} q(x,\zeta,\zeta')$. Let $K\geq K_0$ and $x\in\mathcal{X}$ such that $\norm{x}\leq a/K$ and $x_j>0$. We work under $\mb{P}^{(x,\xi)}$. The event $E'$ implies that the environment $\Xi$ remains equal to $\xi$ on $[0,T]$, and $R'$ implies that there is no recovery on $[0,T]$. What's more, inequality \eqref{minrateinf} entails that each atom of $\scr{N}_X$ in $[0,T]\times [0,(\underline{\alpha}\underline{b})/2]\times \{i\}\times \left\{1\right\}$ arriving at a time $t$ such that $\Xi_t=\xi$, $\norm{X^K_t}\leq a/K$ and $X^K_j(t)>0$ causes a new infection in group $i$. Hence, we see that $I'_{j,\xi}\cap R'\cap E'$ implies $\left\{\bar{\tau}^K_{a/K}\leq T\right\}$. This entails that for all $K$ large enough, for all $x\in\mathcal{X}^K$ such that $\norm{x}\leq a/K$ and $x_j>0$, we have
	\[\mb{P}^{(x,\xi)}\left(\bar{\tau}^K_{a/K}\leq T\right)\geq c'_{j,\xi},
	\]
	where $c'_{j,\xi}=\mb{P}^{(x,\xi)}\left(I'_{j,\xi}\cap R'\cap E'\right)>0$, which does not depend on $x$.  
	Then, $c'=\min c'_{j,\xi}$ satisfies $\eqref{boundhita}$ for all $K$ large enough, which concludes the proof.
	
	\hfill $\square$
	
	\subsection{Proof of Theorem \ref{plyap}}
	
	Let $p\in\mb{R}$. Let $\overline{f}_p:\mb{R}_+\to\mb{R}$ be defined by 
	\[
	\overline{f}_p(t)=\log\sup_{\hat{\theta}\in\Delta\times E}\mb{E}\left[\exp\left(p\int_0^t G_0(\hat{\Theta}^{\hat{\theta}}_s)\mr{d}s\right)\right].
	\]
	The Markov property of the family $\hat{\Theta}^{\hat{\theta}}$ entails that $\overline{f}_p$ is a subadditive function, i.e., $\overline{f}_p(s+t)\leq \overline{f}_p(s)+\overline{f}_p(t)$ for all $s,t\geq 0$. Since $\overline{f}_p$ is Borel measurable and $\overline{f}_p(t)/t\in\left[-|p|\|G_0\|_{\infty},|p|\|G_0\|_{\infty}\right]$ for all $t>0$, the fundamental theorem of subadditive functions \cite[Theorem 7.6.1]{HilPhi} entails that 
	\[
	\frac{\overline{f}_p(t)}{t}\underset{t\rightarrow +\infty}{\longrightarrow} \lambda_p
	\]
	for some $\lambda_p\in\left[-|p|\|G_0\|_{\infty},|p|\|G_0\|_{\infty}\right]$. For all $t>0$ and $\hat{\theta}\in\Delta\times E$ the function $p\mapsto \overline{f}_p(t)$ is a supremum of cumulant generating functions, hence it is convex, and consequently, $p\mapsto \lambda_p$ is convex, being a limit of convex functions. 

	Similarly, the function $\underline{f}_p\colon \mb{R}_+\to\mb{R}$ defined by
	\[
	\underline{f}_p(t)=\log\inf_{\hat{\theta}\in\Delta\times E}\mb{E}\left[\exp\left(p\int_0^t G_0(\hat{\Theta}^{\hat{\theta}}_s)\mr{d}s\right)\right]
	\]
	is Borel-measurable, superadditive, and
	\[
	\frac{\underline{f}_p(t)}{t}\underset{t\rightarrow +\infty}{\longrightarrow} \nu_p
	\]
	for some $\nu_p\in\left[-|p|\|G_0\|_{\infty},|p|\|G_0\|_{\infty}\right]$. 
	
	We must show that $\nu_p=\lambda_p$. In order to do this, we are going to compare the trajectories of $\hat{\Theta}$ starting from different initial conditions. For all $\xi\in E$, let $\varphi^{\xi}:\mb{R}_+\times \Delta \to \Delta$ be the projection of the linear flow induced by $A^\xi$ on the simplex, i.e.
	\[
	\varphi^{\xi}(t,\theta)=\frac{e^{tA^{\xi}}\theta}{\langle \mb{1},e^{tA^{\xi}}\theta\rangle}.
	\]
	Equivalently, one can check that $\varphi^{\xi}$ is the flow associated to the vector field $H^\xi_0$ on the simplex $\Delta$. The fact that $A^{\xi}$ is Metlzer and irreducible entails that for all $t>0$, $e^{tA^{\xi}}$ belongs to $M^d_{++}$, the set of $d\times d$ matrices with positive entries. Indeed,  $e^{tA^{\xi}}=e^{-rt}e^{t(A^{\xi}+rI_d)}$ where $r>0$ can be chosen large enough so that $A^{\xi}+rI_d$ is nonnegative and irreducible in the usual sense. Positive matrices have very nice geometrical properties, with respect to the Hilbert geometry. Let $d_H:\mb{R}^d_{++}\times \mb{R}^d_{++}\to\mb{R}_+$ denote the \emph{Hilbert metric}, defined by
	\[
	d_{H}(x,y)=\log\frac{\max_{i\in\intbk{1,d}}(x_i/y_i)}{\min_{i\in\intbk{1,d}}(x_i/y_i)}.
	\]
	This is in fact a pseudo-metric ($d_H(x,y)=0\Leftrightarrow \exists \lambda>0, x=\lambda y$), which induces a true metric on $\Delta_{++}:=\Delta\cap\mb{R}^d_{++}$. The induced topology is equivalent to the usual one and moreover,
	\begin{align}\label{boundmetric}
	\norm{\theta_1-\theta_2}\leq e^{d_H(\theta_1,\theta_2)}-1
	\end{align}
	for all $\theta_1,\theta_2\in\Delta_{++}$. A theorem of Birkhoff entails that every positive matrix $T$ induces a strict contraction with respect to the Hilbert metric, i.e. there exists $\kappa[T]<1$ such that
	\[
	d_H(Tx,Ty)\leq \kappa[T] d_H(x,y)
	\]
	for all $x,y\in\mb{R}^d_{++}$. This entails that for all $t>0$, $\xi\in E$ and $\theta_1,\theta_2\in\mb{R}^d_{++}$, we have
	\begin{align}\label{contractdH}
	d_H(\varphi^{\xi}(t,\theta_1),\varphi^{\xi}(t,\theta_2))\leq \kappa_{t}\, d_H(\theta_1,\theta_2)
	\end{align}
	with $\kappa_t:=\max_{\xi\in E}\kappa[e^{t A^{\xi}}]<1$.
	Let $\theta_1,\theta_2\in \Delta$, $\xi\in E$, $p\in\mb{R}$, and $\varepsilon>0$. We write $\hat{\theta}_i=(\theta_i,\xi)$ and
	\[
	M_p^{\hat{\theta}_i}(t)=\exp\left(p\int_0^t G_0(\hat{\Theta}^{\hat{\theta}_i}_s)\mr{d}s\right).
	\]
	For all $t>0$ we have
	\[
	M_p^{\hat{\theta}_1}(t)=M_p^{\hat{\theta}_2}(t)\exp\left[p\int_0^t\left(G_0(\Theta^{\hat{\theta}_1}_s,\Xi^{\xi}_s)-G_0(\Theta^{\hat{\theta}_2}_s,\Xi^{\xi}_s)\right)\mr{d}s\right],
	\]	
	hence, using Hölder's inequality, 
	\begin{align}\label{compar}
	\frac{1}{t}\log\mb{E}\left[M_p^{\hat{\theta}_1}(t)\right]\leq \frac{1}{(1+\varepsilon)t}\log \mb{E}\left[M_{p(1+\varepsilon)}^{\hat{\theta}_2}(t)\right]+\frac{\varepsilon}{(1+\varepsilon)t}\log \mb{E}\left[L^{\hat{\theta}_1,\hat{\theta}_2}_{p,\varepsilon}(t)\right]
	\end{align}
	where, setting $C_G=\sup_{\zeta\in E}\norm{G^{\zeta}_0}_{\mr{Lip}}$,
	\[
	L^{\hat{\theta}_1,\hat{\theta}_2}_{\varepsilon}(t)=\exp\left(|p|(1+\varepsilon^{-1})C_G\int_0^t\norm[\big]{\Theta^{\hat{\theta}_1}_s-\Theta^{\hat{\theta}_2}_s}\mr{d}s\right).
	\]
	Let us bound $\mb{E}\left[L^{\hat{\theta}_1,\hat{\theta}_2}_{p,\varepsilon}(t)\right]$. For all $t>0$, $\Theta^{\hat{\theta}_1}_t,\Theta^{\hat{\theta}_2}_t\in\Delta_{++}$ and  \eqref{contractdH} implies that $t\mapsto d_H(\Theta^{\hat{\theta}_1}_t,\Theta^{\hat{\theta}_2}_t)$ is decreasing. Moreover, we can show that, roughly speaking, the rate of decrease is exponential provided that the environmental process $\Xi^{\xi}$ does not switch too often. Let $\delta>0$, and  for all $k\in\mb{N^*}$, set 
	\[
	\chi_k=\mb{1}_{\left\{\scr{N}_{\Xi}\left(](k-1)\delta,k\delta]\times [0,\overline{q}]\times E\right)=0\right\}},
	\]
	where $\overline{q}=\max_{\zeta\neq \zeta'}q(0,\zeta,\zeta')$.
	The variables $\chi_{k}, k\in\mb{N}^*$ are i.i.d. with $\mb{P}(\chi_k=1)=e^{-\delta |E|\overline{q}}$, and have the property that, for all $k\geq 2$,
	\begin{align*}
		\left\{\chi_k=1\right\}&\subset \left\{\forall t\in [(k-1)\delta,k\delta[,\ \Xi^{\xi}_{t}=\Xi^{\xi}_{(k-1)\delta}\right\} \\
		&\subset \left\{\forall i\in\left\{1,2\right\},  \Theta^{\hat{\theta}_i}_{k\delta}=\varphi^{\Xi^\xi_{(k-1)\delta}}\left(\delta,\Theta^{\hat{\theta}_i}_{(k-1)\delta}\right)\right\}\\
		&\subset  \left\{d_H\left(\Theta^{\hat{\theta}_1}_{k\delta},\Theta^{\hat{\theta}_2}_{k\delta}\right)\leq \kappa_{\delta}d_H\left(\Theta^{\hat{\theta}_1}_{(k-1)\delta},\Theta^{\hat{\theta}_2}_{(k-1)\delta}\right)\right\}
	\end{align*}
	using \eqref{contractdH} for the last line. In addition, if $k_0=\inf\left\{k\in\mb{N}^*: \chi_k=1\right\}$, then \[d_H(\Theta^{\hat{\theta}_1}_{k_0\delta},\Theta^{\hat{\theta}_2}_{k_0\delta})\leq D_{\delta}<\infty,\] where $D_{\delta}$ is the $d_H$-diameter of the compact $\left\{\varphi^{\zeta}(\delta,\theta),\ (\theta,\zeta)\in\Delta\times E\right\}\subset \Delta_{++}$. It follows that for all $n\in\mb{N}$, if we set $S_n=\sum_{k=1}^n\chi_k$,
	\[
	\left\{S_n\geq 1\right\}\subset \left\{\forall t\in [n\delta,n+1\delta[,\,d_H\left(\Theta^{\hat{\theta}_1}_{t},\Theta^{\hat{\theta}_2}_{t} \right)\leq \kappa_{\delta}^{S_n-1}D_{\delta}\right\}.
	\]	
	Combining this with \eqref{boundmetric} and using that $e^{u}-1\leq ue^{u}$ for all $u\geq 0$ and $\sup_{\theta,\theta'\in\Delta}\norm{\theta-\theta'}\leq d$, we obtain that for all $n\in\mb{N}$ and $t\in[n\delta,(n+1)\delta[$, 
	\begin{align}\label{cvgtraj}
	\norm[\big]{\Theta^{\hat{\theta}_1}_{t}-\Theta^{\hat{\theta}_2}_{t}}\leq C_{\delta}\kappa_{\delta}^{S_n},
	\end{align}
	where $C_{\delta}=d\vee\left(\kappa_{\delta}^{-1}D_{\delta}e^{D_{\delta}}\right)$.
	
	Set $C_{p,\varepsilon}=|p|(1+\varepsilon^{-1})C_G$, and for all $n\in\mb{N}$, set $L_n=L^{\hat{\theta}_1,\hat{\theta}_2}_{p,\varepsilon}(n\delta)$ and $u_n=\mb{E}\left[L_n\right]$. We have 
	\[
	L_{n+1}\leq L_n\exp\left(C_{p,\varepsilon}\delta C_{\delta}\,\kappa_{\delta}^{n/2}\right)\mb{1}_{\left\{S_n\geq  n/2\right\}}+ \exp(C_{p,\varepsilon}d\delta n)\mb{1}_{\left\{S_n< n/2\right\}},
	\]
	hence
	\begin{align}
	u_{n+1}\leq a_nu_n+b_n \label{ineqrec}
	\end{align}
	where $a_n=\exp\left(C_{p,\varepsilon}\delta C_{\delta}\,\kappa_{\delta}^{n/2}\right)$ and $b_n=\exp(C_{p,\varepsilon}d\delta n)\mb{P}\left[S_n<n/2\right]$. We may assume that $\delta$ is chosen small enough so that $e^{-\delta |E|\overline{q}}=\mb{E}[\chi_1]>1/2$. Then, we can use the Cramér-Chernoff bound 
	\[
	\mb{P}(S_n<n/2)\leq e^{-nI_\delta(1/2)},
	\]
	where $I_{\delta}$ is the Cramér transform of the Bernouilli distribution of parameter $e^{-\overline{q}|E|\delta}$, given by $I_{\delta}(r)=r\log(re^{\overline{q}|E|\delta})+(1-r)\log\left((1-r)/(1-e^{-\overline{q}|E|\delta})\right)$ for $r\in (0,1)$. Since $I_{\delta}(1/2)\rightarrow +\infty$ as $\delta\rightarrow 0$, we may choose $\delta=\delta_{p,\varepsilon}$ small enough so that $B_{p,\varepsilon}:=I_{\delta_{p,\varepsilon}}(1/2)-C_{p,\varepsilon}d\delta_{p,\varepsilon}>0$.
	By an elementary induction \eqref{ineqrec} entails that for all $n\in\mb{N}$,
	\begin{align}
	u_n\leq \left(\prod_{k=0}^{n-1}a_k\right)u_0 + \sum_{j=1}^n\left(\prod_{k=j}^{n-1}a_k\right)b_{j-1}
	\leq \left(\prod_{k=0}^{\infty}a_k\right)\left(1+\sum_{j=0}^{\infty}b_j\right) \leq C'_{p,\varepsilon}
	\end{align}
	where $C'_{p,\varepsilon}=\exp\left(C_{p,\varepsilon}\delta C_{\delta}(1-\sqrt{\kappa_{\delta}})^{-1}\right)\left(1+(1-e^{-B_{p,\varepsilon}})^{-1}\right).$ Thus,
	\[
	\sup_{t\geq 0}\mb{E}\left[L_{p,\varepsilon}^{\hat{\theta}_1,\hat{\theta}_2}(t)\right]= \sup_{n\geq 0}u_n\leq C'_{p,\varepsilon}.
	\]
	Since constant $C'_{p,\varepsilon}$ does not depend on $(\theta_1,\theta_2,\xi)$, by plugging this into \eqref{compar}, we obtain that for all $\xi\in E$ and all $t>0$,  
	\begin{align}\label{compar2}
		\frac{1}{t}\log\sup_{\theta_1\in\Delta}\mb{E}\left[M_p^{(\theta_1,\xi)}(t)\right]\leq \frac{1}{(1+\varepsilon)t}\log \inf_{\theta_2\in\Delta}\mb{E}\left[M_{p(1+\varepsilon)}^{(\theta_2,\xi)}(t)\right]+\frac{\varepsilon}{(1+\varepsilon)t}\log (C'_{p,\varepsilon}).
	\end{align}
	Since $Q(0)$ is irreducible, $c_1:=\inf_{\xi,\xi'\in E}\mb{P}\left[\Xi^{\xi'}_1=\xi\right]>0$. Thus, using the Markov property at time 1 we obtain that for all $t\geq 1$, $\theta'\in E$ and $\xi,\xi'\in E$,
	\begin{align*}
	\mb{E}\left[M^{(\theta',\xi')}_{p(1+\varepsilon)}(t)\right]
	&\geq e^{-|p|(1+\varepsilon)\norm{G_0}_{\infty}} \mb{E}\left[\mb{1}_{\left\{\Xi^{\xi'}_1=\xi\right\}}\exp\left(p\int_1^tG_0(\hat{\Theta}^{(\theta',\xi')}_s)\mr{d}s\right)\right] \\ 
	&\geq e^{-|p|(1+\varepsilon)\norm{G_0}_{\infty}}c_1\inf_{\theta_2\in\Delta}\mb{E}\left[\exp\left(-p(1+\varepsilon)\int_0^{t-1}G_0(\hat{\Theta}^{(\theta_2,\xi)}_s)\mr{d}s\right)\right] \\
	&\geq e^{-2|p|(1+\varepsilon)\norm{G_0}_{\infty}}c_1\inf_{\theta_2\in\Delta}\mb{E}\left[M^{(\theta_2,\xi)}_{p(1+\varepsilon)}(t)\right].
	\end{align*}
	Combining this with \eqref{compar2} yields, for all $t\geq 1$,
	\begin{align*}
	&\ \frac{1}{t}\log\sup_{\hat{\theta}\in\Delta\times E}\mb{E}\left[M_p^{\hat{\theta}}(t)\right]\\&\leq \frac{1}{(1+\varepsilon)t}\left(\log \inf_{\hat{\theta}\in\Delta\times E}\mb{E}\left[M_{p(1+\varepsilon)}^{\hat{\theta}}(t)\right]+\log(e^{2|p|(1+\varepsilon)\norm{G_0}_{\infty}}c_1^{-1})+\varepsilon\log (C'_{p,\varepsilon})\right).
	\end{align*}
	Hence, letting $t\rightarrow +\infty$, we obtain that $\lambda_p\leq \nu_{p(1+\varepsilon)}/(1+\varepsilon)$. Moreover we have $\nu_{p(1+\varepsilon)}\leq\lambda_{p(1+\varepsilon)}$ by definition, thus, for all $q\in\mb{R}$ and $\varepsilon>0$, 
	\[
	(1+\varepsilon)\lambda_{q/(1+\varepsilon)}\leq \nu_q\leq \lambda_q.
	\]
	Since $q\mapsto \lambda_q$ is convex on $\mb{R}$, it is continuous, hence, letting $\varepsilon\rightarrow 0$, we obtain that $\nu_q=\lambda_q =:g(q)$. This ends the proof of item \ref{item1plyap}.
	
We now prove  item \ref{item2bisplyap}. We already know that $p\mapsto\lambda_p$ is convex, hence $g$ is convex. Now let us investigate $\lim_{p\rightarrow +\infty}g(p)/p$. This limit exists due to the convexity of $g$. For all $p>0$, $\xi\in E$ and $t\geq 0$ we have, recalling that $\theta^{\xi}_*$ is an equilibrium of the flow $\varphi^\xi$ and $G_0(\theta^{\xi}_*,\xi)=\Lambda^\xi$, 
 	\begin{align*}
 	\mb{E}\left[\exp\left(p\int_0^tG_0(\hat{\Theta}^{(\theta^{\xi}_*,\xi)}_s)\mr{d}s\right)\right]&\geq \mb{P}\left[\forall s\in [0,t],\ \Xi^{\xi}_s=\xi\right] \exp\left(ptG_0(\theta^{\xi}_*,\xi)\right)\\ &\geq \exp\left(p\Lambda^{\xi}t\right)\exp\left(-|q(0,\xi,\xi)|t\right),
 	\end{align*}
 	hence, when $p>0$, $\overline{f}_p(t)/t\geq p\Lambda^{\xi}-|q(0,\xi,\xi)|$.
 	Dividing by $p$ and letting $t\rightarrow +\infty$ and then $p\rightarrow +\infty$ yields
 	\[
 	\lim_{p\rightarrow +\infty}\frac{g(p)}{p}\geq \max_{\xi\in E}\Lambda^{\xi}.
 	\]
	Let $\hat{\theta} \in  \Gamma( \hat \Theta)$. Since $\Gamma( \hat \Theta)$ is positively invariant by the flow $\varphi^{\xi'}$ for all $\xi' \in E$, we get that $\hat{\Theta}_t^{\hat{\theta}} \in \Gamma( \hat \Theta)$ almost surely for all $t\geq 0$. In particular, for all $p>0$, 
	\[
	\underline{f}_p(t)\leq \log\mb{E}\left[\exp\left(p\int_0^t G_0(\hat{\Theta}^{\hat{\theta}}_s)\mr{d}s\right)\right]\leq pt\max_{\Gamma( \hat \Theta)}G_0,
	\]
	hence $g(p)/p\leq \max_{\Gamma( \hat \Theta)}G_0$.
	The proofs are similar for the bounds when $p$ goes to $- \infty$.
	
	Let us prove item \ref{item2plyap}. Obviously $g(0)=0$. Let us show that $g'(0)$ exists and equals $\Lambda$. We claim that 
	\begin{align}\label{boundI}
	I:=\int_0^\infty\mr{osc}(P^{\hat{\Theta}}_sG_0)\mr{d}s <\infty
	\end{align}
    where $\mr{osc}(f)=\sup f- \inf f$ denotes the oscillation of $f$.
	Let us prove this claim.
   	Let $\hat \theta = (\theta, \zeta)$, $\hat \theta' = (\theta', \zeta')\in \Delta\times E$. For all $t\geq 0$, we have
	\begin{align}\label{eq:oscPtheta}
  	P^{\hat{\Theta}}_tG_0(\hat{\theta})-P^{\hat{\Theta}}_tG_0(\hat{\theta}')= &\ \mb{E}\left[\left( G_0( \Theta^{\hat \theta }_t) - G_0( \Theta^{\hat \theta' }_t) \right) 1_{\{T_c(\zeta,\zeta') > t/2\}} \right] \nonumber \\& +\mb{E}\left[\left( G_0( \Theta^{\hat \theta }_t) - G_0( \Theta^{\hat \theta' }_t) \right) 1_{\{T_c(\zeta,\zeta') \leq t/2\}} \right]. 
 	\end{align}
  where $T_c(\zeta,\zeta')$ denotes the first coalescent time of $\Xi^{\zeta}$ and $\Xi^{\zeta'}$, i.e. $T_c(\zeta,\zeta')=\inf \{t\geq 0 :\Xi^{\zeta}_t=\Xi^{\zeta'}_t\}$. By construction, $\Xi^{\zeta}$ and $\Xi^{\zeta'}$ coincide after $T_c(\zeta,\zeta')$ and moreover, using the irreducibility of $Q(0)$, $\inf_{\zeta_1,\zeta_2\in E}\mb{P}\left(T_c(\zeta_1,\zeta_2)\leq 1\right)>0$. This entails, using the Markov property, that there exists $C,a>0$ such that for all $\zeta_1,\zeta_2\in E$ and $t\geq 0$, $\mb{P}(T_c(\zeta_1,\zeta_2)>t)\leq Ce^{-at}$. Hence, the absolute value of the first term in the right handside of \eqref{eq:oscPtheta} can be bounded by $2 \|G_0\|_{\infty} Ce^{ - a t/2}$.  As for the second term, the Markov property at time $t/2$ yields
  \[
  \left|\mb{E}\left[\left( G_0( \Theta^{\hat \theta }_t) - G_0( \Theta^{\hat \theta' }_t) \right) 1_{\{T_c(\zeta,\zeta') \leq t/2\}} \right]\right| \leq \sup_{\theta_1,\theta_2\in \Delta,\ \xi\in E} \left|\mb{E} \left[ G_0( \hat{\Theta}^{(\theta_1,\xi)}_{t/2}) - G_0( \hat{\Theta}^{(\theta_2,\xi)}_{t/2})  \right] \right|.
  \]	
  Fix some $\delta>0$, and let $\theta_1,\theta_2\in\Delta$, $\xi\in E$. Let $(\chi_k)_{k\in\mb{N}^*}$, $(S_n)_{n\in\mb{N}}$ and $C_{\delta}$ be as above (in the proof of item \ref{item1plyap}). Then for all $u\geq 0$ we have, recalling \eqref{cvgtraj},
  \begin{align*}
  	\mb{E}\left[\norm[\big]{\Theta^{(\theta_1,\xi)}_u-\Theta^{(\theta_2,\xi)}_u}\right]\leq C_{\delta}\mb{E}\left[\kappa_{\delta}^{S_{\lfloor u/\delta\rfloor}}\right]=C_{\delta}\mb{E}\left[\kappa^{\chi_1}\right]^{\lfloor u/\delta\rfloor}\leq C'e^{-a' u}
  \end{align*}
  where $C'=C_{\delta}/\mb{E}\left[\kappa^{\chi_1}\right]$ and $a'=-\log\mb{E}\left[\kappa^{\chi_1}\right]/\delta>0$. It follows that the absolute value of the second term in the right handside of \eqref{eq:oscPtheta} can be bounded by $C_GC'e^{-a't/2}$. Hence, plugging the bounds into \eqref{eq:oscPtheta} yields, for all $t\geq 0$, 
  \[
  \mr{osc}(P^{\hat{\Theta}}_tG_0)\leq 2\norm{G_0}_{\infty}Ce^{-at/2}+C'e^{-a't/2},
  \]
  which finally proves the claim \eqref{boundI} with $I=4\norm{G_0}_{\infty}C/a+2C'/a'$.
  
  A useful consequence of \eqref{boundI}is that	\begin{align}\label{encadpiG0}
  		\sup_{t\geq 0,\ \hat{\theta}\in\Delta\times E}\left|t\Lambda-\int_0^tP^{\hat{\Theta}}_sG_0(\hat{\theta})\mr{d}s\right|\leq I<\infty,
  	\end{align}	
    using that $t\Lambda=t\pi_0(G_0)=\int_0^t\pi_0 P^{\hat{\Theta}}_sG_0\mr{d}s$ and $\left|\pi_0 P^{\hat{\Theta}}_sG_0-P^{\hat{\Theta}}_sG_0(\hat{\theta})\right|\leq \mr{osc}(P^{\hat{\Theta}}_s)$.    
    Let $t\geq 0$ and $\hat{\theta}\in\mb{\Delta}\times E$. Set $Z=\int_0^tG_0(\Theta^{\hat{\theta}}_s)\mr{d}s$ and let $\ell:\mb{R}\mapsto\mb{R}_+$ denote its cumulant generating function, i.e. $\ell(p)=\log\mb{E}\left[e^{pZ}\right]$. The function $\ell$ is of class $\mathcal{C}^{\infty}$ and for all $p\in\mb{R}$ we have $\ell'(p)=\mb{E}_p[Z]$ and $\ell''(p)=\mb{E}_p[Z^2]-\mb{E}_p[Z]^2$, where $\mb{E}_p$ is the expectation operator associated to the  probability distribution of density $e^{pZ}/\mb{E}[e^{pZ}]$ with respect to $\mb{P}$. Hence $0\leq\ell''\leq t^2\norm{G_0}_\infty$ and a second order Taylor expansion at 0 yields, for all $p\in\mb{R}$,
    \[\ell(p)\leq p\int_0^tP^{\hat{\Theta}}_sG_0(\hat{\theta})\mr{d}s+\frac{1}{2}p^2t^2\norm{G_0}_{\infty}.
    \]
    Combining this with \eqref{encadpiG0}, we obtain that
    \[
    \overline{f}_p(t)\leq p\Lambda t+|p|I+\frac{1}{2}p^2t^2\norm{G_0}_{\infty}.
    \]
    Now, by the subadditive theorem we have $g(p)=\lambda_p=\inf_{s>0}\overline{f}_p(s)/s$, hence
    \[
    g(p) \leq p\Lambda+\frac{|p|I}{t}+\frac{1}{2}p^2t\norm{G_0}_{\infty}.
    \]
    Since $g$ is convex, it has a right (resp. left) derivative at $0$, denoted $g'_{+}(0)$ (resp. $g'_{-}(0)$), and the above inequality yields, for all $t>0$,
   	\[
   	g'_{+}(0)=\lim_{\substack{p\rightarrow 0 \\ p>0}}\frac{g(p)}{p}\leq \Lambda+\frac{I}{t},\quad g'_{-}(0)=\lim_{\substack{p\rightarrow 0 \\ p<0}}\frac{g(p)}{p}\geq \Lambda-\frac{I}{t}.
   	\]
 	Letting $t$ go to $+\infty$ and using that $g'_{-}(0)\leq g'_{+}(0)$ by convexity, we conclude that $g$ is differentiable at $0$ with
 	\[
 	g'(0)=\Lambda.
 	\]

%
	
	We now prove item  \ref{item3plyap}. Assume that $\Lambda>0$. Recall that for all $\hat{y}=(y,\xi)\in \mb{R}_+^d\times E$, $Y^{\hat{y}}$ denotes the spatial component of the linear PDMP given by $Y^{\hat{y}}_0=y$ and $\dot{Y}^{\hat{y}}_t=A^{\Xi_t}Y^{\hat{y}}_t$. Assume first that $0 \in \Gamma(Y)$. Then, Proposition 3.14 in \cite{BLMZ} implies that there exist $\delta > 0$, a finite open coverage $\mathcal{O}_1, \ldots, \mathcal{O}_k$ of $\Delta$, and positive numbers $t_1 \leq  \ldots \leq t_k$ such that, for all $\hat y = (y, \xi) \in \Delta \times E$, 
	\[
	y \in \mathcal{O}_k \quad \Rightarrow \quad \mb{P}\left( \|Y^{\hat{y}}_{t_k}\| \leq \frac{1}{2}\right) \geq \delta.
	\] 
By linearity of $ y \mapsto Y^{\hat{y}}_t$, we deduce that there exists a measurable map $T : \mb{R}_+^d \to \mb{R}_+$ such that, for all $\hat y \in \mb{R}^d_+ \times E$,	
\begin{equation}
\label{eq:0acclin}
\mb{P}\left( \|Y^{\hat{y}}_{T(y)}\| \leq \frac{\|y\|}{2}\right) \geq \delta.
\end{equation}
	Moreover, $T(y) \in \{ t_1, \ldots, t_k\}$ for all $y \in \mb{R}_+^d$. Let us fix $\hat{y}\in(\mb{R}_+^d\setminus\{0\})\times E$, and write $Y=Y^{\hat{y}}$. We construct a increasing sequence of stopping times $(T_n)_{n \in\mb{N}}$, defined as follows :
	\[
	T_0 = 0, \qquad T_{n + 1} = T_n + T(Y_{T_n}), \, n \in\mb{N}.
	\]
Note that since $T(y) \geq t_1>0$ for all $y \in \mb{R}_+^d$, $\lim_{n \to \infty} T_n = + \infty$.	We define the event 
	\[
	E_n = \bigcap_{i= 0}^{n-1} \left\lbrace \|Y_{T_{n+1}}\| \leq \frac{\|Y_{T_n}\|}{2} \right\rbrace.
	\]
	The strong Markov property and Equation \eqref{eq:0acclin} imply that $\mb{P}(E_n) \geq \delta^n.$ We let $(N_t)_{t\in\mb{R}_+}$ be the counting process associated to the $(T_n)_{n \in\mb{N}}$ : $N_t = \sup\{ n\in\mb{N} : \: T_n \leq t\},$ and we set $n_t = \lfloor \frac{t}{t_k} \rfloor$ and $m_t = \lfloor \frac{t}{t_1} \rfloor$. For all $t\geq 0$ we have $n_t \leq N_t \leq n'_t$, and, on the event $E_{n'_t}$, 
	\[
	\|Y_t\| \leq \|Y_0\| 2^{ - N_t} e^{ (t - N_t) \|G_0\|_{\infty} } \leq \|Y_0\| 2^{-n_t} M,
	\]
	where $M= e^{t_k \|G_0\|_{\infty}}$. Thus, for all $\hat y \in \Delta \times E$ and all $t\geq 0$,
	\[
	\mb{P}\left( \|Y_t\|\leq 2^{-n_t}M \right) \geq \mb{P}( E_{n'_t} ) \geq \delta^{n'_t}
	\]
	This yields
	\[
	\mb{P}\left( \|Y_t\|\leq 2 M e^{ - c_1 t} \right) \geq e^{ - c_2 t},
	\]
	where $c_1 = \log(2)/t_k$ and $c_2 = - \log(\delta)/t_1$, and thus
	\[
	\mb{E}\left( \|Y_t\|^{-p} \right) \geq (2M)^{-p} e^{(c_1 p - c_2)t}
	\]
	This implies that for all $p > c_2/c_1$, $g(-p) > 0$. Thus, item \ref{item3plyap} follows from the convexity of $g$ and the fact that $g'(0) = \Lambda > 0$. Moreover, $p^* \leq c_2/c_1$.

Assume now that $0$ does not belong to $\Gamma(Y)$. By definition, there exists $\varepsilon > 0$ and $\hat{y} \in \Delta\times E$ such that, $\mb{P}$-almost surely, $\|Y_t^{\hat y}\| \geq \varepsilon$ for all $t \geq 0$. In particular, for all $p > 0$, $\log \mb{E} ( \|Y_t^{\hat y}\|^{-p} ) \leq  - p \log( \varepsilon)$. This implies $g(-p) \leq 0$ for all $p>0$, hence $g(-p)<0$ for all $p>0$ due to the convexity of $g$ and the fact that $g'(0)>0$. Thus, $p^* = + \infty$, which ends the proof of item \ref{item3plyap}.

The proof of item \ref{item4plyap} is very similar, we just give a sketch of it. If $\infty \notin \Gamma(Y)$, there exists $M > 0$ and $y \in \Delta$, such that, almost surely, $\|Y_t^{\hat y}\| \leq M$. This implies that $g(p) < 0$ for all $p > 0$, hence $p_* = + \infty$. If $\infty \in \Gamma(Y)$, we can adapt the proof of Lemma 3.16 (and thus of Proposition 3.14) in \cite{BLMZ} to show that  that there exist $\delta > 0$, a finite open coverage $\mathcal{O}_1, \ldots, \mathcal{O}_k$ of $\Delta$, and positive numbers $t_1 \leq  \ldots \leq t_k$ such that, for all $\hat y = (y, \xi) \in \Delta \times E$, 
	\[
	y \in \mathcal{O}_k \quad \Rightarrow \quad \mb{P}\left( \|Y^{\hat{y}}_{t_k}\| \geq 2 \right) \geq \delta.
	\] 
	(just replace the definition of $\mathcal{O}(\mb{i}, \mb{u}, \beta)$ in the proof of Lemma 3.16 in \cite{BLMZ} by $\mathcal{O}(\mb{i}, \mb{u}, \beta) = \{ x \in \Delta \: : \| \mb{\Phi}^\mb{i}_\mb{u}(x) \| \in [3,4], \: p(x, \mb{i}, \mb{u}) > \beta\}$) . As before, we can deduce from this that for some measurable map $T$, and some $\delta > 0$, $\mb{P} ( \|Y^{\hat y}_{T(y)}\| \geq 2 \|y\|) \geq \delta$. Reasoning like in the proof of item \ref{item3plyap}, we deduce that $g(p) > 0$ for $p$ large enough, and therefore, $p_* < + \infty$. 

	\subsection{Proof of Theorem \ref{mintimeextpers}}\label{sec:mintimeextpers}
	
	\begin{Pro}\label{prolyapXK}
		Let $p  \in \mb{R}$ such that $g(p) < 0$. Let $T>0$ and $\eta>0$ be given by Proposition \ref{prolyapPDMP}. There exist $\gamma\in(0,1)$, $a>0$ and $K_0\in\mb{N}^*$ such that for all $K\geq K_0$ and all $\hat{x}=(x,\xi)\in\hat{\mathcal{X}}_+^K$ satisfying $a/K\leq \norm{x} \leq \eta$, we have
		\begin{align}\label{lyap}
			\mb{E}^{\hat{x}}\left[\norm{X^K_T}^{p}\mb{1}_{\left\{T<\tau_0^K\right\}}\right]\leq \gamma\, \norm{x}^{p}.
		\end{align}
	\end{Pro}
	\prf We first make the proof in the case where $p  > 0 $ is such that $g( - p) < 0$. Let $T,\gamma_0,\eta$ be given by Proposition \ref{prolyapPDMP}. Let $C_T,\varepsilon_T$ be given by Lemma \ref{lemcoupl}, set $M=2e^{C_FT}$ and $m=e^{-C_FT}/2$ and let $C'$ and $K_0$ be given by Lemma \ref{lemtauq}. Let $K\geq K_0$ and $\hat{x}=(x,\xi)\in\hat{\mathcal{X}}^K_+$ such that $\norm{x}\leq \eta$. For all $\varepsilon \in (0,\varepsilon_T]$, we have, combining Proposition \ref{prolyapPDMP}, Lemma \ref{lemcoupl} and Lemma \ref{lemtauq}, 
	\begin{align*}
	&\quad \ \mb{E}^{\hat{x}}\left[\norm{X^K_T}^{-p}\mb{1}_{\left\{T<\tau_0^K\right\}}\right] \\
	&=\mb{E}^{\hat{x}}\left[\norm{X^K_T}^{-p}\mb{1}_{\left\{T<\sigma^K_{\varepsilon}\right\}}\right]+\mb{E}^{\hat{x}}\left[\norm{X^K_T}^{-p}\mb{1}_{\left\{\sigma^K_{\varepsilon}\leq T <\overline{\tau}^{K}_{M\norm{x}}\wedge\underline{\tau}^{K}_{m\norm{x}}\right\}}\right] \\
	&\quad+\mb{E}^{\hat{x}}\left[\norm{X^K_T}^{-p}\mb{1}_{\left\{\sigma^K_{\varepsilon}\vee \left(\overline{\tau}^{K}_{M\norm{x}}\wedge\underline{\tau}^{K}_{m\norm{x}}\right)\leq T <\tau_0^K\right\}}\right] \\
	&\leq (1-\varepsilon)^{-p}\mb{E}^{\hat{x}}\left[\norm{X_T}^{-p}\right] +(m\norm{x})^{-p}\mb{P}^{\hat{x}}\left[\sigma^K_{\varepsilon}\leq T\right]+ (1/K)^{-p}\mb{P}^{\hat{x}}\left[\overline{\tau}^{K}_{M\norm{x}}\wedge\underline{\tau}^{K}_{m\norm{x}}\leq T\right] \\
	&\leq \norm{x}^{-p}\left((1-\varepsilon)^{-p}\gamma_0+2m^{-p}de^{-C_T K\norm{x}\varepsilon^2}\right.\\
	& \qquad \qquad\left. + m^{-p}T|E|\sup_{\xi_1\neq \xi_2,\ \norm{z-y}\leq d\varepsilon}\left|q(z,\xi_1,\xi_2)-q(y,\xi_1,\xi_2)\right|+2(K\norm{x})^{p}e^{-C'K\norm{x}}\right).
	\end{align*}
	Since the functions $q(\cdot,\xi_1,\xi_2)$ are uniformly continuous, we can choose $\varepsilon>0$ small enough so that
	\[
	\gamma_1:=(1-\varepsilon)^{-p}\gamma_0+m^{-p}T|E|\sup_{\xi_1\neq \xi_2,\ \norm{z-y}\leq d\varepsilon}\left|q(z,\xi_1,\xi_2)-q(y,\xi_1,\xi_2)\right|<1,
	\]
	and then, we can choose $a>0$ large enough so that 
	\[
	\gamma:=\gamma_1+2m^{-p}de^{-C_Ta\varepsilon^2}+\sup_{y\geq a}\left(2y^pe^{-C'y}\right)<1.
	\]
	We can conclude : for all $K\geq K_0$ and all $\hat{x}=(x,\xi)\in\hat{\mathcal{X}}^K_+$ such that $a/K\leq \norm{x}\leq \eta$, 
	\[
	\mb{E}^{\hat{x}}\left[\norm{X^K_T}^{-p}\mb{1}_{\left\{T<\tau_0^K\right\}}\right]\leq \gamma x^{-p}.
	\]	
	
In the case where $p >0$ is such that $g(p)< 0$, the proof is similar, and even simpler, with the use of the bound 
 \begin{align*}
 \mb{E}^{\hat{x}}\left[\norm{X^K_T}^{p}\mb{1}_{\left\{T<\tau_0^K\right\}}\right] \leq &(1+\varepsilon)^{p}\mb{E}^{\hat{x}}\left[\norm{X_T}^{p}\right] +(M\norm{x})^{p}\mb{P}^{\hat{x}}\left[\sigma^K_{\varepsilon}\leq T\right]\nonumber\\&+d^p \mb{P}^{\hat{x}}\left[\overline{\tau}^{K}_{M\norm{x}}\wedge\underline{\tau}^{K}_{m\norm{x}}\leq T\right].
 \end{align*}
	\hfill $\square$
	
	\begin{Lem}\label{hittingtime1}
		Assume that $\Lambda>0$ and let $p\in (0,p^*)$. Let $\eta>0$ be given by Proposition \ref{prolyapPDMP}, item \ref{itemprolyap}. There exist $a'>0$, $C_5\geq 1$ and $K_0\in\mb{N}^*$ such that for all $K\geq K_0$, $\rho>0$ and $\hat{x}=(x,\xi)\in\hat{\mathcal{X}}^K_+$ such that $a'/K\leq \rho \leq \norm{x}\leq \eta$, we have
		\[
		\mb{P}^{\hat{x}}\left(\underline{\tau}^K_{\rho}<\overline{\tau}^{K}_{\eta}\right)\leq C_5 (\norm{x}/\rho)^{-p}.
		\]
	\end{Lem}
	\prf Let $p\in(0,p^*)$. Let $T,\eta$ be given by Proposition \ref{prolyapPDMP} and $\gamma,a,K_0$ be given by Proposition \ref{prolyapXK} respectively. Let $K\geq K_0$ and set, for all $n\in\mb{N}$,
	\[
	Z^K_n=\norm[\big]{X^K_{nT}}^{-p}\mb{1}_{\left\{nT<\tau_0^K\right\}}.
	\]
	Let $a'\geq a$, which remains to be fixed later in the proof. Let $\rho>0$ and an initial condition $\hat{x}=(x,\xi)\in\hat{\mathcal{X}}^K_+$ be given for $\hat{X}^K$, such that $a'/K\leq \rho \leq \norm{x} \leq \eta$. Set
	\[
	\scr{T}^K=\left\lceil\frac{\overline{\tau}^K_{\eta}\wedge\underline{\tau}^K_{\rho}}{T}\right\rceil.
	\]
	Then the process $(Z^K_{\scr{T}^K\wedge n})_{n\in\mb{N}}$ is a non-negative, bounded supermartingale, with respect to the filtration $(\scr{F}_{nT})_{n\in\mb{N}}$. Indeed,
	\begin{align*}
		\mb{E}^{\hat{x}}\left[Z^K_{\scr{T}^K\wedge (n+1)}\big| \scr{F}_{nT}\right]&=Z^K_{\scr{T}^K}\mb{1}_{\left\{\scr{T}^K\leq n\right\}}+ \mb{E}^{\hat{x}}\left[Z^K_{n+1}\mb{1}_{\left\{\scr{T}^K>n\right\}}\big| \scr{F}_{nT}\right] \\
		&=Z^K_{\scr{T}^K}\mb{1}_{\left\{\scr{T}^K\leq n\right\}}+ \mb{1}_{\left\{\scr{T}^K>n\right\}}\mb{E}^{\hat{X}^K_{nT}}\left[\norm[\big]{X^K_{T}}^{-p}\mb{1}_{\left\{T<\tau_0^K\right\}}\right]\\
		&\leq Z^K_{\scr{T}^K}\mb{1}_{\left\{\scr{T}^K\leq n\right\}}+ \mb{1}_{\left\{\scr{T}^K>n\right\}}Z^K_n \\
		&= Z^K_{\scr{T}^K\wedge n},
	\end{align*}
	where the second equality comes from the Markov property of $\hat{X}^K$ and the inequality comes from \eqref{lyap}. Consequently,
	\begin{align}\label{ineqZ1}
		\mb{E}^{\hat{x}}\left[Z^K_{\scr{T}^K}\mb{1}_{\left\{\scr{T}^K<\infty\right\}}\right]\leq Z^K_0=\norm{x}^{-p}.
	\end{align}
	Set $M=2e^{C_FT}$ and $m=e^{-C_FT}/2$ and let $C'>0$ be given by Lemma \ref{lemtauq}. We may assume that $K_0$ is large enough so that $\min_{i\in\intbk{1,d}}K_i(K)/K\geq \underline{\alpha}/2$  for all $K\geq K_0$, which implies that  $\norm{X^K_{\underline{\tau}^K_{\rho}}}\geq \rho-2/(\underline{\alpha}K)$ $\mb{P}^{\hat{x}}$-almost surely. Using the strong Markov property of the process $X^K$ at $\overline{\tau}^{K}_{\eta}\wedge\underline{\tau}^{K}_{\rho}$ and Lemma \ref{lemtauq}, we get, for $K\geq K_0$, 
	
	\begin{align}\label{ineqZ2}
		\mb{E}^{\hat{x}}\left[Z^K_{\scr{T}^K}\mb{1}_{\left\{\underline{\tau}^K_{\rho}<\overline{\tau}^{K}_{\eta}\right\}}\right]&\geq \mb{P}^{\hat{x}}\left(\underline{\tau}^K_{\rho}<\overline{\tau}^{K}_{\eta}\right)\inf_{\substack{(y,\zeta)\in \hat{\mathcal{X}}^K\\\rho-2/(\underline{\alpha}K)\leq\norm{y}\leq \rho}}\mb{E}^{(y,\zeta)}\left[\inf_{0\leq t <T}\left(\norm{X^K_{t}}^{-p}\mb{1}_{\left\{t<\tau_0^K\right\}}\right)\right]\nonumber \\ &\geq \mb{P}^{\hat{x}}\left(\underline{\tau}^K_{\rho}<\overline{\tau}^{K}_{\eta}\right)(M\rho)^{-p}\left(1-\sup_{\substack{(y,\zeta)\in \hat{\mathcal{X}}^K\\\norm{y}\geq (a'-2/\underline{\alpha})/K}}\mb{P}^{(y,\zeta)}\left(\overline{\tau}^K_{M\rho}<T\right)\right) \nonumber
		\\ &\geq \mb{P}^{\hat{x}}\left(\underline{\tau}^K_{y}<\overline{\tau}^{K}_{\eta}\right)(M\rho)^{-p}(1-2e^{-C'(a'-2/\underline{\alpha})}).
	\end{align}

	Choosing $a'$ large enough so that $1-2e^{-C'(a'-2/\underline{\alpha})}\geq 1/2$ and combining \eqref{ineqZ1} with \eqref{ineqZ2}, we obtain
	\[
	\mb{P}^{\hat{x}}\left(\underline{\tau}^K_{\rho}<\overline{\tau}^{K}_{\eta}\right)\leq 2\,M^p(\norm{x}/\rho)^{-p},
	\]
	which ends the proof.
	
	\hfill $\square$

For all $\rho\geq 0$, we set $\hat{\mathcal{X}}^K_\rho=\left\{(x,\xi)\in\hat{\mathcal{X}}^K: \norm{x}\geq \rho\right\}$. 

	\begin{Lem}\label{hittingtime2}
	Assume that $\Lambda>0$, and let $p\in (0,p^*)$. Let $\eta>0$ and $a'>0$ be given by Proposition \ref{prolyapPDMP} and Lemma \ref{hittingtime1} respectively. There exist $T',c_1>0$, $\delta\in(0,\eta)$ and $K_0\in\mb{N}^*$ such that for all $K\geq K_0$, $\rho\in[a'/K,\delta]$, $\hat{x}\in \hat{\mathcal{X}}^K_{\eta}$ and $t\geq 0$, we have
	\[
	\mb{P}^{\hat{x}}\left(\underline{\tau}^K_{\rho}>t\right)\geq e^{-c_1\rho^p\lceil t/T'\rceil}.
	\]
	\end{Lem}
	\prf Let $p\in (0,p^*)$. Let $T>0$, $\eta >0$ be given by Proposition \ref{prolyapPDMP} and let $a'>0$, $K_0\in\mb{N}^*$ and $C_5\geq 1$ be given by Lemma \ref{hittingtime1}. Set $M=2e^{C_FT}$. 
	Let $K\geq K_0$ and set
	\[
	\sigma^K=\inf\left\{t\geq \underline{\tau}^K_{\eta/M}: \norm{X^K_t}\geq \eta\right\}.
	\]
	Let $\rho \in [a'/K,\eta/(2M)]$ and define $\phi_\rho^K:\mb{R}_+\to [0,1]$ by
	\[
	\phi^K(t)=\inf_{\hat{x}\in\hat{\mathcal{X}}^K_{\eta}}\mb{P}^{\hat{x}}\left(\underline{\tau}^K_\rho>t\right).
	\]
	Let $\hat{x}=(x,\xi)\in\hat{\mathcal{X}}^K_{\eta}$ and $t\geq 0$. We have
	\begin{align*}
	\mb{P}^{\hat{x}}\left(\underline{\tau}^K_\rho>t\right) &\geq \mb{P}^{\hat{x}}\left(\sigma^K<\underline{\tau}^K_{\rho},\ \underline{\tau}^K_\rho-\sigma^K>t-\sigma^K\right) \\
	&\geq \mb{P}^{\hat{x}}\left(\sigma^K<\underline{\tau}^K_\rho,\ \sigma^K \geq 2T,\ \underline{\tau}^K_\rho-\sigma^K> t-2T\right) \\&\quad +\mb{P}^{\hat{x}}\left(\sigma^K<\underline{\tau}^K_\rho,\ \sigma^K <2T,\ \underline{\tau}^K_\rho-\sigma^K>t\right).
	\end{align*}
	Using the strong Markov property of $X^K$ at time $\sigma^K$ and the definition of $\phi^K$, we get
	\begin{align}\label{survtauy}
	\mb{P}^{\hat{x}}\!\left(\underline{\tau}^K_\rho>t\right) &\geq \mb{P}^{\hat{x}}\!\left(\sigma^K<\underline{\tau}^K_\rho,\ \sigma^K\geq 2T\right)\phi^K(t-2T)+	\mb{P}^{\hat{x}}\!\left(\sigma^K<\underline{\tau}^K_\rho,\ \sigma^K< 2T\right)\phi^K(t) \nonumber\\
	&=\mb{P}^{\hat{x}}\!\left(\sigma^K<\underline{\tau}^K_\rho\right)\phi^K(t-2T)+\mb{P}^{\hat{x}}\!\left(\sigma^K<\underline{\tau}^K_\rho,\sigma^K<2T\right)\left(\phi^K(t)-\phi^K(t-2T)\right).
	\end{align}    
	Let $C'>0$ be given by Lemma \ref{lemtauq}. We may assume that $K_0$ is large enough so that for all $K\geq K_0$, $2e^{-C'K\eta/(2M)}\leq 1/4$ and $\min_{i\in\intbk{1,d}}\left(\eta/M-1/K_i\right)\geq \eta/(2M)$, which implies that $\norm{X^K_{\underline{\tau}^K_{\eta/M}}}\geq \eta/(2M)$ $\mb{P}^{\hat{x}}$-almost surely. Using the strong Markov property at time $\underline{\tau}^K_{\eta/M}$ and Lemma \ref{hittingtime1}, we obtain
	\begin{align*}
	\mb{P}^{\hat{x}}\left(\sigma^K<\underline{\tau}^K_\rho\right)&\geq \inf_{(y,\zeta)\in\hat{\mathcal{X}}^K_+,\ \norm{y}\geq \eta/(2M)}\mb{P}^{(y,\zeta)}\left(\overline{\tau}^{K}_{\eta}<\underline{\tau}^K_\rho\right) \\
	&\geq 1-C'_5\rho^p    
	\end{align*}
	where $C'_5=C_5\left(\eta/(2M)\right)^{-p}$ and, using Lemma \ref{lemtauq},
	\begin{align*}
	\mb{P}^{\hat{x}}\left(\sigma^K<\underline{\tau}^K_y,\ \sigma^K<2T\right)&\leq \mb{P}^{\hat{x}}\left(\sigma^K<2T\right) \\
	&\leq \mb{P}^{\hat{x}}\left(\underline{\tau}^K_{\eta/M}<T\right)+\mb{P}^{\hat{x}}\left(\sigma^K-\underline{\tau}^K_{\eta/M}<T\right) \\
	& \leq 2e^{-C'K \eta} + 2e^{-C'K \eta/(2M)} \\
	& \leq 1/2.
	\end{align*}
	Since the function $\phi^K_\rho$ is non-decreasing, plugging these inequalities into \eqref{survtauy} and taking the infimum on $\hat{x}\in\hat{\mathcal{X}}^K_{\eta}$ yields
	\[
	\phi^K_\rho(t)\geq \left(1-C'_5\rho^p\right)\phi^K_\rho(t-2T)+(1/2)\left(\phi^K_\rho(t)-\phi^K_\rho(t-2T)\right),
	\]
	hence
	\[
	\phi^K_\rho(t)\geq (1-2C'_5\rho^p)\phi^K_\rho(t-2T).
	\]
	Now, let us fix $\delta\in(0,\eta/(2M)]$ small enough so that so that $2C'_5\delta^p\leq 1/2$. For all $\rho\in[a'/K,\delta]$, given that $\phi^K_\rho(0)=1$ and that $\phi^K_\rho$ is non-increasing, we obtain
	\begin{align*}
	\phi^K_\rho(t)&\geq \phi^K_\rho\left(2T\left\lceil \frac{t}{2T}\right\rceil\right) \\
	&\geq \left(1-2C'_5\rho^p\right)^{\lceil t/(2T)\rceil} \\
	&\geq e^{-4C'_5\rho^{p}\lceil t/(2T)\rceil},
	\end{align*}
	where we used that $\log(1-h/2)\geq-h$ for $0\leq h\leq 1$ in the last step. \hfill $\square$
	\medskip
	
	We now turn to the proof of Theorem \ref{mintimeextpers}. Let $\eta,T',c_1,K_0$ be given by Lemma \ref{hittingtime2}. For all $K\geq K_0$, applying this proposition with $\rho=a'/K$ yields, for all $\hat{x}\in \hat{\mathcal{X}}^K_{\eta}$ and all $t\geq 0$, 
	\[
	\mb{P}^{\hat{x}}\left(\tau^K_0>t\right)\geq \mb{P}^{\hat{x}}\left(\underline{\tau}^K_{a'/K}>t\right)\geq \exp\left(-c_1\frac{a'^p}{K^p}\left\lceil \frac{t}{T'}\right\rceil\right)\geq \exp\left(-\frac{c_1a'^p}{K^p}\right)\exp\left(-\frac{C_1t}{K^p}\right)
	\]
	where $C_1=c_1a'^p/T'$. For a general initial condition $\hat{x}\in\hat{\mathcal{X}}^K_+$, the strong Markov property at time $\bar{\tau}^K_{\eta}$ then yields, for all $t\geq 0$, 
	\begin{align}\label{probasurvival}
	\mb{P}^{\hat{x}}\left(\tau^K_0>t\right)\geq \mb{P}^{\hat{x}}\left(\overline{\tau}^K_{\eta}<\tau^{K}_0\right)\exp\left(-\frac{c_1a'^p}{K^p}\right)\exp\left(-\frac{C_1t}{K^p}\right).
	\end{align}
	Let $C_5$ be given by Lemma \ref{hittingtime1} and $a''>a'$ be large enough so that $C_5(a''/a')^{-p}\leq 1/2$. Let $K\geq K_0$ and let $\hat{x}=(x,\xi)\in\hat{\mathcal{X}^K_+}$. First, if $\norm{x}\geq a''/K$, then Lemma \ref{hittingtime1} yields
	\begin{align}\label{probaescape1}
	\mb{P}^{\hat{x}}\left(\overline{\tau}^K_{\eta}<\tau^{K}_0\right)\geq \mb{P}^{\hat{x}}\left(\overline{\tau}^K_{\eta}<\tau^{K}_{a'/K}\right) \geq 1-C_5\left(\frac{\norm{x}}{a'/K}\right)^{-p}\geq \exp\left(-\frac{2C_5a'^p}{(K\norm{x})^p}\right),
	\end{align}
	using that $1-h\geq e^{-2h}$ for all $0\leq h \leq 1/2$.
	Second, if $\norm{x}\leq a''/K$, then 
	\[
	\mb{P}^{\hat{x}}\left(\overline{\tau}^K_{\eta}<\tau^{K}_0\right)\geq \mb{P}^{\hat{x}}\left(\overline{\tau}^K_{a''/K}<\tau^{K}_0\right)\inf_{\hat{y}\in\hat{\mathcal{X}}^K_{a''/K}}\mb{P}^{\hat{y}}\left(\overline{\tau}^K_{\eta}<\tau^{K}_0\right).
	\]
	It follows from Lemma \ref{lemboundary} that one can bound from below the first factor in the right handside by a constant $c'>0$ independent of $x$. As for the second factor, it is greater than or equal to $1/2$, hence 
	\begin{align}\label{probaescape2}
	\mb{P}^{\hat{x}}\left(\overline{\tau}^K_{\eta}<\tau^{K}_0\right)\geq c'/2.
	\end{align}
	Let us choose $C''_1$ large enough so that $C''_1\geq 2C_5a'^p$ and $\exp(-C''_1/(a''^p))\leq c'/2$. Combing \eqref{probaescape1} and \eqref{probaescape2}, we obtain that for all $\hat{x}\in\hat{\mathcal{X}}^K_+$, 
	\begin{align}\label{probaescape3}
	\mb{P}^{\hat{x}}\left(\overline{\tau}^K_{\eta}<\tau^{K}_0\right)\geq \exp\left(-\frac{C''_1}{(K\norm{x})^p}\right).
	\end{align}
	We can conclude by plugging this into \eqref{probasurvival}. Using that $\norm{x}\leq d$ for all $x\in\mathcal{X}$ and setting $C'_1=C''_1+C_1(da')^p$, we obtain, for all $K\geq K_0$, all $\hat{x}=(x,\xi)\in\hat{\mathcal{X}}^K_+$ and all $t\geq 0$, 
	\[
	\mb{P}^{\hat{x}}\left(\tau^K_0>t\right)\geq \exp\left(-\frac{C'_1}{(K\norm{x})^p}\right)\exp\left(-\frac{C_1t}{K^p}\right).
	\]
	Finally, by integrating over time, we obtain, since $K\norm{x}\geq 1$, 
	\[
	\mb{E}^{\hat{x}}\left(\tau^K_0\right)=\int_{\mb{R}_+} \mb{P}^{\hat{x}}\left(\tau^K_0>t\right)\mr{d}t\geq e^{-C'_1}\int_{\mb{R}_+}\exp\left(-\frac{C_1t}{K^p}\right)\mr{d}t=e^{-C'_1}\frac{K^p}{C_1},
	\]
	which ends the proof. \hfill $\square$

		\subsection{Proof of Theorem \ref{majtimeextpers}}\label{sec:majtimeextpers}
	
		\begin{Pro}\label{prolyapinvXK}
	Assume that $\Lambda>0$ and $0\in\Gamma(Y)$, and let $p>p^*$. Let $T>0$ and $\eta>0$ be given by item \ref{itemprolyapinv} of Proposition \ref{prolyapPDMP} applied to $-p$. There exist $\tilde{\gamma}>1$, $a>0$ and $K_0\in\mb{N}^*$ such that for all $K\geq K_0$ and all $\hat{x}=(x,\xi)\in\hat{\mathcal{X}}_+^K$ satisfying $a/K\leq \norm{x}\leq \eta$,
	\begin{align}\label{lyapinv}
		\mb{E}^{\hat{x}}\left[\norm{X^K_T}^{-p}\mb{1}_{\left\{T<\tau_0^K\right\}}\right]\geq \tilde{\gamma}\, \norm{x}^{-p}.
	\end{align}
\end{Pro}
\prf The proof is similar to the proof of Proposition \ref{prolyapXK}.
Let $p\in (p^*,+\infty)$, and let $T,\tilde{\gamma}_0,\eta$ be given by item \ref{itemprolyapinv} of Proposition \ref{prolyapPDMP} applied to $-p$. Let $C_T,\varepsilon_T$ be given by Lemma \ref{lemcoupl}, set $M=2e^{C_FT}$ and $m=e^{-C_FT}/2$ and let $C'$ and $K_0$ be given by Lemma \ref{lemtauq}. Let $K\geq K_0$ and $\hat{x}=(x,\xi)\in\hat{\mathcal{X}}^K_+$ such that $\norm{x}\leq \eta$. For all $\varepsilon \in (0,\varepsilon_T]$, we have, combining Proposition \ref{prolyapPDMP}, Lemma \ref{lemcoupl} and Lemma \ref{lemtauq}, 
\begin{align*}
	\mb{E}^{\hat{x}}\left[\norm{X^K_T}^{-p}\mb{1}_{\left\{T<\tau_0^K\right\}}\right] 
	&\geq\mb{E}^{\hat{x}}\left[\norm{X^K_T}^{-p}\mb{1}_{\left\{T<\sigma^K_{\varepsilon}\right\}}\right]\\
	&\geq (1+\varepsilon)^{-p}\mb{E}^{\hat{x}}\left[\norm{X_T}^{-p}\mb{1}_{\left\{T<\sigma^K_{\varepsilon}\right\}}\right]\\
	&= (1+\varepsilon)^{-p}\left(\mb{E}^{\hat{x}}\left[\norm{X_T}^{-p}\right]-\mb{E}^{\hat{x}}\left[\norm{X_T}^{-p}\mb{1}_{\left\{\sigma^K_{\varepsilon}\leq T\right\}}\right]\right)\\	 
	&\geq (1+\varepsilon)^{-p}\left(\tilde{\gamma}_0 \norm{x}^{-p}-(\norm{x}e^{-\norm{G}_{\infty}T})^{-p}\mb{P}^{\hat{x}}\left[\sigma^K_{\varepsilon}\leq T\right]\right)\\
	&\geq \norm{x}^{-p}\Big((1+\varepsilon)^{-p}\tilde{\gamma}_0-2de^{p\norm{G}_{\infty}T}e^{-C_T K\norm{x}\varepsilon^2}
	\\&\qquad \qquad -e^{p\norm{G}_{\infty}T}T|E|\sup_{\xi_1\neq \xi_2,\ \norm{z-y}\leq d\varepsilon}\left|q(z,\xi_1,\xi_2)-q(y,\xi_1,\xi_2)\right|\!\Big)
\end{align*}
Since the functions $q(\cdot,\xi_1,\xi_2)$ are uniformly continuous, we can choose $\varepsilon>0$ small enough so that
\[
\tilde{\gamma}_1:=(1+\varepsilon)^{-p}\tilde{\gamma}_0- e^{p\norm{G}_{\infty}T}T|E|\sup_{\xi_1\neq \xi_2,\ \norm{z-y}\leq d\varepsilon}\left|q(z,\xi_1,\xi_2)-q(y,\xi_1,\xi_2)\right|>1,
\]
and then, we can choose $a>0$ large enough so that 
\[
\tilde{\gamma}:=\tilde{\gamma}_1-2de^{p\norm{G}_{\infty}T}e^{-C_T a\varepsilon^2}>1.
\]
We can conclude : for all $K\geq K_0$ and all $\hat{x}=(x,\xi)\in\hat{\mathcal{X}}^K_+$ such that $a/K\leq \norm{x}\leq \eta$, 
\[
\mb{E}^{\hat{x}}\left[\norm{X^K_T}^{-p}\mb{1}_{\left\{T<\tau_0^K\right\}}\right]\geq \tilde{\gamma} \norm{x}^{-p}.
\]	\hfill $\square$

We now turn to the proof of Theorem \ref{majtimeextpers}.
Let $p>p^*$. Let $T,\eta,\tilde{\gamma},a,K_0$ be given by Proposition \ref{prolyapinvXK}. Let $K\geq K_0$ and set, for all $n\in\mb{N}$,
\[
Z^K_n=\gamma^{-n}\norm{X^K_{nT}}^{-p}\mb{1}_{\left\{nT<\tau_0^K\right\}},
\]
and $\scr{T}^K=\left\lceil (\overline{\tau}^{K}_{\eta}\wedge\underline{\tau}^K_{a/K})/T\right\rceil$.   Let $\hat{x}=(x,\xi)\in \mathcal{X}^K_+$ such that $a/K\leq \norm{x}\leq \eta$. Using \eqref{lyapinv}, we see that under $\mb{P}^{\hat{x}}$ the process $(Z^K_{n\wedge \scr{T}^K})_{n\in\mb{N}}$ is a bounded submartingale with respect to the filtration $(\scr{F}_{nT})_{n\in\mb{N}}$ (see the proof of Lemma \ref{hittingtime1} for the details of a similar argument). Moreover, $\mb{P}^{\hat{x}}\left(\scr{T}^K<+\infty\right)=1$ due to the irreductibility of the chain $\hat{X}^K$ on $\hat{\mathcal{X}}^K_+$. Consequently, the optional sampling theorem yields
\begin{align}\label{optsamp}
	\mb{E}^{\hat{x}}\left[Z^K_{\scr{T}^K}\right]\geq \mb{E}\left[Z^K_0\right]=\norm{x}^{-p}.
\end{align}
Let $m$ be defined as in Lemma \ref{lemtauq}, and let $C'$ be given by the same lemma. We have
\begin{align}
	&\quad\mb{E}^{\hat{x}}\left[Z^K_{\scr{T}^K}\right]\nonumber\\
	&=\mb{E}\left[Z^K_{\scr{T}^K}\mb{1}_{\left\{\underline{\tau}^K_{a/K}<\overline{\tau}^{K}_{\eta}\right\}}\right]+\mb{E}\left[Z^K_{\scr{T}^K}\mb{1}_{\left\{\overline{\tau}^{K}_{\eta}<\underline{\tau}^K_{a/K}\right\}}\right] \nonumber \\
	&\leq K^p \mb{E}^{\hat{x}}\left[\tilde{\gamma}^{-\underline{\tau}^K_{a/K}/T}\right]+\mb{E}^{\hat{x}}\left[\mb{1}_{\left\{\overline{\tau}^{K}_{\eta}<\underline{\tau}^K_{a/K}\right\}}\sup_{0\leq s<T}\left(\norm{X^K_{\overline{\tau}^{K}_{\eta}+s}}^{-p}\mb{1}_{\left\{\overline{\tau}^{K}_{\eta}+s<\tau_0^K\right\}}\right)\right] \nonumber\\
	& \leq K^p \mb{E}^{\hat{x}}\!\left[\tilde{\gamma}^{-\underline{\tau}^K_{a/K}/T}\right]\!+\!\sup_{\substack{\hat{y}=(y,\zeta)\in\hat{\mathcal{X}}^K\\ \norm{y}\geq \eta}}\!\mb{E}^{\hat{y}}\!\left[\sup_{0\leq s<T}\!\left(\norm{X^K_{s}}^{-p}\mb{1}_{\left\{s<\tau_0^K\right\}}\right)\!\!\left(\!\mb{1}_{\left\{T<\underline{\tau}^K_{m\norm{y}}\right\}}\!+\!\mb{1}_{\left\{\underline{\tau}^K_{m\norm{y}}\leq T\right\}}\!\right)\!\right] \nonumber\\
	& \leq K^p \mb{E}^{\hat{x}}\left[\tilde{\gamma}^{-\underline{\tau}^K_{a/K}/T}\right]+(m\eta)^{-p}+2K^{p}e^{-C'K\eta}, \label{leqEZT}
\end{align}
using the strong Markov property of $X^K$ at time $\overline{\tau}^{K}_{\eta}$ and Lemma \ref{lemtauq}. Let $K'_0\geq K_0$ such that $\sup_{K\geq K_0}(2K^pe^{-C'K\eta})\leq (m\eta)^{-p}$, and fix $\eta'>0$ small enough so that $c_1:=\eta'^{-p}-2(m\eta)^{-p}>0$. Let us suppose that $K\geq K'_0$ and $\norm{x}\leq \eta'$. Combining \eqref{leqEZT} with \eqref{optsamp}, we obtain that 
\[
\mb{E}^{\hat{x}}\left[\tilde{\gamma}^{-\underline{\tau}^K_{a/K}/T}\right]\geq c_1K^{-p}.
\]
This enables us to bound stochastically $\underline{\tau}^K_{a/K}$ from below. For all $t>0$, using that $\tilde{\gamma}^{-\underline{\tau}^K_{a/K}/T}\leq \mb{1}_{\left\{\underline{\tau}^K_{a/K}\leq t\right\}}+\tilde{\gamma}^{-t/T}$ yields
\[
\mb{P}^{\hat{x}}\left[\underline{\tau}^K_{a/K}\leq t\right]\geq c_1K^{-p}-\tilde{\gamma}^{-t/T}.
\]
Set $t(K)=\left(T/\log(\tilde{\gamma})\right)(p\log(K)-\log(c_1)+\log(2))_+$, so that $\tilde{\gamma}^{-t(K)/T}\leq c_1K^{-p}/2$. We obtain that for all $K$ large enough and  $\hat{x}=(x,\xi)\in\hat{\mathcal{X}}^K_+$ with $\norm{x}\leq \eta'$,
\begin{align}\label{taua}
\mb{P}^{\hat{x}}\left[\underline{\tau}^K_{a/K}\leq t(K)\right] \geq c_1K^{-p}/2.
\end{align}
Moreover we claim that there exists constants $c_2,c_3,T'>0$ such that, for $K$ large enough : all $\hat{x}=(x,\xi)\in\hat{\mathcal{X}}^K_+$ : 
\begin{enumerate}[label= (\alph*)]
	\item \label{tau0} for all $\hat{x}=(x,\xi)\in\hat{\mathcal{X}}^K_+$ such that $\norm{x}\leq a/K$, we have $\mb{P}^{\hat{x}}\left[\tau_0^K\leq 1\right]\geq c_2$;
	\item \label{taueta'} for all $\hat{x}\in\hat{\mathcal{X}}^K_+$ such that $\norm{x}>\eta'$, $\mb{P}^{\hat{x}}\left[\underline{\tau}^K_{\eta'}\leq T'\right]\geq c_3$.
\end{enumerate}
Item \ref{tau0} is a consequence of Lemma \ref{lemboundary}. As for \ref{taueta'}, by Lemma \ref{lem:accessLineartoNonLinear} we have $0\in \Gamma(X)$, hence \cite[Proposition 3.14]{BLMZ} entails that there exists $c',T'>0$ such that, for all $\hat{x}\in\hat{\mathcal{X}}^K_+$ with $\norm{x}>\eta'$, 
\[
\mb{P}^{\hat{x}}\left(\exists t\in[0,T'],\ \norm{X_t}<\eta'/2\right)\geq c'.
\]
Let $C_{T'},\varepsilon_{T'}>0$ be given by Lemma \ref{lemcoupl}, and choose $\varepsilon\in (0,\varepsilon_{T'}\wedge(\eta'/2d)]$ small enough so that
\[
T'|E|\sup_{\xi_1\neq \xi_2,\ \norm{z-y}\leq d\varepsilon}{\left|q(z,\xi_1,\xi_2)-q(y,\xi_1,\xi_2)\right|}\leq c'/3.
\]
Then, using Equation \eqref{ineqcoupl3} we obtain that, for all $K$ large enough and all $\hat{x}=(x,\xi)\in\hat{\mathcal{X}}^K_+$ such that $\norm{x}>\eta'$, 
\begin{align*}
\mb{P}^{\hat{x}}\!\left[\underline{\tau}^K_{\eta'}\leq T'\right]\!\geq \mb{P}^{\hat{x}}\!\left(\exists t\!\in\![0,T'],\ \norm{X_t}<\frac{\eta'}{2}\right)- \mb{P}^{\hat{x}}\!\left[\sigma^K_{\varepsilon}\leq T'\right]\!\geq \frac{2c'}{3}- 2de^{-C_{T'}K\eta'\varepsilon^2}\geq \frac{c'}{3}, 
\end{align*}
which ends the proof of \ref{taueta'} with $c_3=c'/3$.

Combining \ref{tau0}, \ref{taueta'} and \eqref{taua}, we obtain that, for all $K$ large enough and $\hat{x}\in\hat{\mathcal{X}}^K_+$, 
\[
\mb{P}^{\hat{x}}\left[\tau_0^K\leq t(K)+T'+1\right]\geq c_4K^{-p}
\]
where $c_4=c_1c_2c_3/2$. Let us assume that $K$ is large enough so that $t(K)\leq c_5 \log(K)$ with $c_5=T(p+1)/\log(\tilde{\gamma})$. Then, for all $n\in\mb{N}$, the Markov property at time $nc_5\log(K)$ entails, for all $\hat{x}\in\hat{\mathcal{X}}^K_+$, 
\begin{align*}
\mb{P}^{\hat{x}}\left[\tau_0^K>(n+1)c_5\log(K)\right]&\leq \mb{P}^{\hat{x}}\left[\tau_0^K>nc_5\log(K)\right]\sup_{\hat{y}\in\hat{\mathcal{X}}^K_+}\mb{P}^{\hat{y}}\left[\tau_0^K>c_5\log(K)\right]\\ &\leq \mb{P}^{\hat{x}}\left[\tau_0^K>nc_5\log(K)\right]\left(1-c_4K^{-p}\right).
\end{align*}
By induction we deduce that, for all $K$ large enough, all $\hat{x}\in\hat{\mathcal{X}}^K_+$ and all $t\geq 0$, 
\[
\mb{P}^{\hat{x}}\left(\tau^K_0>t\right)\leq \left(1-c_4K^{-p}\right)^{\left \lfloor t/\left(c_5\log(K)\right)\right\rfloor}\leq \exp\left(-\frac{c_4}{K^p}\left(\frac{t}{c_5\log(K)}-1\right)\right).
\]
Since this holds for all $p>p^*$, we easily deduce Theorem \ref{majtimeextpers}. The majoration of $\mb{E}^{\hat{x}}\left(\tau^K_0\right)$ follows by integrating the upper bound on the survival function over time.
\hfill $\square$

\subsection{Proof of Theorem \ref{majtimeextnonpers}}\label{sec:majtimeextnonpers}

The proof is similar to the proof of Theorem \ref{majtimeextpers}, but is based on the use of a supermartingale instead of a submartingale. Since $g'(0)=\Lambda<0$ (see Theorem \ref{plyap}), we can choose $p>0$ such that $g(p)<0$. Then, let $T,\eta>0$ be given by Proposition \ref{prolyapPDMP} and $a>0$, $\gamma\in(0,1)$, $K_0\in\mb{N}^*$ be given by Proposition \ref{prolyapXK}, for this choice of $p$. Let $a'>a$, which remains to be fixed later in the proof. Let $K\geq K_0$, set 
\[
Z^K_n=\gamma^{-n}\norm{X^K_{nT}}^{p}\mb{1}_{\left\{nT<\tau_0^K\right\}},
\]
for all	$n\in\mb{N}$, and set $\scr{T}^K=\left\lceil( \overline{\tau}^K_{\eta}\wedge\underline{\tau}^K_{a'/K})/T\right\rceil$. Let $\hat{x}=(x,\xi)\in \mathcal{X}^K_+$ be such that $a'/K\leq \norm{x}\leq \eta$. Using \eqref{lyap}, we see that under $\mb{P}^{\hat{x}}$ the process $(Z^K_{n\wedge \scr{T}^K})_{n\in\mb{N}}$ is a supermartingale with respect to the filtration $(\scr{F}_{nT})_{n\in\mb{N}}$. 	
Hence, the optional sampling theorem yields
\begin{align}\label{optsampbis}
	\mb{E}^{\hat{x}}\left[Z^K_{\scr{T}^K}\right]\leq \mb{E}^{\hat{x}}\left[Z^K_0\right]=\norm{x}^{p}.
\end{align}
First, we can use this inequality to get a lower bound on the probability that $\underline{\tau}^K_{a'/K}<\overline{\tau}^{K}_{\eta}$. Let $m$ be defined as in Lemma \ref{lemtauq}, and let $C'$ be given by the same lemma. Using the strong Markov property of the process $X^K$ at time $\overline{\tau}^{K}_{\eta}$ and Lemma \ref{lemtauq} we get
\begin{align}
	\mb{E}^{\hat{x}}\left[Z^K_{\scr{T}^K}\right] \nonumber
	&\geq \mb{E}^{\hat{x}}\left[\mb{1}_{\left\{\overline{\tau}^{K}_{\eta}<\underline{\tau}^K_{a'/K}\right\}}\inf_{0\leq t <T}\left(\norm{X^K_{\overline{\tau}^{K}_{\eta}+t}}^{p}\right)\mb{1}_{\left\{\overline{\tau}^{K}_{\eta}+T<\tau_0^K\right\}}\right] \nonumber\\
	&\geq \mb{P}^{\hat{x}}\left(\overline{\tau}^{K}_{\eta}<\underline{\tau}^K_{a'/K}\right)\inf_{(y,\zeta)\in \hat{\mathcal{X}}^K,\,\norm{y}\geq \eta}\mb{E}^{(y,\zeta)}\left[\inf_{0\leq t <T}\norm{X^K_t}^{p}\mb{1}_{\left\{T<\tau_0^K\right\}}\right] \nonumber \\
	&\geq \mb{P}^{\hat{x}}\left(\overline{\tau}^{K}_{\eta}<\underline{\tau}^K_{a'/K}\right)(m\eta)^{p}\left(1-2e^{C'K\eta}\right). \label{boundtaueta1}
\end{align}
We may assume that $K_0$ is large enough so that $1-2e^{C'K_0\eta}\geq 1/2$. Combining \eqref{boundtaueta1} with \eqref{optsampbis}, we obtain
\begin{align}\label{boundtaueta2}
	\mb{P}^{\hat{x}}\left(\overline{\tau}^{K}_{\eta}<\underline{\tau}^K_{a'/K}\right)\leq 2 (m\eta)^{-p}\norm{x}^{p}.
\end{align}
Then, we can use \eqref{optsampbis} again to bound stochastically $\underline{\tau}^K_{a'/K}$ from below. We may assume that $K_0$ is large enough so that $\min_{i\in\intbk{1,d}}K_i(K)/K\geq \underline{\alpha}/2$ for all $K\geq K_0$, which implies that $\norm{X^K_{\underline{\tau}^K_{a'/K}}}\geq (a'-2/\underline{\alpha})/K$. We now fix the value of $a'$, large enough so that $a'':=a'-2/\underline{\alpha}>0$ and $1-2e^{-C'a''}\geq 1/2$. Using \eqref{optsampbis} we obtain
\begin{align*}
	\norm{x}^{p} 
	\geq \mb{E}^{\hat{x}}\left[Z^K_{\scr{T}^K}\mb{1}_{\left\{\underline{\tau}^K_{a'/K}<\overline{\tau}^{K}_{\eta}\right\}}\right] \geq \mb{E}^{\hat{x}}\left[\gamma^{- \underline{\tau}^K_{a'/K}/T}\mb{1}_{\left\{\underline{\tau}^K_{a'/K}<\overline{\tau}^{K}_{\eta}\right\}}\right]\left(ma''/K\right)^{p}(1-2e^{-C'a''}),
\end{align*}
hence  
\[
\mb{E}^{\hat{x}}\left[\gamma^{- \underline{\tau}^K_{a'/K}/T}\mb{1}_{\left\{\underline{\tau}^K_{a'/K}<\overline{\tau}^{K}_{\eta}\right\}}\right]\leq CK^{p}
\]
with $C=2 (ma'')^{-p}d^{p}$.
Combining this bound with \eqref{boundtaueta2} and using that \[\mb{P}^{\hat{x}}\left(\left\{\underline{\tau}^K_{a'/K}<\overline{\tau}^{K}_{\eta}\right\}\cup\left\{\overline{\tau}^{K}_{\eta}<\underline{\tau}^K_{a'/K}\right\}\right)=1\] due to the irreducibility of $X^K$ on $\hat{\mathcal{X}}^K_+$ yields, for all $t\geq 0$, 
\begin{align*}
	\mb{P}^{\hat{x}}\left(\underline{\tau}^K_{a'/K}>t\right)&\leq \mb{P}^{\hat{x}}\left(\underline{\tau}^K_{a'/K}>t,\, \underline{\tau}^K_{a'/K}<\overline{\tau}^{K}_{\eta}\right)+\mb{P}^{\hat{x}}\left(\overline{\tau}^{K}_{\eta}<\underline{\tau}^K_{a'/K}\right) \\
	& \leq CK^{p}\gamma^{t/T}+ 2 (m\eta)^{-p}\norm{x}^{p}.
\end{align*}
Let $\eta'\in(0,\eta)$ be small enough so that $2 (m\eta)^{-p}\eta'^{p}\leq 1/4$, and set \[t(K)=(Tp/\log(\gamma^{-1}))\log(4CK),\] so that $CK^{p}\gamma^{t(K)/T}=1/4$. Then, the above inequality yields, for all $\hat{x}=(x,\xi)\in\hat{\mathcal{X}}^K_+$ such that $a''/K\leq \norm{x}\leq \eta'$ :
\[
\mb{P}^{\hat{x}}\left(\underline{\tau}^K_{a'/K}\leq t(K)\right)\geq 1/2.
\]
The end of the proof is similar to the one of Theorem \ref{majtimeextpers}, after obtaining \eqref{taua} : the same arguments entail that there exists $C_3,C'_3>0$ such that for $K$ large enough and for all $\hat{x}\in\hat{\mathcal{X}}^K_+$, 
\[
\mb{P}^{\hat{x}}\left(\tau^K_0>t\right)\leq C'_3\exp\left(-\frac{C_3t}{\log(K)}\right).
\]
The fact that $0\in\Gamma(X)$ is justified by Lemma \ref{lem:accessLineartoNonLinear}. 
\hfill $\square$

	\subsection{Proof of Theorem \ref{thm:qsdpers}}\label{sec:prfqsdpers}
	
	As explained in Section \ref{sec:limitQSD}, the first main ingredient in the proof of Theorem \ref{thm:qsdpers} is given by Proposition \ref{lemphitheta:1}, which we state here again.
	
	\begin{Pro}\label{lemphitheta}
		Assume that $\Lambda>0$ and let $p\in (0,p^*)$. Let $T,a>0$  be given by Proposition \ref{prolyapXK} applied to $-p$, and for all $K\geq d$, let $\varphi^K:\hat{\mathcal{X}}^K_+\to \mb{R}_+^*$ be defined by $\varphi^K(x,\xi)=\norm{x}^{-p}\wedge(a/K)^{-p}$. There exist $\theta\in (0,1)$ and $C_6>0$ such that :
		
		\begin{enumerate}[label=\roman*)]
			\item \label{itemphitheta1}for all $K$ large enough,
			\begin{align}\label{phitheta1}
				\tilde{P}^K_T\varphi^K\leq \theta_1 \varphi^K + C_6\ ;
			\end{align}
			\item \label{tightqsd} $\limsup_{K\rightarrow +\infty}\mu^K\varphi_1^K\leq C_6/(1-\theta_1)$.
		\end{enumerate}
	\end{Pro}
	
	\prf We start by the proof of \ref{itemphitheta1}. Let $p\in (0,p^*)$. Let $T, \eta, \gamma, a$ be given by Proposition \ref{prolyapXK} applied to $-p$. Set $M=2e^{C_FT}$ and $m=e^{-C_FT}/2$, and let $C'$ be given by Lemma \ref{lemtauq}. Let $K\geq d$ and define $\varphi^K:\hat{\mathcal{X}}^K_+\to \mb{R}_+^*$ by $\varphi^K(x,\xi)=\norm{x}^{-p}\wedge(a/K)^{-p}$. Let $\hat{x}=(x,\xi)\in\hat{\mathcal{X}}^K_+$. There are three cases to consider. In the following, "$K$ large enough" means $K$ greater than some $K_0$ independent of $\hat{x}$.
	
	First, suppose that $a/K\leq \norm{x} \leq \eta$. Then, using Proposition \ref{prolyapXK} and Lemma \ref{lemtauq}, we get, for $K$ large enough,
	
	\begin{align}
		\tilde{P}^K_T\varphi^K(\hat{x})&=\mb{E}^{\hat{x}}\left[\left(\norm{X^K_T}^{-p}\wedge(a/K)^{-p}\right)\mb{1}_{\left\{T<\tau_0^K\right\}}\right] \nonumber \\
		&\leq \mb{E}^{\hat{x}}\left[\norm{X^K_{T}}^{-p}\mb{1}_{\left\{T<\tau_0^K\right\}}\right]\nonumber\\ &\leq \gamma \varphi^K(\hat{x}) \label{ineqmiddle}.
	\end{align}	
	
	Second, consider the case $\norm{x}>\eta$. We have
	\begin{align}
		\tilde{P}^K_T\varphi^K(\hat{x})&\leq \mb{P}^{\hat{x}}\left(\underline{\tau}^K_{m\eta}\leq T\right)(a/K)^{-p}+\mb{P}^{\hat{x}}\left(\underline{\tau}^K_{m\eta}> T\right)(m\eta)^{-p} \nonumber \\ &\leq 2e^{-C'K\eta}(a/K)^{-p}+(m\eta)^{-p} \nonumber \\ 
		&\leq C_6 \label{ineqfar}
	\end{align}
	where $C_6:=2a^{-p}\sup_{n\in\mb{N}^*}\left(n^pe^{-C'n\eta}\right)+(m\eta)^{-p}<\infty$.
	
	Finally, we treat the case $\norm{x}<a/K$. Letting $c>0$ be given by Lemma \ref{lemboundary}, we have, for $K$ large enough,
	\begin{align}\label{ineqboundary}
		\tilde{P}^K_T\varphi^K(\hat{x})\leq \mb{P}^{\hat{x}}\left(T<\tau_0^K\right)(a/K)^{-p}\leq (1-c)\varphi^K(\hat{x}).
	\end{align}	
	
	Thus we can end the proof of \ref{itemphitheta1} by combining \eqref{ineqmiddle}, \eqref{ineqfar} and \eqref{ineqboundary} : for $K$ large enough,
	\[
	\tilde{P}^K_T\varphi^K\leq \theta_1 \varphi^K + C_6
	\]
	with $\theta_1:=\gamma\vee (1-c)<1$.	
	
	Finally, let us prove \ref{tightqsd}.
	Integrating \eqref{phitheta1} against $\mu^K$ yields
	\[
	e^{-\lambda^K T}\mu^K\varphi_1^K=\mu^K\tilde{P}^K_{T}\varphi^K\leq \theta_1\mu^K\varphi^K+C_6.
	\]
	Now, by Proposition \ref{pro:OKp} we have $\lambda^K\rightarrow 0$ hence 
	\[
	\limsup_{K\rightarrow +\infty}\mu^K\varphi^K\leq \limsup_{K\rightarrow +\infty} \frac{C_6}{e^{-\lambda^K T}-\theta_1}=\frac{C_6}{1-\theta_1}.
	\]	\hfill $\square$
	
	\

	The other important result for the proof of Theorem \ref{thm:qsdpers} is the following. It is a consequence of \cite[Lemma 6.3]{S19}, but for completeness we give the proof here.
	
	\begin{Pro}
		\label{pro:lambdapositifstationary}
		If $\Lambda > 0$, every limit point of $(\mu^K)_{K \geq d}$ is a stationary distribution of $\hat{X}$.
	\end{Pro}
	
\prf Let $f\in\mathcal{C}(\hat{\mathcal{X}})$ and let $\varepsilon > 0$. By compacity of $\hat{\mathcal{X}}$, there exists $\delta>0$ such that for all $x,y\in\mathcal{X}$ and $\xi\in E$, $\| x - y \| \leq \delta$ implies that $\| f(x,\xi) - f(y,\xi) \| \leq \varepsilon$. Let $(P_t)_{t\in\mb{R}_+}=(P^{\hat{X}}_t)_{t\in\mb{R}_+}$ and  $(P^K_t)_{t\in\mb{R}_+}$ denote the semi-group associated to $\hat{X}$ and $\hat X^K$ respectively. Recall that $(\tilde{P}^K_t)_{t\in\mb{R}_+}$ denotes the killed semi-group of $\hat{X}^K$. Fix $t\in\mb{R}_+$. For all $K\geq d$, we have
\[
| \mu^K P_t f - \mu^K f|\leq | \mu^K P_t f -  \mu^K P_t^K f|  + | \mu^K P_t^K f - \mu^K f |.
\]
The second term of the right handside satisfies
\begin{align*}
| \mu^K P_t^K f - \mu^K f |&\leq | \mu^K P_t^K f - \mu^K \tilde{P}_t^K f  |+| \mu^K \tilde{P}_t^K f - \mu^K f |\\&\leq \left|\mb{E}^{\mu^K}\left[f(\hat{X}^K_t)(1-\mb{1}_{\{\tau^K_0>t\}})\right]\right|+ |(e^{-\lambda^K t}-1)\mu^Kf|\\&\leq 2\| f \|_{\infty} ( 1 - e^{ - \lambda^K t} ).
\end{align*}
As for the first term,  
\begin{align*}
| \mu^K P_t f -  \mu^K P_t^K f| &\leq \int_{\hat{\mathcal{X}}^K} \mb{E}^{\hat x}\left[| f(X_t, \Xi_t) - f(X_t^K,\Xi_t^K)|\right] \mu^K(\mr{d}\hat x) \\
&\leq 2 \|f\|_{\infty} \sup_{\hat x \in \hat{\mathcal{X}}^K} \mb{P}^{\hat x} ( \sigma^K_{\delta}\wedge \eta^K \leq t) + \varepsilon
\end{align*}
where $\sigma^K_{\delta}= \inf \{ s \geq 0 \: : \|X^K_s - X_s \| > \delta\}$ and $\eta^K=\inf\{s\geq 0 : \Xi^K_s\neq \Xi_s\}$. Now, by Proposition \ref{pro:OKp} and Proposition \ref{LLN} respectively we have
\[
 \lambda^K\underset{K\rightarrow +\infty}{\rightarrow}0 \quad \text{and}\quad \sup_{\hat x \in \hat{\mathcal{X}}^K} \mb{P}^{\hat x} ( \sigma^K_{\delta}\wedge \eta^K \leq t) \underset{K\rightarrow +\infty}{\rightarrow} 0.
\]
This implies that $\limsup_{K\rightarrow +\infty}|\mu^KP_t f - \mu^K f| \leq  \varepsilon$. Since this holds for all $\varepsilon>0$, we obtain that $|\mu^KP_t f - \mu^K f|\rightarrow 0$ as $K\rightarrow +\infty$.

By the Feller property of $(P_s)_{s\in\mb{R}_+}$, we have $P_tf\in\mathcal{C}(\hat{X})$. Hence, for every weak limit point $\mu$ of $(\mu^K)_{K\geq d}$, we have $\mu P_t f = \mu f$. Since this holds for all $t\in\mb{R}_+$ and $f\in\mathcal{C}(\hat{X})$, $\mu$ is stationary for $(P_t)_{t \in\mb{R}_+}$. \hfill $\square$

\

Now let us prove Theorem \ref{thm:qsdpers}. Since $\Lambda>0$, $p^*$ is positive and we let $p\in(0,p^*)$. By Proposition \ref{lemphitheta}, \ref{tightqsd} we have 
\[
\limsup_{K\rightarrow +\infty}\mu^K\varphi_1^K\leq C_6/(1-\theta_1)
\]
for some $\theta_1\in(0,1)$ and $C_6>0$, where for all $K\geq d$,  $\varphi^K:\hat{\mathcal{X}_+}\to \mb{R}_+^*$ is defined by $\varphi^K(x,\xi)=\norm{x}^{-p}\wedge(a/K)^{-p}$, for some $a>0$. Therefore, for all $\varepsilon\in(0,d]$, we have 
\[
\limsup_{K\rightarrow +\infty}\mu^K\{(x,\xi)\in\hat{\mathcal{X}}_+ : \norm{x}< \varepsilon\}\leq \limsup_{K\rightarrow +\infty}\mu^K\left(\frac{\varphi^K}{\varepsilon^{-p}\wedge(a/K)^{-p}}\right)\leq \varepsilon^p C_6/(1-\theta_1).
\]
Since the right handside vanishes as $\varepsilon\rightarrow 0$, the sequence $(\mu^{K})_{K\geq d}$ is persistent, i.e. it is tight on $\hat{\mathcal{X}}_+$. 

Therefore, if we let $\mu\in \scr{L}$, then $\mu\in \mathcal{P}(\hat{\mathcal{X}}_+)$.  Moreover, Proposition \ref{pro:lambdapositifstationary} entails that $\mu$ is a stationary distribution of $\hat{X}$, which proves item \ref{lppers}. Let $K_1\geq d$. Given that $\varphi^{K_1}$ is continuous and bounded on $\hat{\mathcal{X}_+}$ and that $\varphi^{K}\geq \varphi^{K_1}$ for all $K\geq K_1$, we have
\[
\mu\varphi^{K_1}\leq \limsup_{K\rightarrow +\infty}\mu^K\varphi^{K_1}\leq \limsup_{K\rightarrow +\infty}\mu^K\varphi^K\leq C_6/(1-\theta).
\]
Furthermore, since $\varphi^{K_1}$ converges increasingly to $(x,\xi)\mapsto \norm{x}^{-p}$ as $K_1$ goes to infinity, we deduce by monotone convergence that 
\[
\int_{\hat{\mathcal{X}}_+}\norm{x}^{-p}\mu(\mr{d}x,\mr{d}\xi)\leq C_6/(1-\theta),
\]
which ends the proof of item \ref{itemqsdpers}.

Finally, let us prove item \ref{cvgqsd}. We assume that \ref{ass:monotonesubhomo} holds and that $x\mapsto Q(x)$ is constant. Then, by Theorem \ref{th:BS19gen} there exists a unique persistent stationary distribution $\mu^*$ of $\hat{X}$. Hence, $\scr{L}=\{\mu^*\}$, which implies that $\mu^K$ converges weakly to $\mu^*$ as $K$ goes to $+\infty$, by compacity of $\mathcal{P}(\hat{\mathcal{X}})$ for the weak topology. \hfill $\square$

\section{Appendix}
	\label{sec:appendix}
	
		\subsection{Trajectorial representations with Poisson random measures}\label{represtraj}
		
		In this section we justify the construction of the continuous-time Markov chains $\hat{X}^K$ and the PDMP $\hat{X}$ given in Section \ref{sec:construction}.
		
		\paragraph{Construction of $\hat{X}^K$.} Let us fix $K\geq d$ and $\hat{x}=(x,\xi)\in\hat{\mathcal{X}}$. The SDE \eqref{eqXK}-\eqref{eqIK} can be written compactly as 
		\begin{align}\label{rep}
		\hat{X}^{K,\hat{x}}_t=(\left\lfloor x\right\rfloor_K,\xi)+\int_{(0,t]\times \mb{R}_+\times Y}G^K(u,y,\hat{X}^{K,\hat{x}}_{s-})\scr{N}(\mr{d}s,\mr{d}u,\mr{d}y)
		\end{align}
 		where $Y=\left(\llbracket 1,d \rrbracket \times\left\{-1,1\right\}\right)\sqcup E$ and $G^K:\mb{R}_+\times Y\times \mathcal{X}^K\times E\to\mb{R}^d\times \mb{R}$ is defined by
 		\[
 		G^K(u,y,z,\zeta)=\begin{cases} \mb{1}_{\left\{u\leq K_i \beta_{he_i}(z,\zeta)\right\}}K_i^{-1}h(e_i,0)&\text{ if }y=(i,h)\in\intbk{1,d}\times \{-1,1\}\\  \mb{1}_{\left\{\xi'\neq \zeta\right\}}\mb{1}_{\left\{u\leq q\left(z,\zeta,\xi'\right)\right\}}(0,\xi'-\zeta)& \text{ if }y=\xi'\in E
 		\end{cases}.
 		\]
 		 		
 		Setting $\nu=\mr{Leb}\otimes \sum_{z\in\left(\llbracket 1,d \rrbracket \times\left\{-1,1\right\}\right)\sqcup E}\delta_{z}$, we can check that $\hat{z}+G^K(u,y,\hat{z})\in\mathcal{X}^K$ for all $(u,y,\hat{z})\in \mb{R}_+\times Y\times \hat{\mathcal{X}}^K$ and that $\nu\left(\left\{(u,y)\in \mb{R}_+\times Y: G^K(u,y,\hat{z})=\hat{z}'\right\}\right)=L^K(\hat{z},\hat{z}')$ for all distinct $\hat{z},\hat{z}'\in \hat{\mathcal{X}}^K$, where we recall that $L^K$ denotes the desired transition rate matrix for $\hat{X}^K$. Recall that $\nu$ is such that $\mr{Leb}\otimes \nu$ is the intensity measure of $\scr{N}$.

	Let $\dagger$ be a cemetery point. We can define inductively a sequence of random variables $(\hat{Z}^{K,\hat{x}}_n,T^{K,\hat{x}}_n)_{n\in\mb{N}}$ with values in $\left(\hat{\mathcal{X}}^K\times \mb{R}_+\right)\cup \{(\dagger,+\infty)\}$ by setting $(\hat{Z}^{K,\hat{x}}_0,T^{K,\hat{x}}_0)=(\left\lfloor x \right\rfloor_K,\xi,0)$ and then for all $n\in\mb{N}$, omitting the exponent $\hat{x}$ :
	

	\begin{enumerate}[label=$\bullet$]
		\item on $\left\{(\hat{Z}^K_{n},T^K_{n})\neq(\dagger,+\infty)\right\}$,
		\[
		T^K_{n+1}:=\inf\left\{t\geq T^K_n: \int_{(T^K_n,t]\times \mb{R}_+\times Y}G^K(u,y,\hat{Z}^K_n)\scr{N}(\mr{d}s,\mr{d}u,\mr{d}y)\neq 0\right\}
		\]
		and then 
		\begin{enumerate}[label=$-$]
			\item on $\left\{T^K_{n+1}<+\infty\right\}$,
			\begin{align*}
				\hat{Z}^K_{n+1}&:=\hat{Z}^K_n+\int_{(T^K_n,T^K_{n+1}]\times \mb{R}_+\times Y}G^K(u,y,\hat{Z}^K_n)\scr{N}(\mr{d}s,\mr{d}u,\mr{d}y) ;
			\end{align*}    		
			\item on $\left\{T^K_{n+1}=+\infty\right\}$, $\hat{Z}^K_{n+1}:=\dagger$ ;
		\end{enumerate}    		
		\item on $\left\{(\hat{Z}^K_{n},T^K_{n})=(\dagger,+\infty)\right\}$,
		$(\hat{Z}^K_{n+1},T^K_{n+1}):=(\dagger,+\infty).$
	\end{enumerate}
	For all $n\in\mb{N}$, $T^K_n$ is a $(\scr{F}^{\circ}_t)$-stopping time, and $\hat{Z}^K_n$ is $\scr{F}^{\circ}_{T^K_n}$-measurable. By strong Markov property of $\scr{N}$, conditional on $\{T^K_n<+\infty\}$ the random measure $\theta_{T^K_n}\scr{N}(\cdot):=\scr{N}\left(\left\{(T^K_n+s,u,y),\, (s,u,y)\in \cdot\right\}\right)$ is a Poisson random measure of the same intensity as $\scr{N}$ with respect to the filtration $(\scr{F}^{\circ}_{T^K_n+t})_{t\in\mb{R}_+}$, independent of $\scr{F}^{\circ}_{T^K_n}$. Hence, using the definition of $(\hat{Z}^K_{n+1},T^K_{n+1})$ we obtain, for all $t\in\mb{R}_+$ and $\hat{z}\in\hat{\mathcal{X}}^K$,
	\begin{align}\label{jumpholdZ}	 \mb{P}\left(\hat{Z}^K_{n+1}=\hat{z},T^K_{n+1}-T^K_n>t\,|\,\scr{F}^{\circ}_{T^K_n}\right)=\frac{L^K(\hat{Z}^K_n,\hat{z})}{\left|L^K(\hat{Z}^K_n,\hat{Z}^K_n)\right|}e^{-|L^K(\hat{Z}^K_n,\hat{Z}^K_n)|t}
	\end{align}
    almost surely on $\{T^K_n<\infty\}$. There is not division by zero because for all $\hat{z}=(z,\zeta)\in\hat{\mathcal{X}}^K$, $|L^K(\hat{z},\hat{z})|\geq |q(z,\zeta,\zeta)|>0$ by irreducibility of $Q(z)$. Equation \eqref{jumpholdZ} shows that the sequence $(Z^K_n,T^K_n)_{n\in\mb{N}}$ is distributed as the embedded chain of a continuous-time Markov chain of rate matrix $L^K$. In particular $\sup_{n\in\mb{N}}T^K_n=+\infty$ almost surely. Hence, setting $N^{K,\hat{x}}_t=\sup\left\{n\in\mb{N}: T^{K,\hat{x}}_n\leq t \right\}$, the process $(\hat{X}^{K,\hat{x}}_t)_{t\in\mb{R}_+}$ defined by
	\begin{align}\label{defXI}
		\hat{X}^{K,\hat{x}}_t=\hat{Z}^{K,\hat{x}}_{N^{K,\hat{x}}_t}
	\end{align}
	is a continuous-time Markov chain of rate matrix $L^K$, which is $(\scr{F}^{\circ}_t)$-adapted. By construction, it solves \eqref{rep}. It is straightforward to see that the embedded chain of any other solution must coincide with $(Z^{K,\hat{x}}_n,T^{K,\hat{x}}_n)_{n\in\mb{N}}$ almost surely, hence the uniqueness of the solution up to indistinguishability. Finally, the strong Markov property of $\hat{X}^{K,\hat{x}}$ with respect to $(\scr{F}^{\circ}_t)_{0\leq t \leq \infty}$ follows from the strong Markov property of $\scr{N}$ and the uniqueness of solutions of \eqref{rep}. 
	
	\paragraph{Contruction of $\hat{X}$.} The justification of the construction of the PDMP is very similar to the previous one. To begin with, modifying $\scr{N}_{\Xi}$ on a $\mb{P}$-negligible set if necessary, we may assume that for all $\omega^{\circ}\in\Omega^{\circ}$ and all $t\in\mb{R}_+$, $\scr{N}_{\Xi}(\omega^\circ)(\left\{t\right\}\times \mb{R}_+\times E)\leq 1$ and $\scr{N}(\omega^{\circ})((0,t]\times [0,\overline{q}]\times E)<+\infty$ where $\overline{q}=\sup_{(z,\zeta)\in\hat{\mathcal{X}}}|q(z,\zeta,\zeta)|<\infty$.

	Let $\hat{x}=(x,\xi)\in\hat{\mathcal{X}}$. We define inductively the sequence $(Z^{\hat{x}}_n,\tilde{\Xi}^{\hat{x}}_n,T^{\hat{x}}_n)_{n\in\mb{N}}$ with values in $(\mathcal{X}\times E\times\mb{R}_+)\cup\left\{(\dagger,\dagger,+\infty)\right\}$ by $(Z^{\hat{x}}_0,\tilde{\Xi}^{\hat{x}}_0,T^{\hat{x}}_0)=(x,\xi,0)$ and then, for all $n\in\mb{N}$, omitting the exponent $\hat{x}$ :
	

	\begin{enumerate}[label=$\bullet$, leftmargin=*]
		\item on $\left\{(Z_{n},\tilde{\Xi}_{n},T_{n})\neq(\dagger,\dagger,+\infty)\right\}$,
		\[
		T_{n+1}:=\inf\left\{t\geq T_n:\!\! \int_{(T_n,t]\times\mb{R}_+\times E}\!(\xi'-\tilde{\Xi}_n)\mb{1}_{\left\{u\leq q\left(\psi^{\tilde{\Xi}_n}_{s-T_n}(Z_n),\tilde{\Xi}_n,\xi'\right)\right\}}\scr{N}_{\Xi}(\mr{d}s,\mr{d}u,\mr{d}\xi')\neq 0\right\}
		\]
		and then
		\begin{enumerate}[label=$-$]
			\item on $\left\{T_{n+1}<+\infty\right\}$,
			\begin{align*}Z_{n+1}&:=\psi^{\tilde{\Xi}_n}_{T_{n+1}-T_{n}}(Z_n), \label{Xchinc}\\
				\tilde{\Xi}_{n+1}&:=\tilde{\Xi}_n+\int_{(T_n,T_{n+1}]\times\mb{R}_+\times E}(\xi'-\tilde{\Xi}_n)\mb{1}_{\left\{u\leq q\left(\psi^{\tilde{\Xi}_n}_{s-T_n}(Z_n),\tilde{\Xi}_n,\xi'\right)\right\}}\scr{N}_{\Xi}(\mr{d}s,\mr{d}u,\mr{d}\xi') ;
			\end{align*}    		
			\item on $\left\{T_{n+1}=+\infty\right\}$, $(Z_{n+1},\tilde{\Xi}_{n+1}):=(\dagger,\dagger)$ ;
		\end{enumerate}    		
		\item on $\left\{(Z_{n},\tilde{\Xi}_{n},T_{n})=(\dagger,\dagger,+\infty)\right\}$, 
		$
		(Z_{n+1},\tilde{\Xi}_{n+1},T_{n+1}):=(\dagger,\dagger,+\infty).
		$
	\end{enumerate}
	The sequence $(T_n)_{n\in\mb{N}}$ is an increasing sequence of $(\scr{F}^{\circ}_t)$-stopping times which tends to $+\infty$ on all $\Omega^{\circ}$ (thanks to the modification of $\scr{N}_{\Xi}$). Setting $N_t=\sup\left\{n\in\mb{N}: T_n\leq t \right\}$, the process defined by
\begin{align}
	(X^{\hat{x}}_t,\Xi^{\hat{x}}_t)=\left(\varphi^{\tilde{\Xi}_{N_t}}_{t-T_{N_t}}(Z_{N_t}),\tilde{\Xi}_{N_t}\right)
\end{align}
for all $t\in\mb{R}_+$ is càdlàg, $(\scr{F}^{\circ}_t)$-adapted, and solution of the SDE 
\begin{align}			
	X^{\hat{x}}_t&= x+\int_0^t F(\hat{X}^{\hat{x}}_s)\mr{d}s\label{eqXc}\\
	\Xi^{\hat{x}}_t&= \xi+\int_{(0,t]\times \mb{R}_+\times E}\mb{1}_{\left\{\xi'\neq \Xi^{\hat{x}}_{s-}\right\}}\left(\xi'-\Xi^{\hat{x}}_{s-}\right)\mb{1}_{\left\{u\leq q\left(\hat{X}^{\hat{x}}_{s-},\xi'\right)\right\}}\scr{N}_{\Xi}\left(\mr{d}s,\mr{d}u,\mr{d}\xi'\right).\label{eqIc}
\end{align}

Regarding the uniqueness, it is straightforward to see that if $(Y^{\hat{x}}_t,J^{\hat{x}}_t)_{t\in\mb{R}_+}$ is another càdlàg solution of this SDE, then its embedded chain $(\tilde{Y}_n,\tilde{J}_n,S_n)_{n\in\mb{N}}$ must coincide almost surely with $(Z_{n},\tilde{\Xi}_{n},T_{n})$. By embedded chain, we mean that $(S_n)_{n\in\mb{N}}$ is the sequence of jump times of $(J_t)_{t\in\mb{R}_+}$ and for all $n\in\mb{N}$, $(\tilde{Y}_n,\tilde{J}_n)$ equals $(Y_{S_n},J_{S_n})$ if $S_n<\infty$ and $(\dagger,\dagger)$ otherwise. Hence, the solution of \eqref{eqXc}-\eqref{eqIc} is unique (up to indistinguishability). 

The recurrence relations on the embedded chain can be written under the form 
\[
(\tilde{X}_{n+1},\tilde{I}_{n+1},T_{n+1})=H(\tilde{X}_{n+1},\tilde{I}_{n+1},T_{n+1},\theta_{T_n}\scr{N}),   	
\]
where $H$ is a measurable map from 
$\left((\mathcal{X}\times E\times\mb{R}_+)\cup\left\{(\dagger,\dagger,+\infty)\right\}\right)\times \Pi$ to $(\mathcal{X}\times E\times\mb{R}_+)\cup\left\{(\dagger,\dagger,+\infty)\right\}$, with $\Pi$ denoting the set of measures on $(\mb{R}_+^2\times E,\scr{B}(\mb{R}_+)^{\otimes 2}\otimes \scr{P}(E))$ taking values in $\mb{N}\cup\left\{+\infty\right\}$, equipped with the sigma-algebra $\scr{G}$ generated by the maps $\pi\mapsto \pi(A)$, $A\in \scr{B}(\mb{R}_+)^{\otimes 2}\otimes \scr{P}(E)$. Given that $(Z_0,\tilde{\Xi}_0,T_0)=(x,\xi,0)$, we deduce that we can write
\[
\hat{X}^{\hat{x}}=\Phi(\hat{x},\scr{N}_{\Xi})
\]
for some measurable map $\Phi:\hat{\mathcal{X}}\times \Pi\to \mathcal{D}(\mb{R}_+,\hat{\mathcal{X}})$. What's more, by construction, for all $\tau\in\mb{R}_+$ we have $(X^{\hat{x}}_{\tau+t})_{t\in\mb{R}_+}=\Phi(X^{\hat{x}}_{\tau},\theta_{\tau}\scr{N}_{\Xi})$.

It remains to show that $\hat{X}^{\hat{x}}$ is a switched dynamical system of local characteristics $((F^{\xi})_{\xi\in E},Q)$. By construction of the embedded chain, almost surely on $\{T_n<+\infty\}$ we have
	\begin{align*}
		\mb{P}\left(T_{n+1}-T_n>t\,|\, \scr{F}^{\circ}_{T_n}\right)=\exp\left(-\int_0^{t}\left|q\left(\psi^{\tilde{\Xi}_n}_s(Z_n),\tilde{\Xi}_n,\tilde{\Xi}_n\right)\right|\mr{d}s\right).
	\end{align*}
	Then, almost surely on $\{T_{n+1}<+\infty\}$ we have
    \begin{align*}Z_{n+1}=\psi^{\tilde{\Xi}_n}_{T_{n+1}-T_{n}}(Z_n)\quad \text{and}\quad\mb{P}\left(\tilde{\Xi}_{n+1}=\cdot\,|\, \scr{F}^{\circ}_{T_n},T_{n+1}\right)=\frac{q(Z_{n+1},\tilde{\Xi}_{n},\cdot)}{|q(Z_{n+1},\tilde{\Xi}_n,\tilde{\Xi}_n)|},
    \end{align*}
	while on $\left\{T_{n+1}=+\infty\right\}$, $(Z_{n+1},\tilde{\Xi}_{n+1})=(\dagger,\dagger)$. This implies that $\hat{X}^{\hat{x}}$ is a switched dynamical system of local characteristics $((F^{\xi})_{\xi\in E},Q)$, see the construction based on the embedded chain in \cite{Dav84}. Since for all finite $(\scr{F}^{\circ}_t)$-stopping time $\tau$ we have $(\hat{X}^{\hat{x}}_{\tau+t})_{t\in\mb{R}_+}=\Phi(\hat{X}^{\hat{x}}_{\tau},\theta_{\tau}\scr{N}_{\Xi})$, the strong Markov property of $\scr{N}$ implies the strong Markov property of $\hat{X}^{\hat{x}}$, with respect to the filtration $(\scr{F}^{\circ}_t)_{0\leq t\leq \infty}$.

	\subsection{A Chernoff bound for Poisson-driven martingales}
	
	\begin{Lem}\label{lemChernoff}
		Let $(U,\scr{B}(U))$ be a Polish space equipped with its Borel sigma-field, and let $\mu$ be a $\sigma$-finite measure on $(U,\scr{B}(U))$. Let $(\Omega,\scr{F},\left(\scr{F}_t)_{0\leq t \leq \infty},\mb{P}\right)$ be a filtered probability space satisfying the usual conditions, equipped with a Poisson random measure $\scr{N}$ on $\mb{R}_+\times U$ of intensity $\mr{Leb}\otimes \mu$. Let $\mathscr{P}$ denote the $(\scr{F}_t)$-predictable sigma-field on $\Omega\times \mb{R}_+$ and let $G:(\omega,s,u)\rightarrow G_{s,u}(\omega)$ be a measurable function from $(\Omega\times \mb{R}_+ \times U,\scr{P}\otimes \scr{B}(U))$ to $\left(\mb{R},\scr{B}(\mb{R})\right)$. Let $A,C,T>0$, let $B\in\scr{B}(U)$ such that $\mu(B)<\infty$, and assume that $\mb{P}$-almost surely,   
		\begin{align} \label{boundG}
		|G_{s,u}|\leq C\mb{1}_B(u)
		\end{align}
		for all $(s,u)\in [0,T]\times U$.
		Then, setting $\tilde{N}=N-\mr{Leb}\otimes \mu$ we have	
		\begin{align}\label{Chernoff}
		\mb{P}\left(\sup_{0\leq t \leq T}\int_{(0,t]\times U}G_{s,u}\,\tilde{N}\left(\mr{d}s,\mr{d}u\right)\geq A \right) \leq \exp\left(-\frac{A}{2}\left(\frac{A}{2C^2\mu(B)T}\wedge\frac{\log(2)}{C}\right)\right).
		\end{align}		
	\end{Lem}
	\prf Let $(M_t)_{0\leq t \leq T}$ be defined by
		\[
		M_t=\int_{(0,t]\times U}G_{s,u}\,\tilde{N}\left(\mr{d}s,\mr{d}u\right).
		\]
		Let $\tau=\inf\left\{t\in[0,T] : M_t \geq A\right\}$ and set $M'_t=M_{t\wedge\tau}$. Note that \eqref{boundG} entails $M'\leq A+C$. For all $\lambda\geq 0$, $e^{\lambda M'}$ is a bounded $(\scr{F}_t)$-local submartingale, hence it is a submartingale and Doob's maximal inequality yields
		\begin{align}
		\mb{P}\left(\sup_{0\leq t \leq T}M_t \geq A\right)=\mb{P}\left(\sup_{0\leq t \leq T}e^{\lambda M'_t} \geq e^{\lambda A}\right)\leq \mb{E}\left(e^{\lambda M'_T}\right)e^{-\lambda A}. \label{Doob}
		\end{align}
		Fix $0\leq \lambda \leq \log(2)/C$. For all $0\leq t \leq T$, we have
		\begin{align}
		e^{\lambda M'_t}&=1+\lambda\int_0^te^{\lambda M'_{s-}}\mr{d}M'_s+\sum_{0<s\leq t}e^{\lambda M'_{s-}}(e^{\lambda \Delta M'_s} -1 -\lambda \Delta M'_{s}) \nonumber\\
		&\leq 1+\lambda\int_0^te^{\lambda M'_{s-}}\mr{d}M'_s +\sum_{0<s\leq t}e^{\lambda M'_{s-}}\lambda^2(\Delta M'_s)^2 \label{ineqexpM'}
		\end{align}
		using \eqref{boundG}. Let $(L_t)_{0\leq t \leq T}$ denote the local-martingale term of the above right handside. Its oblique bracket is given by
		\begin{align*}
		\langle L\rangle_t=\lambda^2\int_{0}^{t\wedge \tau}e^{2\lambda M'_{s-}}G_{s,u}^2\mr{d}s\mu(\mr{d}u) \leq \lambda^2C^2e^{2\lambda A}\mu(B)T,
		\end{align*}
	    hence $(L_t)$ is a square-integrable martingale. Taking the expectation in \eqref{ineqexpM'} yields
	    \begin{align*}
	    e^{\lambda M'_t}&\leq 1+ \lambda^2\mb{E}\left(\int_{(0,t]\times U}\mb{1}_{\left\{s\leq \tau\right\}}e^{\lambda M_{s-}}G_{s,u}^2 N(\mr{d}s,\mr{d}u)\right)
	    \\&\leq 1+ \lambda^2C^2\mu(B)\mb{E}\left(\int_0^t e^{\lambda M'_s}\mr{d}s\right).   
	    \end{align*}
    	Using Fubini's theorem and Grönwall's lemma, we get $e^{\lambda M'_T}\leq e^{\lambda^2C^2\mu(B)T}$, and after plugging this into \eqref{Doob} we obtain
    	\[
    	\mb{P}\left(\sup_{0\leq t \leq T}M_t \geq A\right)\leq e^{\lambda\left(\lambda C^2\mu(B)T-A\right)}.
    	\]
    	The right handside is minimized at
		\[
		\lambda=\frac{A}{2C^2\mu(B)T}\wedge \frac{\log(2)}{C},
		\]
		and \eqref{Chernoff} follows easily. \hfill $\square$

	\bibliographystyle{abbrv}
	\bibliography{biblio_bdpre}
	
\end{document}